\documentclass[11pt]{article}

\usepackage{amsmath,amssymb,amscd,amsthm,latexsym}
\usepackage{hyperref}
\usepackage{color}

\usepackage{ mathrsfs }

\allowdisplaybreaks

\def\scK{{\mathscr K}}

\def\1{{\bf 1}}

\def\nn{\nonumber}

\def\sD {{\cal D}}
 \def\sE {{\cal E}} 
 \def\sF {{\cal F}}
\def\sG {{\cal G}}  
 
 \def\sL {{\cal L}}

 \def\sN {{\cal N}} 
  \def\sR {{\cal R}}

 \def\bN {{\mathbb N}} 
 
 \def\bR {{\mathbb R}}

\def\R {{\mathbb R}}

\newtheorem{thm}{Theorem}[section]
\newtheorem{lemma}[thm]{Lemma}
\newtheorem{defn}[thm]{Definition}

\newtheorem{prop}[thm]{Proposition}
\newtheorem{corollary}[thm]{Corollary}
\newtheorem{remark}[thm]{Remark}
\newtheorem{example}[thm]{Example}
\numberwithin{equation}{section}

\def\qed{{\hfill $\Box$ \bigskip}}

\def\FF{{\mathcal F}}
\def\EE{{\mathcal E}}

\def\R{{\mathbb R}}

\def\E{{\mathbb E}}
\def\P{{\mathbb P}}
\def\N{{\mathbb N}}
\def\eps{\varepsilon}
\def\wh{\widehat}
\def\wt{\widetilde}
\def\pf{\noindent{\bf Proof.} }

\newcommand{\ceil}[1]{\left\lceil #1 \right\rceil}

\DeclareMathOperator*{\esssup}{ess\,sup}

\begin{document}
\allowdisplaybreaks

\title{Heat kernel estimates for symmetric jump processes with mixed polynomial growths
}

\author{ Joohak Bae \quad Jaehoon Kang \quad Panki Kim \quad and\quad
Jaehun Lee}

\date{}
\maketitle

\begin{abstract}
In this paper, we study the transition densities of pure-jump symmetric Markov processes in $\bR^d$, whose jumping kernels are comparable to radially symmetric functions with  mixed polynomial growths. 
Under some mild assumptions on their scale functions, we establish  sharp two-sided estimates  of transition densities (heat kernel estimates) for such processes.
This is the first study on global heat kernel estimates of  jump processes (including non-L\'evy processes) whose weak scaling index is not necessarily strictly less than 2. 
As an application, we proved that the  finite second moment condition on such symmetric Markov process is equivalent to the Khintchine-type  law of iterated logarithm at the infinity.

 \bigskip

\bigskip

\end{abstract}
\allowdisplaybreaks

\medskip

\begin{doublespace}

\tableofcontents

\section{Introduction and main results}

A  principle  notion connecting probability theory with partial differential equation is the heat kernel. 
In probability theory a heat kernel of an operator $\sL$ is the transition density $p(t,x,y)$ (if it exists) of a Markov process $X$, which possesses $\sL$ as its infinitesimal generator. A heat kernel is considered as the fundamental solution for the heat equation $\partial_t u =\sL u$ in the field of partial differential equation. However, except a few special cases, obtaining the explicit expression of  $p(t, x, y)$ is usually impossible. Thus obtaining sharp estimates of $p(t, x, y)$ is a  fundamental issue both in probability theory and partial differential equation.

Although the heat kernels for diffusion processes have been studied for over a century, the heat kernel estimates for the discontinuous Markov process $X$ (equivalently, for the non-local operator $\sL$) have only been studied in recent years.
After pioneering works by \cite{BL02, CK03, Ko00}, obtaining sharp two-sided estimates of heat kernels for various classes of discontinuous Markov processes becomes an active topic in modern probability theory 
(see \cite{BBCK,  BGK, BGR14, Chen,
CHXZ,  CKK08,  CKK,  CK08, CZ16, CZ17, CZ18, HKE, GHH17,GHH18,  GHL14, GS1, GS2, J, KaSz15,  KaSz17, KSV, Kn13,KS12, KR, M, MS1, MS2, Sz11,  Sz17} and references therein.).
In \cite{CK08}, the authors  investigated heat kernel estimates for symmetric discontinuous Markov processes (on a large class of metric measure spaces) whose jumping intensities are comparable to radially symmetric functions of variable order. In particular, heat kernel estimates therein cover a class of  symmetric Markov process $X=(X_t, \P^x, x\in \R^d, t \ge 0)$, without diffusion part, whose jumping kernel $J(x,y)$ satisfies the following conditions:
\begin{equation}\label{e:bm}
 \frac{c^{-1}}{|x-y|^d \phi_1 (|x-y|)} \le J(x, y) \le  \frac{c}{|x-y|^d \phi_1 (|x-y|)}\ , \quad x,y \in \R^d,
\end{equation}
where $\phi_1$ is a non-decreasing function on $[0, \infty)$ satisfying  
\begin{equation}\label{e:bm0}
c_1 (R/r)^{\alpha_1} \le  \phi_1(R)/\phi_1(r) \le c_2 (R/r)^{\alpha_2}, \quad 0<r<R< \infty
\end{equation}  with $\alpha_1, \alpha_2\in (0,2).$
Under the assumptions \eqref{e:bm} and \eqref{e:bm0}, $p(t,x,y)$ has the following estimates: for any $t>0$ and $x,y \in \R^d$,
\begin{equation}\label{e:1.12}
c^{-1}\left(\phi_1^{-1}(t)^{-d}\wedge \frac{t}{|x-y|^d\phi_1 (|x-y|)}\right)
\le p(t,x,y)\le c \left(\phi_1^{-1}(t)^{-d}\wedge
\frac{t}{|x-y|^d\phi_1 (|x-y|)}\right).
\end{equation}
(See \cite[Theorem 1.2]{CK08}.) Here and below, we denote $a \wedge b:=\min\{a, b\}$.
Thus, $\phi_1$ is {\it the scale function}, i.e.,  $|x-y| =\phi_1(t)$ provides 
the borderline for $p(t,x,y)$ to have either 
near-diagonal  estimates or off-diagonal  estimates.
Moreover, it is not difficult to show  from \eqref{e:1.12} that
\begin{equation}\label{e:1.121}
c^{-1}\phi_1(r) \le \E^z[\tau_{B(z,r)}] \le c \phi_1(r) \quad  \text{ for all }\,\, z \in \R^d, \, r>0,
\end{equation}
where $\tau_A$ is the first exit time from $A$ for the process $X$. (See \cite{BGK} and \cite[Section 4.3]{HKE}.) Here, the function $\phi_1$ commonly appears throughout \eqref{e:bm}, \eqref{e:1.12} and \eqref{e:1.121}.
Thus, under the assumptions \eqref{e:bm} and \eqref{e:bm0}, 
\begin{equation}\label{e:EphcomJ}
\frac{c^{-1}}{ J(x,y) r^d} \le \E^z[\tau_{B(z,r)}]\le \frac{c }{ J(x,y) r^d},  \text{ for all } r>0 \text{ and } x,y,z \in \R^d \text{ with } |x-y|=r.
\end{equation}

In this paper, we investigate the estimates of 
transition densities of pure-jump symmetric Markov processes in $\bR^d$, whose jumping kernels satisfying 
\eqref{e:bm} with general mixed polynomial growths, i.e., 
 $\phi_1$  satisfies 
\eqref{e:bm0} with $\alpha_1, \alpha_2 \in (0,\infty)$. 
As a corollary of the main result, 
we have the global sharp two-sided estimates when $\alpha_1$ is greater than $1$ (see Corollary \ref{c:main}).
Unlike heat kernel estimates in \eqref{e:1.12},  $\phi_1$ may not be the scale function for the heat kernel in general (see \eqref{d:Phi} and Theorem \ref{t:main}). For instance, when process $X$ is a subordinate Brownian motion, 
Ante Mimica \cite{M} { established} the heat kernel estimates for the case that the scaling order of characteristic exponent of $X$ may not be strictly
below 2 (see \cite{Sz17} for some partial generalization to L\'evy processes). 
We are strongly motivated by the research done in \cite{M} and consider the case when $\Phi$ in  \eqref{d:Phi}, which is a scale function for the heat kernel, satisfies (local) lower weak scaling condition whose scaling index is greater than 1. Under this assumption, we establish two-sided heat kernel estimates of symmetric jump processes in $\R^d$.
Our result pioneered evidence of sharp heat kernel estimates covering  non-L\'evy processes whose weak scaling index is not necessarily strictly less than 2, which has been a major open problem in this area (c.f., \cite{KaSz15, Sz17}).

In our settings, \eqref{e:EphcomJ} does not hold in general and we only have
\begin{equation}\label{e:EphncomJ}\E^z[\tau_{B(z,r)}]\lnapprox  \frac{c}{J(x,y) r^d},  \text{ for all } r>0 \text{ and } x,y,z \in \R^d \text{ with } |x-y|=r.
\end{equation}
(See (2.1) and Lemma \ref{l:4.3} below.)

In \cite{HKE}, the authors considered heat kernel estimates for mixed-type  symmetric jump processes of on metric measure spaces under a general volume doubling condition. Using 
variants of cut-off Sobolev inequalities and the Faber-Krahn inequalities, they established stability of heat kernel estimates. In particular, they have established heat kernel estimates for $\alpha$-stable-like processes even with 
$\alpha  \ge  2$ when the underlying spaces have walk dimensions larger than $2$ (see \cite{GHH18, HL, MS1, MS2} also).
Note that Euclidean space has the walk dimension $2$;  thus, the results in  \cite{HKE} does not cover our results and, in fact,  a general version of  \eqref{e:EphcomJ} does hold in \cite{HKE}.
 By contrast, some results in  \cite{HKE, PHI} are applicable to our study and we have used several main results in  \cite{HKE, PHI} to show parabolic Harnack inequality and the near-diagonal lower bound of $p(t,x,y)$.

Specifically describing the results of the paper, we start with a description of the setup of this paper.

\begin{defn}
{\rm Let $g:(0,\infty) \to (0,\infty)$, and $a\in(0, \infty]$,  $\beta_1, \beta_2>0$, and $0< c\le1\le C$.
\begin{enumerate}
\item[(1)] For $a \in (0,\infty)$, we say that $g$ satisfies $L_a(\beta_1, c)$ (resp. $L^a(\beta_1, c)$) if
$$ \frac{g(R)}{g(r)} \geq c \left(\frac{R}{r}\right)^{\beta_1} \quad \text{for all} \quad r\leq R< a\;(\text{resp.}\;a\le r\leq R).$$
We also say that $g$ satisfies the weak lower scaling condition near $0$ (resp. near $\infty$) with index $\beta_1$.
\item[(2)] We say that $g$ satisfies $U_a(\beta_2, C)$ (resp. $U^a(\beta_2, C)$) if
$$ \frac{g(R)}{g(r)} \leq C \left(\frac{R}{r}\right)^{\beta_2} \quad \text{for all} \quad r\leq R< a\;(\text{resp.}\;a\le r\leq R).$$
We also say that $g$ satisfies the weak upper scaling condition near $0$ (resp. near $\infty$) with index $\beta_2$.
\item[(3)] 
When $g$ satisfies $U_a(\beta, C)$ (resp. $L_a(\beta, c)$) with $a=\infty$, then we say that $g$ satisfies the global weak upper scaling condition $U(\beta,C)$. (resp. the global weak lower scaling condition $L(\beta,c)$.) 
\end{enumerate}
}\end{defn}

{\it Throughout this paper except Subsections \ref{s:prop_psiPhi} and \ref{s:prop_FFinf}, we will assume that  $\psi:(0,\infty)\to(0,\infty)$ is a 
non-decreasing
 function satisfying $L(\beta_1,C_L)$, $U(\beta_2,C_U)$, and 
\begin{align}\label{e:intcon}
\int_0^{1} \frac{s}{\psi(s)} \,ds<\infty. 
\end{align}
}

Denote $diag=\{(x, x): x\in \R^d\}$ by the diagonal set and let $J:\R^d\times \R^d\setminus{diag}\rightarrow [0,\infty)$ be a symmetric function satisfying 
\begin{equation}\label{e:J_psi}
 \frac{{\bar C}^{-1}}{|x-y|^d\psi(|x-y|)}\leq  J(x,y)\leq \frac{\bar C}{|x-y|^d\psi(|x-y|)},\quad (x,y) \in \R^d \times \R^d \setminus diag
\end{equation}
for some $\bar C\geq 1$.
Note that 
\eqref{e:intcon} combined with \eqref{e:J_psi} and $L(\beta_1,C_L)$ on $\psi$ is a natural assumption to ensure that 
\begin{align}\label{e:levym}
\sup_{x \in \R^d} \int_{\R^d} \big(|x-y|^2 \wedge 1\big) J(x,y) dy 
\leq c \left(\int_0^{1} \frac{sds}{\psi(s)} +\int_1^{\infty} \frac{ds}{s\psi(s)}  \right)<\infty.
\end{align}
For $u, v\in L^2(\bR^d, dx)$, define
\begin{align}\label{e:DF}
\sE(u,v):=\int_{\bR^d\times \bR^d}(u(x)-u(y))(v(x)-v(y)) J(x,y) d x d y
\end{align}
and 
$
\sD(\sE)=\{f\in L^2(\bR^d):\sE(f, f)<\infty\}.$
By applying {the lower scaling assumption} $L(\beta_1,C_L)$ on $\psi$,  \eqref{e:J_psi} and \eqref{e:levym} to \cite[Theorem 2.1]{SU1} and \cite[Theorem 2.4]{SU2},   we observe that $(\sE,\sF)$ is a regular Dirichlet form on $L^2(\bR^d, d x)$. Thus, there is a Hunt process $X$ associated with $(\sE,\sF)$ on $\bR^d$, starting from quasi-everywhere point in $\bR^d$.  Furthermore, by \eqref{e:levym} and \cite[Theorem 3.1]{MU} $X$ is conservative. 

For the first main theorem, we define our scale function as 
\begin{align}\label{d:Phi}
\Phi(r):=\frac{r^2}{2\int_0^{r} \frac{s}{\psi(s)} \,ds}.
\end{align}
In general, the function $\Phi$ is less than or equal to $\psi$ (see \eqref{comp1} below). 
However, these functions may not be comparable unless $\beta_2<2$ (see Lemma \ref{L:comp2} and Section \ref{s:Example} below). 
We remark here that the function $\Phi$ has been observed as the correct scale function recently (see \cite{GS1, GS2, KL, M, SZ}). 
Let $\Phi^{-1}$ be the generalized inverse $\Phi^{-1}(t):= \inf\{s\geq0: \Phi(s)> t   \}$ (with the convention $\inf\emptyset=\infty$).

\begin{thm}\label{t:exp-1}
Let $\psi$ be a non-decreasing function satisfying $L(\beta_1,C_L)$ and $U(\beta_2,C_U)$. Assume that conditions \eqref{e:intcon} and \eqref{e:J_psi} hold. Then, there is a conservative Feller process $X=(X_t, \P^x, x\in \R^d, t \ge 0)$ associated with $(\sE, \sF)$ that starts every point in $\bR^d$. Moreover,  $X$ has a continuous transition density function $p(t,x,y)$ on $(0,\infty)\times\bR^d\times\bR^d$, with the following estimates: there exist $a_U, C, \delta_1>0$ such that
	\begin{equation}\label{e:exp-1}
	p(t,x,y) 
	\leq C\left(\Phi^{-1}(t)^{-d}\wedge \left(\frac{t}{|x-y|^d\psi(|x-y|)} + \Phi^{-1}(t)^{-d} \exp{\left(-\tfrac{a_U |x-y|^2}{\Phi^{-1}(t)^2}\right)}\right) \right)
	\end{equation}
	and 
\begin{equation}\label{e:exp-11}
{p(t,x,y) \ge  C\Phi^{-1}(t)^{-d}{\bf 1}_{\{ |x-y| \le \delta_1\Phi^{-1}(t)\}}+ \frac{Ct}{|x-y|^d \psi(|x-y|)}{\bf 1}_{\{ |x-y| \ge \delta_1\Phi^{-1}(t)\}}}.
\end{equation}
\end{thm}
\noindent The proofs of \eqref{e:exp-1} and \eqref{e:exp-11} are given in Section \ref{s:ODUE} and Proposition \ref{l:LHK-psi}, respectively.

\begin{remark}
\rm
Theorem \ref{t:exp-1}, in particular, implies that if $ \psi(r) \asymp \Phi (r)$ for all large $r>1$,
we have sharp two-sided estimates for $t>1$;  there exists $C>0$ such that for all $t>1$ and $x, y \in \R^d$,
\begin{equation}\label{e:exp-111}
C^{-1} \left(\Phi^{-1}(t)^{-d} \wedge \frac{t}{|x-y|^d \Phi(|x-y|)}\right) \le p(t,x,y) \le  C\left(\Phi^{-1}(t)^{-d} \wedge \frac{t}{|x-y|^d \Phi(|x-y|)}\right).\end{equation}
The condition $ \psi(r) \asymp \Phi (r)$ for $r>1$ is equivalent to that $\Phi$ satisfies $U^1(\delta_0,\wt C_U)$ with $\delta_0<2$. (See Lemma \ref{L:comp2}.) Such estimates for large time in \eqref{e:exp-111} under the weak scaling condition at the infinity are even unknown before.
\end{remark}

 Let $$G(x,y)=\int_0^\infty p(t,x,y)dt$$ be the Green function for $X$.  As a corollary of Theorem \ref{t:exp-1}, we observe  the two-sided sharp estimate for Green function. Recall that $\psi$ satisfies $L(\beta_1, C_L)$ and $U(\beta_2,C_U)$. 
\begin{corollary}\label{c:Green}
Suppose that  $d >\beta_2 \wedge 2$. Then there exists $c\ge1$ such that for any $x,y \in \R^d$,
\begin{align*}
c^{-1} \Phi(|x-y|)|x-y|^{-d} \le G(x,y)\le c \,\Phi(|x-y|)|x-y|^{-d}.
\end{align*}
\end{corollary}
\noindent Using our scale function $\Phi$, we define 
\begin{align}
\scK(s):=\sup_{b\le s}\frac{\Phi(b)}{b}.\label{d:F}
\end{align} 
If $\Phi$ satisfies $L_a(\delta,\wt C_L)$ with $\delta>1$, then $\scK(0)=0$ and $\scK$ is non-decreasing. 
Thus, the generalized inverse $\scK^{-1}(t):= \inf\{s\geq0: \scK(s)> t   \}$ is well defined on $[0,\sup_{b<\infty}\frac{\Phi(b)}{b})$.
Now we define
\begin{align*}
\sG(c, t, r)&:= \frac{t}{r^d\psi(r)} +  \Phi^{-1}(t)^{-d} \exp\left(-c\frac{r}{\scK^{-1}(t/r)} \right) .
\end{align*} 
If $\Phi$ satisfies $L^a(\delta, \wt C_L)$ with $\delta>1$,
\begin{align}
\scK_\infty(s)=\scK_{\infty,a}(s):=\begin{cases}
\displaystyle \sup_{a\le b\le s}\frac{\Phi(b)}{b} &\text{if }s\ge a,\\
a^{-2}{\Phi(a)}s&\text{if } 0<s< a.
\end{cases}\label{d:Finf}
\end{align}
Then, $\scK_\infty$ and the generalized inverse $\scK_\infty^{-1}$ are well-defined and non-decreasing on $[0, \infty)$. Similarly, we define 
\begin{align*}
\sG_{\infty}(c, t, r)&:= \frac{t}{r^d\psi(r)} +  \Phi^{-1}(t)^{-d} \exp\left(-c\frac{r}{\scK_{\infty}^{-1}(t/r)} \right).
\end{align*} 
Some properties of $\scK$ and $\scK_\infty$ are shown in Subsection \ref{s:prop_FFinf}. Here is the main result of this paper.

\begin{thm}\label{t:main}
Let $\psi$ be a 
non-decreasing function satisfying $L(\beta_1,C_L)$ and $U(\beta_2,C_U)$. Assume that conditions \eqref{e:intcon} and \eqref{e:J_psi} hold, and $\Phi$ satisfies $L_a(\delta,\wt C_L)$ or $L^a(\delta,\wt C_L)$ for some $a >0$ and $\delta>1$. 
Then, the following estimates hold:
\begin{enumerate}
\item[(1)] When $\Phi$ satisfies $L_a(\delta,\wt C_L)$: \\
For every $T>0$, there exist positive constants $c_1=c_1(T, a, \delta, \beta_1, \beta_2, \wt C_L, C_L, C_U ) \ge1$ and $a_U \le a_L$ such that for any $(t,x,y)\in(0,T)\times\bR^d\times\bR^d$,
\begin{align}\label{e:HKE0}
c_1^{-1} \left(\Phi^{-1}(t)^{-d}\wedge \sG(a_L, t, |x-y|)\right) 
\le p(t,x,y) \le  c_1\left(\Phi^{-1}(t)^{-d}\wedge \sG(a_U, t, |x-y|)\right).
 \end{align}
 Moreover, if $\Phi$ satisfies $L(\delta,\tilde{C}_L)$, then \eqref{e:HKE0} holds for all $t \in (0,\infty)$.
 
 \item[(2)] When $\Phi$ satisfies $L^a(\delta,\wt C_L)$: \\
For every $T>0$, there exist positive constants $c_2=c_2(T, a, \delta, \beta_1, \beta_2, \wt C_L, C_L, C_U )\ge1$ and $a'_{U} \le a'_{L}$ such that for any $(t,x,y)\in[T,\infty)\times\bR^d\times\bR^d$,
\begin{align*}
c_2^{-1} \left(\Phi^{-1}(t)^{-d}\wedge \sG_{\infty}(a'_{L}, t, |x-y|)\right) 
\le p(t,x,y) \le  c_2\left(\Phi^{-1}(t)^{-d}\wedge \sG_{\infty}(a'_{U}, t, |x-y|)\right).
 \end{align*}
In particular, if $\delta=2$, then  $c^{-1}t\le  \scK_{\infty}^{-1}(t)\le c t$ for $t \ge T$.
\end{enumerate}
\end{thm}
Theorem \ref{t:main}(2) covers \cite[Corollary 3.11]{SW} where $\psi(r)=r^{2+\eps}$, $r>1$ and $\eps>0$, is considered.
Using Theorems \ref{t:exp-1} and \ref{t:main}(2), we will show in Section \ref{s:App} that the finite second moment condition is equivalent to the Khintchine-type law of iterated logarithm at the infinity.
 In  \cite{G}, Gnedenko proved  this result for the L\'evy process  (see also \cite[Proposition 48.9]{Sat}).
The equivalence between the law of iterated logarithm and the  finite second moment condition for the non-L\'evy process  has been a long standing open problem since  the work done in \cite{G}.

\begin{remark}
\rm
\begin{enumerate}
\item[(1)] The assumption that $\Phi$ satisfies the local weak lower scaling condition with index $\delta>1$ is used only in proving off-diagonal estimates of the transition density function.
\item[(2)] Although we use 
$\scK_\infty$ in Theorem \ref{t:main}(2) instead of $\scK$, neither the value $a \in (0, \infty)$ nor the behavior of $\scK_\infty$ near zero is irrelevant. 
See Remarks \ref{mwsc} and \ref{r:replace}.
\end{enumerate}
\end{remark}

\noindent If $\psi$ satisfies $L_a(\delta,\wt C_L)$, then $\delta<2$ and $\Phi$ satisfies $L_a(\delta,\wt C_L)$ by Lemma \ref{l:rel1}. Thus, we have the following: 
\begin{corollary}\label{c:main}
Let $\psi$ be a 
non-decreasing function that satisfies $L(\beta_1,C_L)$ and $U(\beta_2,C_U)$. Assume that conditions \eqref{e:intcon} and \eqref{e:J_psi} hold and $\psi$ satisfies $L_a(\delta,\wt C_L)$ for some $a>0$ and $\delta>1$. Then, for any $T>0$, there exist positive constants $c=c( \beta_1, \beta_2, \delta,  C_L, C_U, \wt C_L, T)\ge1$ and $a_{U} \le a_{L}$ such that \eqref{e:HKE0} holds for all $(t,x,y)\in(0,T)\times\bR^d\times\bR^d$. Moreover, if $\beta_1>1$, then \eqref{e:HKE0} holds for all $t \in (0,\infty)$.
\end{corollary}

A non-negative $C^\infty$ function $\phi$ on $(0,\infty)$ is called a Bernstein function if
$(-1)^n \phi^{(n)} (\lambda) \le 0$
for every $n\in \mathbb N$ and $\lambda>0$.
The exponent $({r}/{\Phi^{-1}(t)})^2$ in \eqref{e:exp-1} is not comparable to ${r}/{\scK^{-1}(t/r)}$ in general (see Lemma \ref{l:F} and Corollary \ref{c:example} below).
However, the following corollary indicates that we can replace ${r}/{\scK^{-1}(t/r)}$ with a simpler function $({r}/{\Phi^{-1}(t)})^2$  if we additionally assume that $r \mapsto \Phi(r^{-1/2})^{-1}$ is a Bernstein function.

\begin{corollary}\label{c:main2}
Let $\psi$ be a non-decreasing function satisfying $L(\beta_1,C_L)$ and $U(\beta_2,C_U)$. Assume that conditions \eqref{e:intcon} and \eqref{e:J_psi} hold, $\Phi$ satisfies $L_a(\delta,\wt C_L)$ some $a>0$ and $\delta>1$, and $r \mapsto \Phi(r^{-1/2})^{-1}$ is a Bernstein function. Then, for any $T>0$, there exist positive constants $c\ge1$ and $a_{U} \le a_{L}$ such that for all $(t,x,y)\in(0,T)\times\bR^d\times\bR^d$,
\begin{align}
\begin{split}\label{sbm}
&c^{-1} \left(\Phi^{-1}(t)^{-d}\wedge \left(\frac{t}{|x-y|^d\psi(|x-y|)} +  \Phi^{-1}(t)^{-d} \exp\left(-a_L\tfrac{|x-y|^2}{\Phi^{-1}(t)^2} \right) \right)\right) \\
&\le p(t,x,y) \le  c\left(\Phi^{-1}(t)^{-d}\wedge \left(\frac{t}{|x-y|^d\psi(|x-y|)} +  \Phi^{-1}(t)^{-d} \exp\left(- a_U\tfrac{|x-y|^2}{\Phi^{-1}(t)^2} \right) \right)\right).
\end{split}
 \end{align}
Moreover, if $\Phi$ satisfies $L(\delta, C_L)$ with $\delta>1$, \eqref{sbm} holds for all $t \in (0,\infty)$.
\end{corollary}

For given function $f:(0,\infty) \to (0,\infty)$, we say that $f$ varies regularly at $0$ (resp. at $\infty$) with index $\delta_0 \in [0,\infty)$ if 
$$ \lim_{ x \to 0} \frac{f(\lambda x)}{f(x)} = \lambda^{\delta_0} \quad (\text{resp.}\lim_{ x \to \infty} \frac{f(\lambda x)}{f(x)} = \lambda^{\delta_0}) $$
for any $\lambda>0$. Especially, we say that $f$ is slowly varying at $0$ (resp. at $\infty$) if $f$ varies regularly at $0$ (resp. at $\infty$) with index $0$.
Note that if  $f$ is non-increasing and is regularly varying at $0$ (resp. at $\infty$) with index $\delta_0$, then for any $a>0$ and $0<\underline{\delta} < \delta_0 < \bar{\delta}$, there is $C_U,C_L>0$ such that $f$ satisfies both $U_a(\bar{\delta}, C_U)$ and $L_a(\underline{\delta},C_L)$ (resp. $U^a(\bar{\delta}, C_U)$ and $L^a(\underline{\delta},C_L)$).

Recall that $Y=(Y_t)_{t\geq 0}$ is  a pure-jump isotropic unimodal L\' evy process in $\R^d$, $d\geq 1$ if its characteristic exponent $\xi \mapsto \Psi(|\xi|)$ is 
\begin{align*}
    \Psi(|\xi|)=\int_{\R^d}(1-e^{i\xi\cdot x})\,j(|x|)dx,
\end{align*}
where $j(|x|)$ is the L\'evy  kernel of $Y$ and $r \mapsto j(r)$ is non-increasing. Clearly, if $r \mapsto r^{-d}/j(r)$ is comparable to a non-decreasing function  satisfying $L(\beta_1,C_L)$ and $U(\beta_2,C_U)$, then Theorem \ref{t:main} can be applied to get the estimates of the transition density of $Y$. 
We also can check $\Psi$  directly and, if $\Psi(\lambda)-\frac{\lambda}{2}\Psi'(\lambda)$ varies regularly, we can 
show that the estimates in Corollary \ref{c:main2} holds true.

Given a characteristic exponent $\xi \mapsto \Psi(|\xi|)$, let 
$ \phi_\Psi(\lambda)=\int_{\R^d}(4\pi \lambda)^{-d/2}\exp({-\frac{|\xi|^2}{4\lambda}})\Psi(|\xi|)dx,$ which is 
a Bernstein function (see \cite[Remark 3.2]{KM}).
\begin{corollary}\label{c:levy}
Suppose $Y=(Y_t)_{t\geq 0}$ is  a pure-jump isotropic unimodal L\' evy process in $\R^d$, $d\geq 1$ with the characteristic exponent $\xi \mapsto \Psi(|\xi|)$
 and that 
$q(t,|x-y|)$ is the transition density of $Y$.
Suppose that  $g(\lambda):=\Psi(\lambda)-\frac{\lambda}{2}\Psi'(\lambda)$ varies regularly at $0$ with index $\alpha_1 \in 
(0, 4)$  and varies regularly at $\infty$ with index $\alpha_2 \in 
(1, 2]$.  Then, there exist $c_i \in (0, \infty)$ such that with $H(\lambda):=\phi_\Psi(\lambda)-\lambda\phi_\Psi'(\lambda) $,
$$
\lim_{\lambda \to 0} \frac{ \phi_\Psi(\lambda) }{\Psi(\sqrt{\lambda})}= c_1(\alpha_1, d), \quad 
\lim_{\lambda \to  \infty} \frac{ \phi_\Psi(\lambda) }{\Psi(\sqrt{\lambda})}= c_2(\alpha_2, d),    $$
 $$
 \lim_{\lambda \to 0} \frac{ H(\lambda) }{g(\sqrt{\lambda})}= c_3(\alpha_1, d), \quad
\lim_{\lambda \to  \infty} \frac{ H(\lambda) }{g(\sqrt{\lambda})}= c_4(\alpha_2, d)  $$
and, 
for every $T>0$, there exist positive constants $c_1 \ge1$ and $a_U \le a_L$ such that for any $(t,x)\in(0,T)\times\bR^d$,
\begin{align}
\begin{split}\label{e:HKE0L}
&c_1^{-1} \left(\phi_\Psi^{-1}(t^{-1})^{d/2}\wedge \left(\frac{tH(|x|^{-2})}{|x|^d} +  \phi_\Psi^{-1}(t^{-1})^{d/2} \exp\left(-a_L{|x|^2}{\phi_\Psi^{-1}(t^{-1})} \right) \right)\right) \\
&\le q(t,|x|) \le  c_1\left(\phi_\Psi^{-1}(t^{-1})^{d/2}\wedge \left(\frac{tH(|x|^{-2})}{|x|^d} +  \phi_\Psi^{-1}(t^{-1})^{d/2}\exp\left(-a_U{|x|^2}{\phi_\Psi^{-1}(t^{-1})} \right) \right)\right).
\end{split}
 \end{align}
 Moreover, if $\alpha_1 >1$, then \eqref{e:HKE0L} holds for all $t \in (0,\infty)$.
\end{corollary}
Corollary \ref{c:levy} is a direct consequence of Corollary \ref{c:main2} and \cite[(1.7), Propositions 3.3(1), 3.4 and 3.7]{KM}.

The remainder of the paper is organized as follows. Section \ref{s:pre} describes the properties of $\psi$ and $\Phi$, and verifies some relationships between them. Moreover, the properties of $\scK$ and $\scK_{\infty}$ under the lower weak scaling assumption on $\Phi$ are verified.
 Section \ref{s:NDE} proves a preliminary  upper bound and near-diagonal estimates of transition density function. Subsection \ref{s:PI} presents the  Poincar\'e inequality, which is the first step to find a correct scale function. Using this Poincar\'e inequality, in Subsection \ref{s:Nash}, we show that Nash inequality holds, which  yields the existence of the transition density function and its near-diagonal upper bound. Subsection \ref{s:UHK1} uses  scaled processes of $X$ to obtain  an upper bound of the transition density function (see Theorem \ref{t:UHK-Phi}). Although this upper bound is not sharp, it is  adequate  to get the lower bound of survival probability and $\text{CSJ}(\Phi)$ condition defined in \cite{HKE}. The Poincar\'e inequality and the lower bound of survival probability provide the upper and lower bounds of the mean exit times of balls. Using the  mean exit time estimates of  balls, Subsection \ref{s:PHI} shows the near-diagonal lower bound of the transition function, parabolic Harnack inequality, and parabolic H\"older regularity by applying the results in \cite{HKE, PHI}.  
 Section \ref{s:ODE} describes the proof of off-diagonal estimates of the transition density function. Subsection \ref{s:ODUE} and Subsection \ref{s:ODLE} prove the off-diagonal upper bound and lower bound of the transition density function, respectively. 
As an application of the main result, in Section \ref{s:App} we show that the finite second moment condition is equivalent to the Khintchine-type law of iterated logarithm at the infinity.
 Section 6 provides examples covered by the main result.

\vspace{3mm}

\noindent
{\it Notations} : Throughout this paper, the constants $C_1$, $C_2$, $C_3$, $C_4$, $C_5$, $C_6$, $C_L$, $C_U$, $\wt C_L$, $\beta_1$, $\beta_2$, $\delta$, $\delta_1$ will remain the same, whereas  $C$, $c$, and  $c_0$, $a_0$, $c_1$, $a_1$, $c_2$, $a_2$, $\ldots$  represent 
constants having insignificant values that may be changed  from one
appearance to another. All these constants are positive finite.
The labeling of the constants  $c_0$, $c_1$, $c_2$, $\ldots$ begins anew in the proof of
each result.  $c_i=c_i(a,b,c,\ldots)$, $i=0,1,2,  \dots$, denote generic constants depending on $a, b, c, \ldots$. The dependence on the dimension $d \ge 1$ and the constant $\bar C$ in \eqref{e:J_psi}
may not
be explicitly mentioned.
Recall that we use the notation $a \wedge b = \min \{ a, b\}$. We also denote
 $a \vee b := \max \{ a, b\}$, $\bR_+:=\{r\in \bR:r>0\}$, and $B(x,r):=\{y\in\bR^d: |x-y|<r\}$. We use the notation $f\asymp g$ if the quotient $f/g$ remains bounded between two positive constants.

\section{Preliminary}\label{s:pre}

\subsection{Basic properties of $\psi$ and $\Phi$}\label{s:prop_psiPhi}
In this subsection, we will observe  some elementary properties of $\psi$ and $\Phi$. Since $\psi$ is non-decreasing, we see that
\begin{align}\label{comp1}
\Phi(r)
=\frac{r^2}
{2\int_0^{r} \frac{s}{\psi(s)} \,ds} \le \frac{r^2}
{2\int_0^{r} \frac{s}{\psi(r)} \,ds}=\psi(r).
\end{align}
Thus, under \eqref{e:J_psi}, we obtain that for any $x,y\in \bR^d$,
\begin{align}\label{e:JPhile}
J(x,y)\le \frac{c}{|x-y|^d\Phi(|x-y|)}.
\end{align}
Since
$
(1/2\Phi(r))'={2}r^{-3} \int_0^{r}s( \frac{1}{ \psi(r)} -\frac{1}{\psi(s)} ) \,ds  \le 0,$
$\Phi(s)$ is also non-decreasing. 
Note that, since ${r^2}/{\Phi(r)} = 2\int_0^r \frac{s}{\psi(s)}ds
$
is increasing in $r$, we have that for any $0<r \le R$,
\begin{equation}\label{e:Phi}
\frac{\Phi(R)}{\Phi(r)} \le \frac{R^2}{r^2}.
\end{equation}
From this, we see that if $\Phi$ satisfies 
$L^a(\beta,  c)$, then $\beta\le2$.

\begin{remark}\label{mwsc}
\rm Suppose $g:(0,\infty)\to(0,\infty)$ is non-decreasing.
If $g$ satisfies $L_a(\beta,  c)$, then $g$ satisfies $L_b(\beta,  c(ab^{-1})^{\beta})$ for any $b > a$. Indeed, for $r\leq a \leq R\leq b$, 
$$
g(R)\geq g(a)\geq c\left(\frac{a}{r}\right)^{\beta}g(r)\geq c\left(\frac{a}{b}\right)^{\beta}\left(\frac{R}{r}\right)^{\beta}g(r)
$$
and for $a \leq r\leq R\leq b$,
$$
g(R)\geq g(r)\geq c\left(\frac{a}{b}\right)^{\beta}\left(\frac{R}{r}\right)^{\beta}g(r).
$$
Similarly, if $g$ satisfies $L^a(\beta, c)$, then $g$ satisfies $L^b(\beta,  c(a^{-1}b)^{\beta})$ for $b<a$.
\end{remark}

The next result is straightforward. We skip the proof. 
\begin{lemma}\label{lem:inverse}
Let $g:(0,\infty) \rightarrow (0,\infty)$ be a 
non-decreasing function with $g(\infty)=\infty.$
\begin{enumerate}
\item[(1)] If $g$ satisfies $L_a(\beta, c)$ (resp. $U_a(\beta, C)$),  then $g^{-1}$ satisfies $U_{g(a)}(1/\beta,c^{-1/\beta})$\\  (resp. $L_{g(a)}(1/\beta,C^{-1/\beta})$). 
\item[(2)] If $g$ satisfies $L^a(\beta, c)$  (resp. $U^a(\beta, C)$),  then $g^{-1}$ satisfies $U^{g(a)}(1/\beta,c^{-1/\beta})$\\ (resp. $L^{g(a)}(1/\beta,C^{-1/\beta})$).
\end{enumerate}
\end{lemma}
The following lemma will be used in the paper several times.
\begin{lemma}\label{l:int_outball}
Assume that $\psi$ satisfies $L(\beta, c)$ and  $U(\wh \beta, C)$. Then, 
for any $x\in \R^d$ and $r>0$,
\begin{align}\label{int_outball}
\int_{B(x,r)^c}\frac{1}{|x-y|^d\psi(|x-y|)}dy= c_0(d) \int_{r}^{\infty} \frac1{s\psi(s)}  ds \asymp \frac{1}{\psi(r)}.
\end{align}
\end{lemma}
\pf Using $L(\beta, c)$, 
 \begin{align*}
 \int_{r}^{\infty} \frac1{s\psi(s)}  ds=  \frac{1}{\psi(r)} \int_{ r}^{\infty} \frac{\psi(r)}{s\psi(s)}  ds \le   \frac{c^{-1}}{\psi(r)} \int_{r}^{\infty} \frac{r^{\beta}}{s^{1+\beta}}  ds 
\le  \frac{c_1}{\psi(r)}. 
    \end{align*}
  On the other hand, by $U(\wh \beta, C)$ we have
 \begin{align*}
 \int_{r}^{\infty} \frac1{s\psi(s)}  ds=  \frac{1}{\psi(r)} \int_{ r}^{\infty} \frac{\psi(r)}{s\psi(s)}  ds \ge   \frac{C^{-1}}{\psi(r)} \int_{r}^{\infty} \frac{r^{\wh \beta}}{s^{1+\wh \beta}}  ds 
\ge  \frac{c_2}{\psi(r)}.  
    \end{align*}
\qed

The next lemma shows that the index in the weak scaling condition for $\Phi$ is always in $(0, 2]$. 

\begin{lemma}\label{l:rel1}
Let $a\in(0,\infty]$, $0<\beta\le \wh \beta$, $0<c\le1\le C$.
\begin{enumerate}
\item[(1)] If $\psi$ satisfies $U_a(\wh \beta, C)$, then $\Phi$ satisfies $U_a(\wh \beta \wedge 2, C)$. 
\item[(2)] If $\psi$ satisfies \eqref{e:intcon} and $L_a(\beta, c)$, then $\beta<2$ and $\Phi$ satisfies $L_a(\beta, c)$.
\end{enumerate}
\end{lemma}
\pf
(1) Note that $
\Phi(\lambda r)={r^2}/
{\int_0^{r} \frac{t}{\psi(\lambda t )} \,dt}.
$ Using $U_a(\wh \beta,C)$, we have $\psi(\lambda t ) \le C \psi( t )\lambda^{\wh \beta}$ for any $0< t \le \lambda t \le a$ . Thus, for any $0<r \le \lambda r \le a$ we have
$$
\Phi(\lambda r)=\frac{r^2}
{\int_0^{r} \frac{t}{\psi(\lambda t )} \,dt} \le \frac{r^2}
{\int_0^{r} \frac{t}{ C \psi( t )\lambda^{\wh \beta}} \,dt}=C \lambda^{\wh \beta}\frac{r^2}
{\int_0^{r} \frac{t}{  \psi( t )} \,dt}=C \lambda^{\wh \beta} \Phi( r),
$$
which concludes $U_a(\wh \beta,C)$ of $\Phi$. Combining with \eqref{e:Phi} we conclude $U_a(\wh \beta \land 2, C)$ for $\Phi$.

(2) Similarly, using $L_a(\beta,c)$ we have
$$
\Phi(\lambda r)=\frac{r^2}
{\int_0^{r} \frac{t}{\psi(\lambda t )} \,dt} \ge \frac{r^2}
{\int_0^{r} \frac{t}{ c \psi( t )\lambda^{\beta}} \,dt}=c \lambda^{\beta}\frac{r^2}
{\int_0^{r} \frac{t}{  \psi( t )} \,dt}=c \lambda^{\beta} \Phi( r)
$$
for $0<r \le \lambda r \le a$, and this is equivalent to  $L_a(\beta,c)$ for $\Phi$. 

Now assume that $\psi$ satisfies $L_a(\beta,c)$ with $\beta \ge 2$. Then
\begin{align*}
\int_0^a \frac{s}{\psi(s)}ds = \int_0^a \frac{s}{\psi(a)} \frac{\psi(a)}{\psi(s)}ds \ge \frac{c}{\psi(a)} \int_0^a \frac{1}{s}ds = \infty,
\end{align*}
which conflicts with \eqref{e:intcon}. This finishes the proof. \qed

The comparability of $\psi$ and $\Phi$ is equivalent to that the index of the weak upper scaling condition is strictly less than 2.

\begin{lemma}\label{L:comp2}
Let  $a \in (0,\infty)$.

(1) There exists $c \in (0,1]$ such that $c \psi(r) \le \Phi(r)$ for all $r<a$ (resp. $r<\infty$), if and only if there exist $\beta \in (0,2)$ and $C\ge1$ such that $\psi$ satisfies $U_{a}(\beta,C)$ (resp. $U(\beta,C)$).

(2) There exists $c \in (0,1]$ such that $c \psi(r) \le \Phi(r)$ for all $r \ge a$, if and only if there exist $\beta \in (0,2)$ and $C \ge 1$ such that $\psi$ satisfies $U^{a}(\beta,C)$.

\end{lemma}
\pf Let  $\Psi(r):=r^2/\psi(r)$.
By the zero version of 
\cite[Corollary 2.6.4]{BGT} (\cite[Corollary 2.6.2]{BGT}, respectively)
the condition $\Psi(r) \asymp \int_0^r \Psi(s)s^{-1}ds$ for all $r<a$ ($r \ge a$ respectively)
is equivalent to that $\Psi$ satisfies the
weak lower scaling condition near zero (near the infinity, respectively) with index $\delta_0>0$, which is also equivalent to that  both $\psi$ and $\Phi$ satisfy the weak upper scaling condition near $0$ (near $\infty$, respectively) with index $2-\delta_0<2$.
\qed
\begin{prop}\label{p:regva}
Suppose that
$$
\lim_{\lambda \to 0} \frac{\Phi(\lambda)}{ \psi(\lambda)} =0. \quad
(\lim_{\lambda \to \infty} \frac{\Phi(\lambda)}{ \psi(\lambda)} =0, \text{respectively}.) 
$$
Then $\Phi$ varies regularly at $0$ (at $\infty$, respectively) with index 
 $2$. 
\end{prop}
\pf
Without loss of generality, we can assume that $\Phi(1)=1$. 
Since $$
\frac{d}{d \lambda}\log (\lambda^2/\Phi(\lambda))=
2\frac{\Phi(\lambda)}{\lambda^2}\frac{d}{d \lambda} \int_{0}^\lambda \frac{sds}{\psi(s)}=2
\frac{\Phi(\lambda)}{\lambda \psi(\lambda)},
$$
we have 
$$
\frac{\lambda^2}{\Phi(\lambda)}=\exp \left( \log (\lambda^2/\Phi(\lambda)) \right)=
\exp \left( 2 \int_1^\lambda \frac{\Phi(s)}{s \psi(s)}ds \right).
$$
Thus, 
using   \cite[Theorem 1.3.1 and (1.5.2)]{BGT} and the zero version of it, we conclude from the above display that $\frac{\lambda^2}{\Phi(\lambda)}$
varies slowly at $0$ if $\lim_{\lambda \to 0} \frac{\Phi(\lambda)}{ \psi(\lambda)} =0$ (at $\infty$ if $\lim_{\lambda \to \infty} \frac{\Phi(\lambda)}{ \psi(\lambda)} =0$, respectively).  
\qed
\subsection{Basic properties of $\scK$ and $\scK_{\infty}$}\label{s:prop_FFinf}

In this subsection, under the assumption that $\Phi$ satisfies $L_a(\delta,\widetilde{C}_L)$ or $L^a(\delta,\widetilde{C}_L)$ with $\delta>1$, we establish some basic properties of $\scK$ and $\scK_{\infty}$ defined in \eqref{d:F} and \eqref{d:Finf}. We remark here that by Proposition \ref{p:regva}, such assumptions hold true if 
${\Phi(\lambda)}/{ \psi(\lambda)} \to 0$ as $\lambda \to 0$ (resp. as $\lambda \to \infty$).

\begin{lemma}\label{l:Fprop}
If $\Phi$ satisfies $L_a(\delta,\wt C_L)$ with $\delta>1$ and $a \in (0,\infty]$, then
\begin{equation}
\label{e:F1}\frac{\Phi(t)}{t} \le \scK(t) \le \wt C_L^{-1}\frac{\Phi(t)}{t},\quad\text{for}\;\;t< a,
\end{equation}
and  
\begin{align}\label{e:wscF}
\wt C_L^2\left(\frac{t}{s}\right)^{\delta-1}\leq \frac{\scK(t)}{\scK(s)}\leq \wt C_L^{-1}\frac{t}{s},\quad\text{for}\;\;s\leq t< a.
\end{align}
\end{lemma}
\pf The first inequality in \eqref{e:F1} immediately follows from the definition of $\scK$. Since $\Phi$ satisfies $L_a(\delta,\wt C_L)$ with $\delta>1$,  we have that for any $b\le t<a$,
$
\wt C_L(t/b)<\wt C_L(t/b)^{\delta}\le {\Phi(t)}/{\Phi(b)},
$
which implies the second inequality in \eqref{e:F1}.

The inequality  \eqref{e:wscF} follows from   \eqref{e:Phi}, \eqref{e:F1} and the  $L_a(\delta,\wt C_L)$ for $\Phi$ as 
\begin{align*}
\wt C_L^{-1}\frac{t}{s} \ge \wt C_L^{-1}\frac{\Phi(t)}{t}\frac{s}{\Phi(s)} \ge  \frac{\scK(t)}{\scK(s)}\ge \wt C_L \frac{\Phi(t)}{t}\frac{ s}{\Phi(s)}\ge\wt C_L^2\left(\frac{t}{s}\right)^{\delta-1}. 
\end{align*}
\qed \\
Using Remark \ref{mwsc} and \eqref{e:F1}, we see that under $L_a(\delta,\wt{C}_L)$ for $\Phi$, we have that for any $0<t<b$,
\begin{align*}
\scK(t)\leq \wt C_L^{-1}(b/a)^\delta \frac{\Phi(t)}{t}.
\end{align*} 
For notational  convenience, we introduce  an auxiliary function 
$\wt\Phi_a(s):=\frac{\Phi(a)}{a^2}s^2\1_{\{0<s< a\}}+\Phi(s)\1_{\{s\ge a\}}$ so that 
\begin{align*}
\scK_\infty(s)=\scK_{\infty,a}(s)=\sup_{b\le s}\frac{\wt\Phi_a(b)}{b}.
\end{align*}
The following lemma shows the relation between $\Phi$ and $\wt \Phi_a$.

\begin{lemma}\label{PhiwtPhi}
\begin{enumerate}
\item[(1)] For  any $t>0$, $\wt\Phi_a(t)\le \Phi(t)$ and for $t\ge c>0$, $\wt \Phi_a(t)\ge ((c/a)^2\wedge 1) \Phi(t)$.  
\item[(2)] For $0<s<t$, 
${\wt \Phi_a(t)}/{\wt \Phi_a(s)}\le {t^2}/{s^2}.
$
\item[(3)] Suppose $\Phi$ satisfies $L^a(\delta, \wt C_L)$ with some $\delta \le 2$. Then, $\wt\Phi_a$ satisfies $L(\delta, \wt C_L)$. 
\end{enumerate}
\end{lemma}
\pf 
(1) Since $\wt \Phi_a(t)=\Phi(t)$ for $t\ge a$, it suffices to prove the case $t<a$. By \eqref{e:Phi},  $\Phi(a)\le  (a/t)^2\Phi(t)$, which implies that
$\wt \Phi_a(t)=\frac{\Phi(a)}{a^2}t^2\le \Phi(t)$.
Since $\Phi$ and $\wt \Phi_a$ are non-decreasing, we see that for $c\le t< a$,
$\wt\Phi_a(t)\ge \wt\Phi_a(c)\ge(c/a)^2\Phi(t).$ 

(2) By \eqref{e:Phi} and definition of $\wt{\Phi}_a$, we only need to verify the case $0<s<a\le t$. Using \eqref{e:Phi} again, for any  $0<s<a\le t$, 
$$\frac{\wt \Phi_a(t)}{\wt \Phi_a(s)}=\frac{\Phi(t)}{\Phi(a)s^2/a^2}\le\frac{a^2}{s^2}\left(\frac{t}{a}\right)^{2}=\left(\frac{t}{s}\right)^{2}.$$

(3) Clearly, for $0<s<t< a$, we have 
${\wt \Phi_a(t)}/{\wt \Phi_a(s)}=({t}/{s})^{2}\ge({t}/{s})^{\delta}.$ 
For $0<s<a\le t$, 
$$\frac{\wt \Phi_a(t)}{\wt \Phi_a(s)}=\frac{\Phi(t)}{\Phi(a)s^2/a^2}\ge\frac{a^2}{s^2}\wt C_L\left(\frac{t}{a}\right)^{\delta}\ge \wt C_L\left(\frac{a}{s}\right)^{\delta}\left(\frac{t}{a}\right)^{\delta}=\wt C_L\left(\frac{t}{s}\right)^{\delta}.$$
\qed

By Lemma \ref{PhiwtPhi}(1) and \eqref{comp1}, $\wt \Phi_a(t) \le \psi(t)$ for all $t>0$ and $a>0$.
In the following lemma, we will see some properties of $\scK_{\infty, a}$ which is similar to Lemma \ref{l:Fprop}.

\begin{lemma}\label{l:Finfprop}
Let $a \in (0,\infty)$. If $\Phi$ satisfies $L^a(\delta, \wt C_L)$ with $\delta>1$, then 
\begin{align}\label{e:F_11}
\frac{\wt\Phi_a(t)}{t}\leq \scK_{\infty, a}(t)\leq \wt C_L^{-1}\frac{\wt\Phi_a(t)}{t},\quad\text{for}\;\;t>0,
\end{align}
and 
\begin{align}\label{e:wscF_1}
\wt C_L^2\left(\frac{t}{s}\right)^{\delta-1}\leq \frac{\scK_{\infty, a}(t)}{\scK_{\infty, a}(s)}\leq \wt C_L^{-1}\frac{t}{s},\quad\text{for}\;\;t>s>0.
\end{align}
Moreover, for any $c_1>0$, there exists $c_2=c_2(c_1, a, \delta, \wt C_L)\ge1$ such that for any $t\ge c_1$,
\begin{align}\label{e:FcompFinf}
c_2^{-1}\sup_{c_1\le b\le t}\frac{\Phi(b)}{b}\leq \scK_{\infty, a}(t)\leq c_2\sup_{c_1\le b\le t}\frac{\Phi(b)}{b}.
\end{align}
\end{lemma}
\pf \eqref{e:F_11} and \eqref{e:wscF_1} follows from 
Lemma \ref{PhiwtPhi}(3),  \eqref{e:F1} and \eqref{e:wscF}. 

We now prove \eqref{e:FcompFinf}. Without loss of generality, we assume that $c_1<a$.  Let $f(t):=\sup_{c_1\le b\le t}\frac{\Phi(b)}{b}$ for $t\ge c_1$. By Remark \ref{mwsc}, $\Phi$ satisfies $L^{c_1}(\delta, (c_1/a)^{\delta}\wt C_L)$. Thus, for $c_1\le b\le t$,
$
(c_1/a)^{\delta}\wt C_L(t/b)<(c_1/a)^{\delta}\wt C_L(t/b)^{\delta}\le {\Phi(t)}/{\Phi(b)}.
$
This together with Lemma \ref{PhiwtPhi}(1) implies that for $t\ge c_1$, 
\begin{align*}
f(t)\le (a/c_1)^{\delta}\wt C_L^{-1}\frac{\Phi(t)}{t}\le (a/c_1)^{\delta+2}\wt C_L^{-1}\frac{\wt\Phi_a(t)}{t}.
\end{align*}
On the other hand, by Lemma \ref{PhiwtPhi}(1) again,
$
f(t)\ge {\Phi(t)}/{t}\ge {\wt\Phi_a(t)}/{t}
$ for $t\ge c_1$.
Thus, \eqref{e:FcompFinf} follows from \eqref{e:F_11}.
\qed

Suppose $\Phi$ satisfies $L^a(\delta, \wt C_L)$. Then, by Remark \ref{mwsc},  $\Phi$ satisfies $L^1(\delta, \wt C_L')$ for some $\wt C_L'=\wt C_L'(\wt C_L, a, \delta)>0$. Thus, if $\Phi$ satisfies the weak lower scaling property at infinity, we will assume that $\Phi$ satisfies $L^1(\delta, \wt C_L)$ instead of $L^a(\delta, \wt C_L)$. Now we further assume that $\delta>1$. Then, $\scK_{\infty}$ defined in \eqref{d:Finf} is $\scK_{\infty}(s)=\sup_{b\le s}\frac{\wt\Phi(b)}{b}$ with 
$\wt \Phi(t):=\wt \Phi_1(t)=\Phi(1)t^2\1_{\{0<t< 1\}}+\Phi(t)\1_{\{t\ge 1\}}$. $\scK_{\infty}$ is non-decreasing function with $\scK_{\infty}(0)=0$ and $\lim_{t\to\infty}\scK_{\infty}(t)=\infty$. 

In the following lemma, we show some inequalities between $\Phi^{-1}$ and  $\scK^{-1}$,  and between $\Phi^{-1}$ and $\scK_{\infty}^{-1}$.
\begin{lemma}	\label{l:F}
(1) Suppose  $\Phi$ satisfies $L_a(\delta,\wt C_L)$ with $\delta>1$ and for some $a> 0$.
For any $T>0$ and $b >0$ there exists a constant $c_1=c_1(b, \wt C_L, a, \delta,T)>0$ such that 
		\begin{equation}
		\label{e:F3}
		{\Phi^{-1}(t)} \le c_1{\scK^{-1}\left(\frac{t}{b\Phi^{-1}(t)}\right)} \quad \mbox{for all } \, t \in (0,T),
		\end{equation}
and  there exists a constant $C_3=C_3(a,\wt C_L, \delta,T)\geq1$ such that for every $t,r>0$ satisfying $t< \Phi(r)\wedge T$,
\begin{equation}\label{e:F4}
			\left(\frac{r}{\Phi^{-1}(t)}\right)^2 \le \frac{r}{\scK^{-1}(t/r)} \le C_3 \left(\frac{r}{\Phi^{-1}(t)}\right)^{\delta/(\delta-1)}.
			\end{equation}
Moreover, if $a=\infty$, then \eqref{e:F3} and \eqref{e:F4} hold with $T=\infty$. In other words, \eqref{e:F3} holds for all $t<\infty$ and \eqref{e:F4} holds for $t<\Phi(r)$.

(2) Suppose $\Phi$ satisfies $L^1(\delta, \wt C_L)$ with $\delta>1$.  For any $T>0$ and $b>0$ there exists a constant $c_2=c_2(T, b, \wt C_L, \delta)\geq 1$ such that for $t\geq T$, 
\begin{align}\label{e:F_13}
{\Phi^{-1}(t)} \le c_2{\scK_{\infty}^{-1}\left(\frac{t}{b\Phi^{-1}(t)}\right)}, 
\end{align}
and for any $T>0$ there exists a constant $C_4=C_4(a,\wt C_L,\delta,T)\geq1$ such that for every $t,r>0$ satisfying $T\leq t\leq \Phi(r)$,
 \begin{align}\label{e:F_14}
C_4^{-1}\left(\frac{r}{\Phi^{-1}(t)}\right)^2 \le \frac{r}{\scK_{\infty}^{-1}(t/r)} \le C_4 \left(\frac{r}{\Phi^{-1}(t)}\right)^{\delta/(\delta-1)}.
\end{align}
\end{lemma}
\pf (1) By Remark \ref{mwsc}, we may and do assume that $a\geq \Phi^{-1}(T)$. Let $c_1= (b/\wt C_L^2)^{\frac{1}{\delta-1}}\vee 1$. For $t< T\leq \Phi(a)$, using $L^a(\delta,\tilde{C}_L)$ and $c_1 \ge 1$ we have
$$\frac{\Phi\left(\frac{\Phi^{-1}(t)}{c_1}\right)}{t} = \frac{\Phi\left(\frac{\Phi^{-1}(t)}{c_1}\right)}{\Phi(\Phi^{-1}(t))} \le \wt C_L^{-1}c_1^{-\delta}\le\frac{\wt C_L}{bc_1}. $$
Thus, by \eqref{e:F1} we obtain for $t<T$,
$$\scK\left(\frac{\Phi^{-1}(t)}{c_1}\right) \le \wt C_L^{-1}\frac{c_1\Phi\left(\frac{\Phi^{-1}(t)}{c_1}\right)}{\Phi^{-1}(t)} \le \frac{t}{b\Phi^{-1}(t)}, $$
which implies \eqref{e:F3}. 

Now we prove the first inequality in \eqref{e:F4}.  Since $t< \Phi(r)\wedge T\le \Phi(r)$, by \eqref{e:Phi},
\begin{align*} 
\frac{\Phi^{-1}(t)^2}{r^2} \le  \frac{\Phi(\Phi^{-1}(t) \cdot \Phi^{-1}(t)/r)}{\Phi(\Phi^{-1}(t))}.  
\end{align*}
Thus, combining above inequality and \eqref{e:F1} we have 
$$ \frac{t}{r} \le \frac{\Phi( {\Phi^{-1}(t)^2}/{r})}{{\Phi^{-1}(t)^2}/{r}} \le \scK\left( \frac{\Phi^{-1}(t)^2}{r}\right). $$
This concludes the first inequality in \eqref{e:F4}.

To prove the second inequality, let $ C_3:=\wt C_L^{-2/(\delta-1)}\ge1$.  Using the condition $L_a(\delta, \wt C_L)$ on $\Phi$, we have that for $t<T\wedge \Phi(r)$,
\begin{align*}
 C_3^{-1} \left(\frac{\Phi^{-1}(t)}{r}\right)^{\delta/(\delta-1)} 
&\ge  C_3^{\delta-1} \wt C_L\frac{\Phi( C_3^{-1}\Phi^{-1}(t) (\frac{\Phi^{-1}(t)}{r})^{1/(\delta-1)})}{\Phi(\Phi^{-1}(t))} \\
&= \wt C_L^{-1}\frac{\Phi( C_3^{-1} \Phi^{-1}(t)^{\delta/(\delta-1)} r^{-1/(\delta-1)})}{t}.
\end{align*}
Thus, by using the above inequality and \eqref{e:F1} we have that for $t< T\wedge \Phi(r)$,
\begin{equation*}
\frac{t}{r} \ge \wt C_L^{-1} \frac{\Phi(C_3^{-1}\Phi^{-1}(t)^{\delta/(\delta-1)}r^{-1/(\delta-1)}  )}{C_3^{-1}\Phi^{-1}(t)^{\delta/(\delta-1)} r^{-1/(\delta-1)}} \ge  \scK(C_3^{-1}\Phi^{-1}(t)^{\delta/(\delta-1)} r^{-1/(\delta-1)}),
\end{equation*}
which implies 
 the second inequality in \eqref{e:F4}. Since we only assumed $a \ge \Phi^{-1}(T)$ on $T$ and $c_1,C_3$ are independent of $T$, \eqref{e:F3} and \eqref{e:F4} holds with $T=\infty$ when $a=\infty$.

(2) Fix $T_1 \in (0,\infty)$. By Lemma \ref{PhiwtPhi}(3), $\wt \Phi$ satisfies $L(\delta, \wt C_L)$. Now the function $\wt \Phi$ satisfies the assumption of Lemma \ref{l:F}(1), thus \eqref{e:F_13} and \eqref{e:F_14} with functions $\wt \Phi$ and $\scK_{\infty}$ hold with $T=\infty$. Now lemma follows from $\Phi^{-1}(t) \asymp \wt \Phi^{-1}(t)$ for $t \ge T_1$. 
\qed

\section{Near-diagonal estimates and preliminary  upper bound} \label{s:NDE}
\subsection{Functional inequalities}\label{s:PI}

Here we will prove (weak) Poincar\'e inequality with respect to our jumping kernel $J$. We start from a simple calculus.
\begin{lemma}
	\label{l:PI}
	For $r>0$, let $g:(0,r] \rightarrow \R$ be a continuous and non-increasing function satisfying
	\begin{equation*}
	\int_0^r sg(s)ds \ge 0
	\end{equation*}
	and $h:[0,r] \rightarrow [0, \infty)$ be a subadditive measurable function with $h(0)=0$, i.e., 
	\begin{equation}\label{e:pi2}
	h(s_1)+h(s_2) \ge h(s_1 +s_2), \quad \text{ for } 0<s_1, s_2 <r \text{ with } s_1 + s_2 <r. 
	\end{equation} 
	Then,
	\begin{equation}
	\label{e:pi3}
	\int_0^r h(s)g(s)ds \ge 0. 
	\end{equation}
\end{lemma}
\pf Using \eqref{e:pi2} we observe that for any $s \in (0,r)$,
\begin{equation}\label{e:h1}
h(s) = \frac{1}{s}\int_0^s h(s)dt \le \frac{1}{s}\int_0^s \big(h(t) + h(s-t)\big)dt = \frac{2}{s}\int_0^s h(t)dt.
\end{equation}
Let $H(s):= \frac{1}{s^2} \int_0^s h(t)dt$. Then by \eqref{e:h1} 
$$H'(s) = \frac{1}{s^2} \left( h(s)- \frac{2}{s} \int_0^s h(t)dt \right) \le 0.$$
Thus, $H(s)$ is non-increasing. Using this, we have that for any $0<r_1 \le r_2 <r_3 \le r$,
\begin{equation}
\label{e:h5}
\frac{1}{r_1^2} \int_0^{r_1} h(t)dt \ge \frac{1}{r_2^2} \int_0^{r_2} h(t)dt \ge \frac{1}{r_3^2-r_2^2} \int_{r_2}^{r_3} h(t)dt.
\end{equation}
If $g(r) \ge 0$, then \eqref{e:pi3} is trivial since $g$ is non-increasing. Assume $g(r)<0$ and
let $r_0:=\inf \{s \le r: g(s)<0  \}$.
Let $0<k:=\frac{1}{r_0^2} \int_0^{r_0} h(s)ds$. 
By using the continuity of $g$, $g(r_0)=0$, and  the integration by parts, we have
\begin{align*}
\int_0^{r_0} h(s)g(s)ds = -\int_0^{r_0} \int_0^s h(t)dt dg(s) \ge -k \int_0^{r_0} s^2 dg(s) = k\int_0^{r_0} 2sg(s)ds
\end{align*}
and
$$ -\int_{r_0}^r h(s)g(s)ds = -\int_{r_0}^r \int_s^r h(t) dt dg(s) \le -k\int_{r_0}^r (r^2-s^2)dg(s) = -k\int_{r_0}^r 2sg(s)ds. $$
Thus,
\begin{align*}
\int_0^r h(s)g(s)ds = \int_0^{r_0} h(s)g(s)ds + \int_{r_0}^r h(s)g(s)ds \ge k\left(\int_0^{r_0}2sg(s)ds + \int_{r_0}^r 2sg(s)ds\right) \ge 0,
\end{align*}
which completes the proof. \qed

By applying the above lemma, we have the following (weak) Poincar\'e inequality. 
\begin{prop}\label{p:PI}
	There exists $C>0$ such that for every bounded and measurable function $f$, $x_0 \in \R^d$ and $r>0$,
	\begin{equation}
	\label{e:PI}
	\frac{1}{r^d \Phi(r)} \int_{B(x_0,r) \times B(x_0,r)} (f(y)-f(x))^2 dx dy \le C\int_{B(x_0, 3 r) \times B(x_0, 3 r) } (f(y)-f(x))^2 J(x,y)dx dy.
	\end{equation}
\end{prop}
\pf Denote $B(r):=B(x_0,r)$. For $0<s<2r$, let 
\begin{align*}
h(s)&:= \int_{B(3 r -s)} \int_{|z|=s} \big(f(x+z)-f(x) \big)^2 \frac{1}{s^{d}} \sigma(dz)dx, 
\end{align*}
where $\sigma$ is surface measure of the ball. 
We observe that the left hand side of \eqref{e:PI} has the following upper bound:
\begin{align*}
&\frac{1}{r^d \Phi(r)} \int_{B(r) \times B(r)} \big(f(y)-f(x)\big)^2 dx dy\nn \\
&\le \frac{c_1}{r^d \Phi(r)} \int_{B(r)} \int_0^{2r} \int_{|z|=s} \big(f(x+z)-f(x)\big)^2 \sigma(dz) ds dx \nn \\
&\le \frac{c_1}{r^d \Phi(r)}\int_0^{2r} \int_{B(3 r -s)} \int_{|z|=s} \big(f(x+z)-f(x)\big)^2 \sigma(dz) dx ds \nn \\
& = \frac{c_1}{r^d\Phi(r)} \int_0^{2r} h(s)s^d ds \le \frac{2^dc_1 }{\Phi(r)} \int_0^{2r} h(s) ds \le \frac{2^{d+2}c_1}{\Phi(2r)} \int_0^{2r} h(s) ds,
\end{align*}
where  the last inequality follows from \eqref{e:Phi}. 

On the other hand,  the right hand side of \eqref{e:PI} has the following lower bound:
\begin{align*}
& \int_{B(3 r) \times B(3 r) } (f(y)-f(x))^2 J(x,y)dx dy \nn \\
&\ge c_2 \int_{B(3 r) \times B(3 r) } (f(y)-f(x))^2 \frac{1}{|x-y|^d \psi(|x-y|)}dx dy \nn \\
&\ge c_2\int_{B(2r)} \int_{B(3 r - |z|)} (f(x+z)-f(x))^2 \frac{1}{|z|^d \psi(|z|)} dx dz
\nn \\ 
&= c_3 \int_0^{2r} \frac{1}{\psi(s)} \int_{B(3 r -s)} \int_{|z|=s} \big(f(x+z)-f(x)\big)^2  \frac{1}{s^d} \sigma(dz) dx ds \nn \\ 
&= c_3 \int_0^{2r} h(s) \frac{1}{\psi(s)} ds.
\end{align*}
Thus, it suffices to show that 
\begin{equation}
\label{e:PI3}
\int_0^r h(s) \frac{1}{\Phi(r)}  ds \le \int_0^r h(s)  \frac{1}{\psi(s)}ds. 
\end{equation}
To show \eqref{e:PI3}, we will use Lemma \ref{l:PI}.
First note that 
$ \int_0^r \frac{s}{\Phi(r)} ds = \frac{r^2}{2\Phi(r)} = \int_0^r \frac{s}{\psi(s)}ds. $
Let $g(s)= \frac{1}{\psi(s)} - \frac{1}{\Phi(r)}$. Then, $g(s)$ is continuous, non-increasing and
\begin{align}\label{e:pi1check}
\int_0^r sg(s)ds = \int_0^r \frac{s}{\psi(s)} -\frac{s}{\Phi(r)} ds =0.
\end{align}
Now, we show that $h$ is a subadditive function. For $s_1, s_2>0$ with $s_1+s_2<2r$,  
\begin{align}
&h(s_1+s_2) = \int_{B(3 r -s_1-s_2)} \int_{|z|=s_1+s_2} \big( f(x+z)-f(x)\big )^2 \frac{1}{(s_1+s_2)^{d}} \sigma(dz)dx \nn \\
&=  \int_{B(3 r -s_1-s_2)} \int_{|\xi|=1} \frac{\big(f(x+(s_1+s_2)\xi)-f(x)\big)^2}{s_1 + s_2}  \sigma(d\xi)dx \nn \\
&\le \int_{B(3 r -s_1-s_2)} \int_{|\xi|=1}   \frac{\big(f(x+(s_1+s_2)\xi)-f(x+s_2 \xi)\big)^2}{s_1}  +  \frac{\big(f(x+s_2 \xi)-f(x)\big)^2}{s_2}   \sigma(d\xi)dx \nn \\
&\leq \quad \int_{B(x_0+s_2 \xi, 3 r -s_1-s_2)} \int_{|\xi|=1}  \frac{\big(f(x+s_1 \xi)-f(x)\big)^2}{s_1} 
\sigma(d\xi)dx +h(s_2)\nn\\
&  \le \quad  \int_{B(3 r -s_1)}\int_{|\xi|=1}  \frac{\big(f(x+s_1 \xi)-f(x)\big)^2}{s_1}  
\sigma(d\xi)dx+h(s_2)\nn\\
&= h(s_1) + h(s_2), \label{e:h2}
\end{align}
where 
the first inequality follows from the inequality
$\frac{(b_1+b_2)^2}{s_1 +s_2} \le \frac{b_1^2}{s_1}+ \frac{b_2^2}{s_2}$.

Thus, by \eqref{e:pi1check} and \eqref{e:h2}, the functions $g(s)$ and $h(s)$ satisfy the assertions of Lemma \ref{l:PI}. Therefore, by Lemma \ref{l:PI} we have \eqref{e:PI3} which implies \eqref{e:PI}. \qed

\begin{corollary}\label{c:PI}
	There exists a constant $C>0$ such that for any bounded $f \in \FF$ and $r>0$,
	\begin{equation}
	\label{e:cPI}
	\frac{1}{r^d} \int_{\R^d}\int_{B(x,r)} (f(x)-f(y))^2 dydx \le C \Phi(r) \EE(f,f).
	\end{equation}
\end{corollary}

\pf Fix $r>0$ and let $\{x_n\}_{n \in \N}$ be a countable set in $\R^d$ satisfying $\bigcup_{n =1}^\infty B(x_n,r) = \R^d$ and $\sup_{y \in \R^d}  \left| \{ n:y \in B(x_n,6 r) \} \right| \le M$. Then by Proposition \ref{p:PI}, we have
\begin{align*} 
&\frac{1}{r^d}\int_{\R^d}\int_{B(x,r)} (f(x)-f(y))^2 dydx \\
&\le \sum_{n=1}^\infty \frac{1}{r^d}\int_{B(x_n,r)}\int_{B(x,r)} (f(x)-f(y))^2 dydx \\
&\le \sum_{n=1}^\infty \frac{1}{r^d} \int_{B(x_n, 2r) \times B(x_n,2r)} (f(x)-f(y))^2 dydx \\
&\le c_1\sum_{n=1}^\infty \Phi(r) \int_{B(x_n, 6 r) \times B(x_n, 6 r)} (f(x)-f(y))^2 J(x,y) dydx \\
&\le c_1M \Phi(r) \int_{\R^d} \int_{B(x,12 r)} (f(x)-f(y))^2J(x,y) dydx \le c_1M\Phi(r)\EE(f,f). 
\end{align*}
This finishes the proof. \qed
\subsection{Nash's inequality and near-diagonal upper bound in terms of $\Phi$}\label{s:Nash}

In this subsection, using \eqref{e:Phi} and \eqref{e:cPI}, we prove 
Nash's inequality  for $(\sE, \sF)$ and the near-diagonal upper bound of $p(t,x,y)$ in terms of $\Phi$.
The proofs in this subsection are almost identical to the corresponding ones \cite[Section 3]{CK08}.
We provide some details for completeness.

\begin{thm}\label{t:nash}
	There is a positive constant $c>0$ such that for every $u \in \FF$ with $\| u \|_1 =1$, we have
	\begin{equation}\label{e:nash}
	\vartheta(\| u \\|_2^2) \le c \,\EE(u,u),      \qquad \text{ where } \,\, \vartheta(r):= \frac{r}{\Phi(r^{-1/d})}.
	\end{equation}
\end{thm}

\pf 
For $r>0$ and $x \in \R^d$, define
$ u_r(x):= r^{-d} \int_{B(x,r)} u(z) dz$.
Using \eqref{e:cPI} and the inequality 
$
\| u_r \|_2^2 \le \| u_r \|_\infty \| u_r \|_1 \le c_1 r^{-d} \| u \|^2_1$, we have that for $u \in \FF$ with $\| u \|_1 =1$,
\begin{align*}
\| u \|^2_2 &\le 2 \| u -u_r \|^2_2 + 2\| u_r \|^2_2 \\
&\le \frac{2}{ r^d} \int_{\R^d} \int_{B(x,r)} (u(x)-u(y))^2 dydx + \frac{2c_1\| u \|^2_1}{r^d} \le c_2 \left( \Phi(r) \EE(u,u) + \frac{1}{r^d} \right).
\end{align*}
Since $\lim_{r \rightarrow 0}\frac{1}{\Phi(r)r^d} = \infty$ and $\lim_{r \rightarrow \infty}\frac{1}{\Phi(r)r^d} = 0$, 
there exists a constant $r_0$ satisfying $ \EE(u,u) = 1/(\Phi(r_0)r_0^d)$.
Then, by taking $r=r_0$, we have $r_0 \le (2c_2)^{1/d} \| u \|_2^{-2/d}$.  Thus, 
by \eqref{e:Phi}, we get
$$\Phi(r_0)r_0^d \le (2c_2)^{(d+2)/d} \| u \|_2^{-2} \Phi(\| u \|_2^{-2/d}) = \frac{(2c_2)^{(d+2)/d}}{\vartheta(\| u \|_2^2)}, $$
which implies
$ \vartheta(\| u \|_2^2) \le (2c_2)^{(d+2)/d} \EE(u,u). $
\qed

Recall that  $X$ is the  Hunt process corresponding to our Dirichlet form $(\sE, \sF)$ defined in \eqref{e:DF} with  jumping kernel $J$ satisfying \eqref{e:J_psi}.
By using our Nash's inequality \eqref{e:nash} and \cite[Theorem 3.1]{BBCK}, we now show that 
$X$
has a density function $p(t, x, y)$ with respect to  Lebesgue measure, which is  quasi-continuous,
and that the upper bound estimate holds quasi-everywhere.

\begin{thm}\label{t:uhkd}
	There is a properly exceptional set $\sN$ of $X$, a positive symmetric kernel $p(t,x,y)$ defined on $(0,\infty)\times (\bR^d\setminus\sN)\times(\bR^d\setminus\sN)$, and positive constants $C$ depending on $\bar C$ in \eqref{e:J_psi} and $\beta_1, C_L$, 
	such that $p_tf(x):=\E^x[f(X_t)]= \int_{\R^d}p(t,x,y)f(y)dy$, and
	\begin{equation}\label{e:uhkd}
	p(t,x,y) \le \frac{C}{\Phi^{-1}(t)^d}\quad\text{for every}\;\;x,y\in\bR^d\setminus\sN\;\text{and for every}\;t>0.
	\end{equation}
	Moreover, for every $t>0$, and $y\in\bR^d\setminus\sN$,  $x \mapsto p(t,x,y)$ is quasi-continuous on $\R^d$.	
\end{thm}
\pf 
Let  $\vartheta(r):= \frac{r}{\Phi(r^{-1/d})}$.
Since $\Phi$ satisfies $L(\beta_1, C_L)$ by Lemma \ref{l:rel1}, the function $r \mapsto 1/\vartheta(r)$ is integrable at $r=\infty$. Thus by [13, Proposition II.1], Theorem \ref{t:nash} implies that
$
\|P_tf\|_{\infty} \leq m(t)\|f\|_1$ where $m(t)$ is the inverse function of $
h(t):=\int^{\infty}_t\frac{1}{\vartheta(x)}dx.
$
Using $L(\beta_1, C_L)$  of $\Phi$ and following the the proof of \cite[Theorem 3.2]{CK08}, we see that $h(t) \leq  c_1 \Phi(t^{-1/d})$. Now, we apply \cite[Theorem 3.1]{BBCK} and we obtain the thereom.\qed

\subsection{An upper bound of heat kernel using scaling}\label{s:UHK1}
In this section, we observe that the off-diagonal upper bound in \cite[Section 4.1--4.4]{CK08}
holds without the condition (1.14) in \cite{CK08}. We provide full details for reader's convenience.

Recall that  $X$ is the  Hunt process corresponding to our Dirichlet form $(\sE, \sF)$ defined in \eqref{e:DF} with  jumping kernel $J$ satisfying \eqref{e:J_psi}. Fix $\rho>0$ and define  a bilinear form $(\sE^{\rho}, \sF)$ by
\begin{align*}
\sE^{\rho}(u,v)=\int_{\bR^d\times \bR^d} (u(x)-u(y))(v(x)-v(y)){\bf 1}_{\{|x-y|\le \rho\}}\, J(x, y)dxdy.
\end{align*}
Clearly, the form $\sE^{\rho}(u,v)$ is well defined for $u,v\in \sF$, and $\sE^{\rho}(u,u)\le \sE(u,u)$ for all $u\in \sF$. Since $\psi$ satisfies $L(\beta_1, C_L)$ and $U(\beta_2, C_U)$,  for all $u\in \sF$,
\begin{equation}\label{EErhocomp}\begin{split}
\sE(u,u)-\sE^{\rho}(u,u)&= \int(u(x)-u(y))^2{\bf 1}_{\{|x-y|>\rho\}}\,J(x, y)dxdy\\
&\le 4\int_{\bR^d}u^2(x)\,dx\int_{B(x,\rho)^c}J(x,y)\,dy\le \frac{c_0\|u\|_{2}^2 }{\psi(\rho)}.
\end{split}\end{equation}
Thus, $\sE_1 (u, u):= \sE(u,u) + \| u \|_2^2$ is equivalent to $\sE^{\rho}_1(u,u):=
\sE^{\rho}(u,u)+ \|u\|_2^2$ for every $u\in \sF$, which implies that  $(\sE^{\rho},
\sF)$ is a regular Dirichlet form on $L^2(\bR^d, dx)$.
We call $(\sE^{\rho}, \sF)$ the $\rho$-truncated
Dirichlet form. The Hunt process associated with $(\sE^{\rho}, \sF)$ which will be denoted by $X^{\rho}$  can be identified in distribution with the Hunt process of the original Dirichlet form $(\sE, \sF)$ by removing those jumps of size larger than  $\rho$. We use $p^{\rho}(t,x,y)$ to denote the transition density function of $X^{\rho}$.

Note that although the function $\psi$ may not be a correct scale function in our setting, we will still use  $\psi$ to define scaled processes. For $\eta>0$, we define $(X^{(\eta)})_t:=\eta^{-1} X_{\psi(\eta) t}$. Then, $X^{(\eta)}$ is a Hunt process in $\R^d$. We call $X^{(\eta)}$ the $\eta$-scaled process of $X$. 
Let 
\begin{equation*}
\psi^{(\eta)}(r):= \frac{\psi(\eta r)}{\psi(\eta)}, \quad \Phi^{(\eta)}(r):= \frac{r^2}{2\int_0^r \frac{s}{\psi^{(\eta)}(s)}ds},\quad J^{(\eta)}(x,y):= \psi(\eta)\eta^d J(\eta x, \eta y).
\end{equation*}
We emphasize once more that $\psi$ satisfies \eqref{e:intcon}, $L(\beta_1,C_L)$ and $U(\beta_2,C_U)$. Furthermore, by Lemma \ref{l:rel1} we have $\beta_1 < 2$. By definition, $\psi^{(\eta)}$ satisfies $L(\beta_1,C_L)$ and $U(\beta_2,C_U)$ for any $\eta>0$. Also, $J^{(\eta)}$ satisfies 
\begin{align}\label{e:Jeta}
\frac{\bar C^{-1}}{|x-y|^d\psi^{(\eta)}(|x-y|)}\leq  J^{(\eta)}(x,y)\leq \frac{\bar C}{|x-y|^d\psi^{(\eta)}(|x-y|)},\quad x,y\in \R^d,\,\, x\not= y,
\end{align}
where the constant $\bar C>0$ is that of \eqref{e:J_psi}.
Thus, Theorem \ref{t:nash} holds for $\eta$-scaled process $X^{(\eta)}$ with the same constants as $X$. i.e., all constants are independent of $\eta$.  
Furthermore, since
$
\Phi^{(\eta)}(r)= {\Phi(\eta r)}/{\psi(\eta)},
$
Lemma \ref{l:rel1} enables that both $\Phi$ and $\Phi^{(\eta)}$ satisfies $L(\beta_1, C_L)$ and $U(2 \land \beta_2,C_U)$.

Since $J^{(\eta)}(x,y)$ is the jumping kernel of $X^{(\eta)}$, the Dirichlet form $(\sE^{(\eta)}, \sF)$ associated with  $X^{(\eta)}$ satisfies
\begin{equation*}
\sE^{(\eta)}(u,v)=\int_{\R^d \times \R^d} (u(x)-u(y))(v(x)-v(y))J^{(\eta)}(x,y)dxdy.
\end{equation*}
Also, since
$
\P^x(X^{(\eta)} \in A) = \P^x(\eta^{-1} X_{\psi(\eta)t} \in A) = \P^{\eta x} (X_{\psi(\eta) t} \in \eta A),
$
we have 
\begin{equation*}
{p^{(\eta)}}(t,x,y) = \eta^d p(\psi(\eta)t,\eta x,\eta y), \quad \mbox{for a.e. } x,y \in \R^d,
\end{equation*}
where $ p^{(\eta)}(t,x,y)$ is a transition density of $X^{(\eta)}$.
For $\rho>0$, let 
$$ J^{(\eta,\rho)}(x,y)= J^{(\eta)}(x,y) \1_{\{|x-y| \le \rho\}}, \quad J^{(\eta)}_{\rho} (x,y):= J^{(\eta)}(x,y) \1_{\{|x-y| > \rho\}}.$$
Then,
$$ \sE^{(\eta,\rho)}(u,v):=\int_{\R^d \times \R^d} (u(x)-u(y))(v(x)-v(y))J^{(\eta, \rho)}(x,y)dxdy $$
is a $\rho$-truncated Dirichlet form for $X^{(\eta)}$. We use $X^{(\eta,\rho)}$ to denote a Hunt process corresponding to Dirichlet form $(\sE^{(\eta,\rho)},\sF)$ and $p^{(\eta,\rho)}(t,x,y)$ to denote the transition density function of $X^{(\eta,\rho)}$. By the same argument as in \eqref{EErhocomp}, there exists $c>0$ such that any $u\in\sF$
\begin{equation}\label{e:EE}
c \Big(\sE^{(\eta)}(u,u)+\|u\|^{2}_{2}\Big) \le \sE^{(\eta,\rho)}(u,u) +\|u\|^{2}_{2}\le \sE^{(\eta)}(u,u)+\|u\|^{2}_{2}.
\end{equation} 
%
Without loss of generality, we assume that $\psi(1)=1$. Then
$X^{\rho}=X^{(1,\rho)}$,   $J^{\rho}=J^{(1,\rho)}$, $\sE^{\rho}(u,v)=\sE^{(1,\rho)}(u,v)$, and $ p^{\rho}(t,x,y)=p^{(1,\rho)}(t,x,y)$.

Since the constants $\bar C$ in \eqref{e:J_psi} and \eqref{e:Jeta} are same, using \cite[Lemma 3.1]{BGK} we have the following. 
\begin{lemma}\label{l:4.24} 
There exists $c>0$ such that	 for any $\rho>0$, $\eta>0$ and $x,y \in \R^d$,
	$$ p^{(\eta)}(t,x,y) \le  p^{(\eta,\rho)}(t,x,y) + \frac{c\,t}{\rho^d \psi^{(\eta)}(\rho)}.$$
\end{lemma}
In the following we give an upper estimate of $ p^{(\eta,\rho)}(t,x,y)$. It is the counterpart of \cite[Lemma 4.3]{CK08}.
\begin{lemma}
	\label{l:scaling}
	There exists a constant $C>0$, independent of $\eta,\lambda>0$, such that
	\begin{equation}
	\label{e:scaling}
	p^{(\eta,\rho)}(t,x,y) \le  \frac{ Ct}{|x-y|^d \Phi^{(\eta)}(|x-y|)}
	\end{equation}
	for every $\eta>0$, $0<t \le \Phi^{(\eta)}(1)=\Phi(\eta)/\psi(\eta)$, $x,y\in \bR^d\setminus \sN$ with $|x-y| \ge 1$ and $\rho = \frac{\beta_1}{3(d+\beta_1)} |x-y|$.
\end{lemma}

\pf Define $\gamma:=\frac{\beta_1}{3(d+\beta_1)}$ and $\vartheta^{(\eta)}(r):= {r}/{\Phi^{(\eta)}(r^{-1/d})}$ and let $ \EE_1^{(\eta)}(u,u):=\sE^{(\eta)}(u,u)+\|u\|^{2}_{2} $.\
Using \eqref{e:nash} for the process $X^{(\eta)}$ and \eqref{e:EE}, there exists a $c_1>0$ independent of $\eta, \rho$ such that
\begin{align*}
\vartheta^{(\eta)}(\| u \|_2^2) \le c_1 \EE^{(\eta)}(u,u) \le c_1 \EE_1^{(\eta)}(u,u) \le c_1 \left(1+\frac{c_2}{\psi^{(\eta)}(\lambda)}\right)\EE_1^{(\eta,\lambda)}(u,u) \le c_3 \EE_1^{(\eta,\lambda)}(u,u)
\end{align*}
for every $\eta>0$ and $\rho>\gamma$. Note that $\EE^{(\eta,\rho)}_1$ is Dirichlet form with respect to 1-subprocess of $X_t^{(\eta,\rho)}$, i.e. this process has exp$(-t)$ killing. 
We have by \cite[Theorem 3.1]{BBCK} and the same way as that for \cite[Theorem 3.2]{CK08} using the above Nash-type inequality for $\EE_1^{(\eta,\rho)}$, there exists constant $c_4>0$ such that
\begin{equation*}
p^{(\eta,\rho)}(t,x,y) \le \frac{c_4 e^{c_5}}{(\Phi^{(\eta)})^{-1}(t)^d} \quad \mbox{for} \quad 0\le t \le \Phi^{(\eta)}(1), \quad x,y \in \R^d\setminus \sN,
\end{equation*}
for every $\eta>0$ and $\rho \ge \gamma$, since  $\Phi^{(\eta)}(1)={\Phi(\eta)}/{\psi(\eta)} \le c_5$.

 On the other hand, by the condition $L(\beta_1,C_L)$ on $\Phi^{(\eta)}$, we have that for $0 \le t \le \Phi^{(\eta)}(1)$,
\begin{equation}
\label{e:q2}
\frac{1}{[(\Phi^{(\eta)})^{-1}(t)]^d} \le \frac{C_L^{-d/\beta_1}}{[(\Phi^{(\eta)})^{-1}(\Phi^{(\eta)}(1))]^d} \left(\frac{t}{\Phi^{(\eta)}(1)}\right)^{-d/\beta_1} = C_L^{-d/\beta_1}\left(\frac{t}{\Phi^{(\eta)}(1)}\right)^{-d/\beta_1}.
\end{equation} 
Define
\begin{equation*}
\Psi(z):= \frac{s}{3} (|z-x| \land |x-y|) \quad \mbox{for} \quad z \in \R^d,
\end{equation*}
where $s>0$ is a number to be chosen later, and 
$$\Gamma^\eta_\rho[v](x):= \int_{|x-y| \le \rho} (e^{v(x)-v(y)}-1)^2 J^{(\eta)}(x,y)dy. $$
\eqref{e:q2} together with \cite[Theorem 3.2]{BBCK} and \cite[Theorem 3.25]{CKS} implies that  there exist constants $c_6, c_7>0$ such that
\begin{equation}
\label{e:q3}
p^{(\eta,\rho)}(t,x,y) \le c_6 \left(\frac{t}{\Phi^{(\eta)}(1)}\right)^{-d/\beta_1} \mbox{exp}\Big( -|\Psi(x)-\Psi(y)|+c_7t\big(\| \Gamma^\eta_\rho[{\Psi}] \|_\infty \lor \| \Gamma^\eta_\rho[{-\Psi}] \|_\infty\big)  \Big)
\end{equation}
for all $t \in [0,\Phi^{(\eta)}(1)]$ and $x,y \in \R^d\setminus \sN$.
 Observe that 
$(1-e^b)^2 \le b^2 e^{2|b|}$, $b \in \bR$ and $|\Psi(z_1)-\Psi(z_2)| \le \frac{s}{3} |(|z_1-x|-|z_2-x|)| \le \frac{s}{3} |z_1 -z_2|$ for all $z_1, z_2 \in \R^d$. Thus, by the above observation and \eqref{e:Jeta}, we have that for every $z\in \R^d$
\begin{align*}
\Gamma^\eta_\rho[{\Psi}](z) &= \int_{|z-w| \le \rho} (1- e^{\Psi(z)-\Psi(w)})^2 J^{(\eta)}(z,w) dw \\	
&\le \int_{|z-w| \le \rho} (\Psi(z)-\Psi(w))^2 e^{2|\Psi(z)-\Psi(w)|} J^{(\eta)}(z,w)dw \\
&\le \left(\frac{s}{3}\right)^2 e^{2s\rho/3} \int_{|z-w| \le \rho} |z-w|^2 J^{(\eta)}(z,w)dw \\
&\le c_8 s^2 e^{2s\rho/3} \int_0^\rho \frac{t}{\psi^{(\eta)}(t)}dt = c_8 s^2 e^{2s\rho/3} \frac{\rho^2}{\Phi^{(\eta)}(\rho)} \le c_{9} \frac{e^{s\rho}}{\Phi^{(\eta)}(\rho)},
\end{align*}
where the last inequality follows from the inequality $x^2e^{2x/3}\leq 9e^x$ for $x>0$.
Thus, $$ \big\| \Gamma^\eta_\rho[{\Psi}] \big\|_\infty \lor \big\| \Gamma^\eta_\rho[{-\Psi}] \big\|_\infty \le c_{9} \frac{e^{s\rho}}{\Phi^{(\eta)}(\rho)}$$
and for any $\eta>0$, $\rho>\gamma$ and $x,y \in \R^d\setminus \sN$ satisfying $|x-y| \ge 1$ the right hand side of \eqref{e:q3} is bounded by
\begin{equation}
\label{e:q4}
 p^{(\eta,\rho)}(t,x,y) \le c_7 \left(\frac{t}{\Phi^{(\eta)}(1)}\right)^{-d/\beta_1} \mbox{exp} \left( -\frac{s|x-y|}{3} + c_{9} \frac{e^{s\rho} t}{\Phi^{(\eta)}(\rho)} \right).
\end{equation}
Now, take $\rho= \gamma |x-y| \ge \gamma$ and $s = \frac{1}{\gamma |x-y|} \log (\frac{\Phi^{(\eta)}(|x-y|)}{t})$. Then since $\gamma<1$, using \eqref{e:Phi} we have
\begin{align}\label{e:q5}
-\frac{s|x-y|}{3} + c_{9} \frac{e^{s\rho} t}{\Phi^{(\eta)}(\rho)} &= \frac{1}{3\gamma} \log\left(\frac{t}{\Phi^{(\eta)}(|x-y|)}\right) + c_{9} \frac{\Phi^{(\eta)}(|x-y|)}{\Phi^{(\eta)}(\rho)}\nn \\
&\le \frac{1}{3\gamma} \log\left(\frac{t}{\Phi^{(\eta)}(|x-y|)}\right) + \frac{c_{9}}{\gamma^2}.
\end{align}

Using \eqref{e:q4}, \eqref{e:q5}, and the condition $L(\beta_1,C_L)$ on $\Phi^{(\eta)}$, we obtain that for any $\eta>0$, $\rho=\gamma|x-y|$ and $x,y \in \R^d\setminus \sN$ satisfying $|x-y| \ge 1$,
\begin{align*}
p^{(\eta,\rho)}(t,x,y) &\le c_{10}\left(\frac{t}{\Phi^{(\eta)}(1)}\right)^{-d/\beta_1} \left( \frac{t}{\Phi^{(\eta)}(|x-y|)} \right)^{1/3\gamma}\nn \\
&=  c_{10}\left(\frac{t}{\Phi^{(\eta)}(1)}\right)^{-d/\beta_1} \left( \frac{t}{\Phi^{(\eta)}(|x-y|)} \right)^{1+d/\beta_1} \nn \\
&= \frac{c_{10}t\Phi^{(\eta)}(1)^{d/\beta_1}}{\Phi^{(\eta)}(|x-y|)\Phi^{(\eta)}(|x-y|)^{d/\beta_1} } \le \frac{c_{10}C_L^{-1}t}{ \Phi^{(\eta)}(|x-y|)|x-y|^d}.
\end{align*}
\qed

Although we used $\psi$ in scaled process, in the next theorem we are able to obtain an upper bound in terms of  $\Phi$.
\begin{thm}
	\label{t:UHK-Phi}
There exists a constant $C>0$ such that for any $t>0$ and $x,y \in \R^d\setminus \sN$,
	\begin{equation}\label{e:UHK-Phi}
	p(t,x,y) \le C \left( \frac{1}{\Phi^{-1}(t)^d} \land \frac{t}{\Phi(|x-y|)|x-y|^d} \right)
	\end{equation}
\end{thm}
\pf Note that \eqref{e:UHK-Phi} holds when $t,x,y$ satisfies $t \ge \Phi(|x-y|)$ by Theorem \ref{t:uhkd}. Thus, it suffices to show the case $t \le \Phi(|x-y|)$.
By \cite[Lemma 7.2(1)]{HKE}, for every $0<t\le \Phi^{(\eta)}(1)$ and $x,y \in \R^d\setminus \sN$ with $|x-y| \ge 1$,
\begin{equation}\label{e:q7}
p^{(\eta)}(t,x,y) \le p^{(\eta,\rho)}(t,x,y) + t \| J_\rho^{(\eta)} \|_\infty \le p^{(\eta,\rho)}(t,x,y) + c_1 \frac{t}{\psi^{(\eta)}(\rho) \rho^d  }.
\end{equation}
Applying \eqref{e:scaling} in \eqref{e:q7}, and using the condition $U(\beta_2,C_U)$ on $\psi^{(\eta)}$ and the inequality $\Phi^{(\eta)}\leq \psi^{(\eta)}$, we get
\begin{align}\label{e:q8}
p^{(\eta)}(t,x,y) &\le c_2 \frac{ t}{|x-y|^d \Phi^{(\eta)}(|x-y|)} + c_1 \frac{t}{\psi^{(\eta)}(\gamma |x-y|) (\gamma |x-y|)^d}\nn\\
& \le c_3 \frac{ t}{|x-y|^d \Phi^{(\eta)}(|x-y|)},
\end{align}
where $\gamma= \frac{\beta_1}{3(d+\beta_1)}$ is the constant in the proof of Lemma \ref{l:scaling}.

 Taking $\eta=|x-y|$, $t/\psi(\eta) \le \Phi(\eta)/\psi(\eta)= \Phi^{(\eta)}(1)$ and $\eta^{-1}|x-y|=1$. Thus by \eqref{e:q8}, we obtain 
\begin{align*}
p(t,x,y) &= \eta^{-d} p^{(\eta)}(t/\psi(\eta),\eta^{-1}x, \eta^{-1}y)\\
&\le c_3 \eta^{-d} \frac{t/\psi(\eta)}{(\eta^{-1}|x-y|)^d \Phi^{(\eta)}(\eta^{-1}|x-y|)}
 = c_3 \frac{t}{|x-y|^d \Phi(|x-y|)},
\end{align*}
which concludes the proof. \qed

\subsection{Consequences of Poincar\'e inequality and Theorem \ref{t:UHK-Phi}
}\label{s:PHI}
Recall that we always assume that $\psi$ satisfies \eqref{e:intcon}, $L(\beta_1, C_L)$ and $U(\beta_2, C_U)$. The upper bound in \eqref{e:UHK-Phi} may not be sharp. However, there are several important consequences which are induced from \eqref{e:UHK-Phi}.
In this subsection we will apply recent results in 
\cite{HKE, PHI} to \eqref{e:UHK-Phi}.

Using \eqref{e:J_psi}, we immediately see that 
there is a constant $c>0$ such that for all $x, y\in \bR^d$ with $x\not= y$,
\begin{equation}\label{ujs}
J(x,y)\le  \frac{c}{r^d}\int_{B(x,r)}J(z,y)\,dz
\quad\hbox{for every }
0<r\le   |x-y| /2.
\end{equation}
(See \cite[Lemma 2.1]{CKK}).
Such property in \eqref{ujs} is called 
$\text{(UJS)}$ in \cite{CKK}.

\begin{lemma}\label{l:4.2}There exists a constant $C>0$ such that $\P^x(\tau_{B(x,r)} \le t) \le \frac{C t}{\Phi(r)}$ for any $r>0$ and $x \in \R^d\setminus\sN$.
	\end{lemma}
\pf Since we have the upper heat kernel estimates in \eqref{e:UHK-Phi},  the condition $L(\beta_1,C_L)$ on $\Phi$, and conservativeness of $X$, the lemma follows from the same argument as in the proof of \cite[Lemma 2.7]{HKE}. 
\qed
\begin{lemma} \label{l:4.3}
There exist constants $c_1,c_2>0$ such that $c_1 \Phi(r) \le \E^x[\tau_{B(x,r)}] \le c_2 \Phi(r)$ for any $r>0$ and $x \in \R^d\setminus\sN$.

	\end{lemma}
\pf By Lemma \ref{l:4.2}, there exists $b>0$ such that $ \P^x(\tau_{B(x,r)} \le b \Phi(r) ) \le {1}/{2}. $
Thus, $$\E^x[\tau_{B(x,r)}] \ge b \Phi(r) \P^x(\tau_{B(x,r)} > b \Phi(r)) \ge \frac{b}{2} \Phi(r).$$
Observe that for any integer $k>0$ and $t \geq 2^k\Phi(r) $, by using \eqref{e:Phi}, we have $r\Phi^{-1}(t)^{-1}\leq 2^{-2k}$. Thus, using Theorem \ref{t:UHK-Phi} and the above observation, we have
\begin{align*}
&\E^x[\tau_{B(x,r)}]=\int_0^\infty \P^x(\tau_{B(x,r)}>t)dt\leq \int_0^\infty\int_{B(x,r)}p(t,x,y)dydt\\
&=\int_0^{\Phi(r)}\int_{B(x,r)}p(t,x,y)dydt+\sum_{k=0}^{\infty}\int_{2^k\Phi(r)}^{2^{k+1}\Phi(r)}\int_{B(x,r)}p(t,x,y)dydt\\
&\leq \Phi(r) +\sum_{k=0}^{\infty}\int_{2^k\Phi(r)}^{2^{k+1}\Phi(r)}c_4r^d\Phi^{-1}(t)^{-d}dt\leq \Phi(r)+2c_4\Phi(r)\sum_{k=0}^{\infty}\left(\frac{1}{2}\right)^{(2d-1)k} \leq c_5\Phi(r).
\end{align*}
\qed

 For any  open set ${D} \subset \bR^d$,
let $\sF_D:=\{u\in\sF: u=0 \;\text{  q.e. in  }  D^c\}$. Then, $(\sE, \sF_{D})$ is also a regular Dirichlet form. We use $p^{D}(t,x,y)$ to denote the transition density function corresponding to $(\sE, \sF_{D})$.

Recall that $(\sE, \sF)$ is a conservative Dirichlet form. Thus, by Theorem \ref{t:UHK-Phi} and \cite[Theorem 1.15]{HKE}, we see that 
$\text{CSJ}(\Phi)$ defined in \cite{HKE} holds. Thus, by $\text{CSJ}(\Phi)$, $\eqref{ujs}$, \eqref{e:JPhile} and Proposition \ref{p:PI}, we have (7) in  \cite[Theorem 1.19]{PHI}.

Therefore, by \cite[Theorem 1.19]{PHI}, following joint H\"older regularity holds for parabolic functions.
Note that, by a standard argument, we now can take the continuous version of  parabolic functions (for example, see \cite[Lemma 5.12]{GHH18}).
We refer \cite[Definition 1.13]{PHI} for the definition of parabolic functions. Let $Q(t, x,r,R):=(t, t+r)\times B(x,R)$.

\begin{thm}\label{t:PHR}
There exist constant $c>0$, $0<\theta<1$ and $0<\epsilon<1$  such that for all $x_0\in \bR^d$, $t_0\geq 0$, $r>0$ and for every bounded measurable function $u=u(t,x)$ that is parabolic in $Q(t_0,x_0,\Phi(r),r)$,  the following parabolic H\"older regularity holds:
\begin{align*}
|u(s,x)-u(t,y)| \leq c\left(\frac{\Phi^{-1}(|s-t|)+|x-y|}{r} \right)^{\theta}\sup_{[t_0, t_0+\Phi(r)]\times \R^d}|u|
\end{align*}
for every $s,t\in (t_0,t_0+\Phi(\epsilon r))$ and $x,y \in B(x_0,\epsilon r)$.
\end{thm}

Since $p^{D}(t,x,y)$ is parabolic,
from now on, we assume $\sN=\emptyset$ and  take  the joint continuous versions of $p(t,x,y)$ and $p^{D}(t,x,y)$. (c.f., \cite[Lemma 5.13]{GHH18}).)

Again,  by \cite[Theorem 1.19]{PHI} we have the interior near-diagonal lower bound of $p^{B}(t, x,y )$ and parabolic Harnack inequality. 

\begin{thm}\label{t:NDL} 
There exist $\eps\in (0,1)$ and $c_1>0$ such that for any $x_0\in \bR^d$, $r>0$, $0<t\le \Phi(\eps r)$ and $B=B(x_0,r)$,
\begin{align*}
p^{B}(t, x,y )\ge \frac{c_1}{\Phi^{-1}(t)^d},\quad x ,y\in B(x_0,\eps\Phi^{-1}(t)).
\end{align*}
\end{thm}

\begin{thm}\label{t:PHI} 
There exist constants $0<c_1<c_2<c_3<c_4$,  $0<c_5<1$ and  $c_6>0$ such that for every $x_0 \in \bR^d $, $t_0\ge 0$, $R>0$ and for every non-negative function $u=u(t,x)$ on $[0,\infty)\times \bR^d$ that is parabolic on  $Q(t_0, x_0,c_4\Phi(R),R)$,
\begin{align*}
  \sup_{Q_- }u\le c_6 \,\inf_{Q_+}u,
 \end{align*} where $Q_-:=(t_0+c_1\Phi(R),t_0+c_2\Phi(R))\times B(x_0,c_5R)$ and $Q_+:=(t_0+c_3\Phi(R), t_0+c_4\Phi(R))\times B(x_0,c_5R)$.
 \end{thm}

\section{Off-diagonal estimates}\label{s:ODE}

\subsection{Off-diagonal upper heat kernel estimates}\label{s:ODUE}

Recall from the previous section that for $\rho>0$, $(\sE^{\rho}, \sF)$ is  $\rho$-truncated
Dirichlet form of $(\sE, \sF)$. Also, the Hunt process associated with $(\sE^{\rho}, \sF)$ is denoted by  $X^{\rho}$, and $p^{\rho}(t,x,y)$ is the transition density function of $X^{\rho}$.

For any  open set ${D} \subset \bR^d$,
let 
 $\{P_t^{D}\}$ and $\{Q_t^{\rho,{D}}\}$ be the semigroups of $(\sE, \sF_{D})$ and $(\sE^{\rho}, \sF_{D})$, respectively. We write $\{Q_t^{\rho,\bR^d}\}$ as $\{Q_t^{\rho}\}$ for simplicity. We also use $\tau_D^{\rho}$ to denote the first exit time of the process $\{X^{\rho}_t\}$ in $D$.

The following lemma is essential to prove the main result in this subsection.

\begin{lemma}[{\cite[Lemma 5.2]{HKE}}] \label{l:5.2} 
There exist constants $c, C_1, C_2>0$ such that 
	for any $t, \rho>0$ and $x,y \in \R^d$,
	$$p^{\rho}(t,x,y) \le c\Phi^{-1}(t)^{-d} \exp \left( C_1 \frac{t}{\Phi(\rho)} -C_2 \frac{|x-y|}{\rho} \right).$$
\end{lemma}
\pf 
Note that by Lemma \ref{l:rel1}, $\Phi$ satisfies $U(\beta_2 \land 2, C_U)$ and $L(\beta_1, C_L).$
 By Theorem \ref{t:uhkd}, \eqref{e:JPhile}, and Lemma \ref{l:4.3}, the assumptions of \cite[Lemma 5.2]{HKE} are satisfied. Thus, the lemma follows.
\qed 

The next lemma was proved in \cite[Lemma 7.11]{HKE}  and \cite[Theorem 3.1]{GHL14}  under the assumption that 
$\phi(r,\cdot)$ is  non-decreasing for all $r>0$. We will prove the lemma without such assumption.  

\begin{lemma}\label{l:7.11}
Let  $r, t, \rho>0$. Assume that 
		\begin{align}\label{l:7.11ass}\P^w(\tau_{B(x,r)}^{\rho} \le t) \le \phi(r,t) \quad \mbox{for all }x\in \R^d, \; w \in B(x,r/4),\end{align}
where $\phi$ is a non-negative measurable function on $\R_+\times\R_+ $. Then, for any integer $k\geq 1$,	
$$Q_t^{\rho} \1_{B(x,k(r+\rho))^c}(z) \le \phi(r,t)^k \quad \mbox{for all }x\in \R^d,  z \in B(x,r/4). $$
\end{lemma}
\pf
Assume that $r, t,\rho>0$ satisfy \eqref{l:7.11ass} and fix $x \in \R^d$. Note that $X^{\rho}_{\tau_{B(x,r)}^{\rho}}=X^{\rho}({\tau_{B(x,r)}^{\rho}})\in
B(x,r+\rho)$, and $|w-y|\ge |x-w|-|y-x|\ge k(r+\rho)$ for any $w\notin B(x,(k+1)(r+\rho))^c$ and $y\in
B(x,r+\rho)$. Thus by the strong Markov property, for all $s\leq t$ and $z\in B(x,r/4)$ we have
\begin{align*}
Q_s^{\rho}\1_{B(x,(k+1)(r+\rho))^c}(z)&=\E^z\left[\1_{\{\tau^{\rho}_{B(x,r)}<s \}}\P^{X^{\rho}({\tau^{\rho}_{B(x,r)})}}\Big(X^{\rho}({s-\tau^{\rho}_{B(x,r)}}) \notin B(x,(k+1)(r+\rho))\Big)  \right]\\
&\leq \P^z(\tau^{\rho}_{B(x,r)}<s)\sup_{y\in B(x,r+\rho),s_1\leq s}Q_{s_1}^{\rho}\1_{B(y,k(r+\rho))^c}(y)\\
&\leq \P^z(\tau^{\rho}_{B(x,r)} \le t)\sup_{y\in B(x,r+\rho),s_1\leq s}Q_{s_1}^{\rho}\1_{B(y,k(r+\rho))^c}(y)\\
&\leq \phi(r,t)\sup_{y \in \R^d, s_1\leq t}Q_{s_1}^{\rho}\1_{B(y,k(r+\rho))^c}(y).
\end{align*}
Thus, by using the above step $k-1$ times we conclude
\begin{align*}
Q_t^{\rho} \1_{B(x,k(r+\rho))^c}(z)  &\leq \phi(r,t)\sup_{y_1 \in \R^d,  s_1\leq t}Q_{s_1}^\rho\1_{B(y_1,(k-1)(r+\rho))^c}(y_1)\\
&\leq \phi(r,t)^2\sup_{y_2 \in \R^d, s_2\leq t}Q_{s_2}^\rho\1_{B(y_2,(k-2)(r+\rho))^c}(y_2)\\
&\leq \cdots\\
&\le\phi(r,t)^{k-1}\sup_{y_{k-1} \in \R^d, s_{k-1}\leq t}Q_{s_{k-1}}^{\rho}\1_{B(y_{k-1},r+\rho)^c}(y_{k-1})\\
&\leq \phi(r,t)^{k-1}\sup_{y_{k-1} \in \R^d, s_{k-1}\leq t}\P^{y_{k-1}}(\tau^{\rho}_{B(y_{k-1},r)}\leq t) \leq \phi(r,t)^k
\end{align*}
for all $z \in B(x,r/4)$.
\qed

The following lemma is a key to obtain upper bound of transition density function and will be used in several times.  

\begin{lemma}\label{l:tail}
Let $f:\bR_+\times\bR_+\to\bR_+$ be a measurable function satisfying that $t\mapsto f(r, t)$ is non-increasing for all $r>0$ and that  $r\mapsto f(r, t)$ is non-decreasing for all $t>0$. Fix $T\in(0,\infty]$. Suppose that  the following hold:
\begin{enumerate}
\item[(i)]  For each $b>0$, $\sup_{ t \le T} f(b\Phi^{-1}(t), t) < \infty$ (resp., $\sup_{ t \ge T} f(b\Phi^{-1}(t), t) < \infty$);
\item[(ii)] there exist $\eta\in(0,\beta_1]$, $a_1>0$ and $c_1>0$ such that
\begin{equation}\label{e:tail}
\int_{B(x,r)^c} p(t, x,y)dy\le c_1\bigg(\frac{\psi^{-1}(t)}{r}\bigg)^\eta+c_1\exp{\Big(-a_1f(r, t)\Big)}
\end{equation} 
for all $t\in (0,T)$ (resp. $t\in [T,\infty)$) and any $r>0$, $x \in \R^d$.
\end{enumerate}
 Then, there exist constants $k, c_0>0$ such that 
\begin{equation*}
p(t,x,y) \le \frac{c_0\, t}{|x-y|^d\psi(|x-y|)} + c_0\,\Phi^{-1}(t)^{-d} \exp{\Big(-a_1k f(|x-y|/(16k), t)\Big)} 
\end{equation*}
for all $t\in (0,T)$ (resp. $t\in [T,\infty)$) and $x, y\in \R^d$.\end{lemma}

\pf
 Since the proofs for the case $t \in (0, T)$ and the case $t \in [T,\infty)$ are similar, we only prove for $t\in(0,T)$. For $x_0 \in \R^d$, let $B(r)=B(x_0,r)$.  By the strong Markov property, \eqref{e:tail}, and the fact that $t\mapsto f(r, t)$ is non-increasing,
 we have that for $x \in B(r/4)$ and $t\in (0,T/2)$,
\begin{align}
\begin{split}\label{cac1}
\P^x(\tau_{B(r)} \le t) 
&=\P^x(\tau_{B(r)} \le t, X_{2t} \in B(r/2)^c) +\P^x(\tau_{B(r)} \le t, X_{2t} \in B(r/2))\\ 
&\le \P^x(X_{2t} \in B(r/2)^c) + \sup_{z \in B(r)^c, s \le t} \P^z(X_{2t-s} \in B(z,r/2)^c) \\
&\le \P^x(X_{2t} \in B(x, r/4)^c) + \sup_{s \le t} \P^z(X_{2t-s} \in B(z,r/4)^c) \\
&\le c_1\bigg(\frac{\psi^{-1}(2t)}{r/4}\bigg)^\eta+c_1\exp{\Big(-a_1f(r/4, 2t)\Big)}.
\end{split}
\end{align}
From this and Lemma \ref{lem:inverse}, we have that for  $x\in B(r/4)$ and $t\in(0, T/2)$,
\begin{equation}\label{e:tail2}
1-P_t^B\1_B(x)= \P^x(\tau_B \le t) \le 
c_2\bigg(\frac{\psi^{-1}(t)}{r}\bigg)^\eta+c_1\exp{\Big(-a_1f(r/4, 2t)\Big)}.
\end{equation}
By \cite[Proposition 4.6]{GHL14} and \eqref{int_outball}, letting $\rho=r$ we have
\begin{align*}
  \Big|P_t^{B(r)} \1_{B(r)}(x)-Q_t^{r,{B(r)}}\1_{B(r)}(x)\Big|
 &\le2t\,\esssup_{z\in \bR^d}\int_{B(z,r)^c}J(z,y)\,dy\le \frac{c_3 t}{\psi(r)}.
 \end{align*}
Combining this with \eqref{e:tail2}, we see that for all $x\in  B(r/4)$ and $t\in (0,T/2)$,
\begin{align}\label{e:tail3}
\P^x(\tau^{r}_{B(r)} \le t) &= 1-Q_{t}^{r,B(r)} \1_{B(r)}(x)
 \le 1-P_t^{B(r)}\1_{B(r)}(x) + \frac{c_3t}{\psi(r)}\nn\\
&\le c_2\bigg(\frac{\psi^{-1}(t)}{r}\bigg)^\eta+c_1\exp{\Big(-a_1f(r/4, 2t)\Big)}+\frac{c_3t}{\psi(r)}
=: \phi_1(r,t).
\end{align}
Applying Lemma \ref{l:7.11} with $r=\rho$ to \eqref{e:tail3}, we see that for $t\in (0,T/2)$
\begin{equation}\label{e:tail4}
\int_{B(x,2kr)^c} p^{r}(t,x,y)dy =Q_t^{r} \1_{B(x,2kr)^c}(x)  \le \phi_1(r,t)^k.
\end{equation}
Let $k=\ceil{ \frac{\beta_2+d}{\eta}}$. 
For $t\in(0, T)$ and $x,y \in \R^d$ satisfying $4k\Phi^{-1}(t) \ge |x-y|$, by using that $r \mapsto f(r,t)$ is non-decreasing and the assumption $(i)$, we have
$f(|x-y|/(16k), t)\le f(\Phi^{-1}(t)/4,t)\le M<\infty $. Thus, by Theorem \ref{t:uhkd}, 
\begin{equation}
\label{e:uhk1}
p(t,x,y) \le c_5\Phi^{-1}(t)^{-d} \le c_5 e^{a_1kM} \Phi^{-1}(t)^{-d} \exp{\Big(-a_1kf(|x-y|/(16k), t)\Big)}.
\end{equation} 

For the remainder of the proof, assume $t\in(0, T)$ and $4k\Phi^{-1}(t) < |x-y|$, and let $r=|x-y|$ and $\rho={r}/{(4k)}$. By \eqref{e:tail4} and Lemmas \ref{lem:inverse}, \ref{l:rel1} and \ref{l:5.2}, we have
\begin{align}
p^{\rho}(t,x,y) &= \int_{\R^d} p^{\rho}(t/2,x,z) p^{\rho}(t/2,z,y)dz\nn \\
&\le \left( \int_{B(x,r/2)^c} + \int_{B(y,r/2)^c} \right) p^{\rho}(t/2,x,z) p^{\rho}(t/2,z,y)dz \nn\\
&\le 2\sup_{z \in \R^d} p^{\rho}(t/2,z,y) \int_{B(x,2k\rho)^c} p^{\rho}(t/2,x,z)dz\nn\\
&\le 
c\Phi^{-1}(t/2)^{-d} \exp \left( C_1 \frac{t}{2\Phi(\rho)} \right)
 \phi_1(\rho,t/2)^k\nn\\
 &\le 
c\Phi^{-1}(t/2)^{-d} \exp \left( C_1 \frac{\Phi(r/4k)}{2\Phi(r/2k)} \right)
 \phi_1(\rho,t/2)^k\nn\\
 &\le c_6\Phi^{-1}(t)^{-d} \phi_1(\rho,t/2)^k. \label{e:rfee}
\end{align}
Note that $k\beta_1  \ge k\eta  \ge \beta_2 + d$, and $\rho \ge \Phi^{-1}(t) > \psi^{-1}(t)$. Thus, by $L(\beta_1, C_L)$ on $\psi$, 
\begin{align*}
\left(\frac{\psi^{-1}(t)}{\rho}\right)^{\eta k}+\left(\frac{t}{\psi(\rho)}\right)^k  
 \le c_7 \left(\left(\frac{\psi^{-1}(t)}{\rho}\right)^{\beta_2 +d} +\left(\frac{\psi^{-1}(t)}{\rho}\right)^{k \beta_1} \right)  
 \le c_8\left(\frac{\psi^{-1}(t)}{r}\right)^{\beta_2 +d}.
\end{align*}
Using Lemmas \ref{lem:inverse} and \ref{l:rel1},
\begin{align*}
 \Phi^{-1}(t)^{-d}\left(\left(\frac{\psi^{-1}(t)}{\rho}\right)^{\eta k}+\left(\frac{t}{\psi(\rho)}\right)^k  \right)
 &\le \frac{c_{8}}{r^d} \frac{\psi^{-1}(t)^{d}}{\Phi^{-1}(t)^d} \left(\frac{\psi^{-1}(t)}{r} \right)^{\beta_2}\\
 &\le \frac{c_{8}}{r^d}  \left(\frac{\psi^{-1}(t)}{\psi^{-1}(\psi(r))} \right)^{\beta_2}\le\frac{c_{9} t}{r^d\psi(r)}.\end{align*}
Applying to \eqref{e:rfee} we have 
\begin{align*}
p^{\rho}(t,x,y)&\le c_{10} \Phi^{-1}(t)^{-d} \left(\left(\frac{\psi^{-1}(t)}{\rho}\right)^{\eta k}+\exp{\Big(-a_1kf(\rho/4, t)\Big)}+\left(\frac{t}{\psi(\rho)}\right)^k \right) \\
  &\le \frac{c_{11} t}{r^d\psi(r)} +  c_{10}\Phi^{-1}(t)^{-d}\exp{\Big(-a_1kf(r/(16k), t)\Big)}.
\end{align*} 
Thus, by Lemma \ref{l:4.24} and $U(\beta_2, C_U)$ on $\psi$, 
we have 
\begin{align}
\label{e:uhk2}
p(t,x,y) &\le p^{\rho}(t,x,y) + \frac{c_{12}t}{\rho^d \psi(\rho)}\nn\\
& \le \frac{c_{13} t}{|x-y|^d\psi(|x-y|)} +  c_{13}\Phi^{-1}(t)^{-d}\exp{\Big(-a_1kf(r/(16k), t)\Big)}.
\end{align}
Now the lemma follows immediately from \eqref{e:uhk1} and \eqref{e:uhk2}.
\qed

The following inequality will be used several times in the proofs of this section:
For any $c_0>0$ and $\alpha \in (0,1)$, there exists $c_1=c_1(c_0, \alpha)>0$ such that $2n \le \frac{c_0}{2d} 2^{n(1-\alpha)} + c_1$ holds for every $n\ge0$. Thus, for any $n\ge0$ and $\kappa \ge1$,
\begin{align}\label{e:3.4.1}
& 2^{nd} \exp\left(-c_02^{n(1-\alpha)}\kappa \right)\le  2^{-nd}\exp\left(2nd-c_02^{n(1-\alpha)}\kappa \right)\nn \\ 
 &\le 2^{-nd}\exp\left( \left(\frac{c_0}{2d} 2^{n(1-\alpha)} + c_1 \right)d-c_02^{n(1-\alpha)}\kappa \right) \le 2^{-nd}\exp\left( \frac{c_0}{2} 2^{n(1-\alpha)}\kappa  + c_1d -c_02^{n(1-\alpha)}\kappa \right) \nn\\
& = e^{c_1d} 2^{-nd} \exp \left(- \frac{c_0}{2}\kappa  \right).
 \end{align}
The next proposition is an intermediate step toward to Theorem \ref{t:exp-1}.
\begin{prop}\label{p:exp}
	There exist constants $a_1,C>0$ and $N \in \N$ such that 
	\begin{equation}\label{e:expN}
	p(t,x,y) 
	\le \frac{C\, t}{|x-y|^d\psi(|x-y|)} + C\,\Phi^{-1}(t)^{-d} \exp{\left(-\frac{a_1|x-y|^{1/N}}{\Phi^{-1}(t)^{1/N}}\right)}, 
	\end{equation}
	for all $t>0$ and $x,y\in\bR^d$.
\end{prop}

\pf Fix $\alpha \in (d/(d +\beta_1),1)$ and let  $N:=\ceil{\frac{\beta_1 + d}{\beta_1}}+1$, and $\eta:=\beta_1-(\beta_1+d)/N>0$.  

 We first claim that there exist $a_2>0$ and $c_1>0$ such that
\begin{equation}\label{e:tailN}
\int_{B(x,r)^c} p(t, x,y)dy\le c_2\left(\frac{\psi^{-1}(t)}{r}\right)^\eta+c_1\exp{\left(-\frac{a_2r^{1/N}}{\Phi^{-1}(t)^{1/N}}\right)},
\end{equation} 
 for any $t,r>0$ and $x \in \R^d$.

When $r \le \Phi^{-1}(t)$, we immediately obtain \eqref{e:tailN} by letting $c_1 = \exp(a_2)$. Thus, we will only consider the case $r>\Phi^{-1}(t)$. 
For any $\rho,t>0$ and all $x,y \in \R^d$, by Lemmas \ref{l:4.24} and \ref{l:5.2} we have
\begin{align}\label{e:p-tru}
p(t,x,y)\leq c_2\Phi^{-1}(t)^{-d} \exp \left( C_1 \frac{t}{\Phi(\rho)} - C_2 \frac{|x-y|}{\rho} \right) +\frac{c_2t}{\rho^d\psi(\rho)}.
\end{align}
Define
$$
\rho_n=\rho_n(r,t)=2^{n\alpha}r^{1-1/N}\Phi^{-1}(t)^{1/N}, \quad n\in \N.
$$
Since $r>\Phi^{-1}(t)$ , we have
$
\Phi^{-1}(t) < \rho_n \leq 2^nr.
$
In particular, $t < \Phi(\rho_n)$. Thus, by \eqref{e:p-tru} we have that for every $t>0$ and $2^nr \le |x-y| < 2^{n+1}r$,
\begin{align*}
p(t,x,y) &\le c_2\Phi^{-1}(t)^{-d} \exp \left( C_1 \frac{t}{\Phi(\rho_n)} - C_2 \frac{|x-y|}{\rho_n} \right) +\frac{c_2t}{\rho_n^d\psi(\rho_n)} \\
&\le c_3 \Phi^{-1}(t)^{-d} \exp\left(-C_2 \frac{2^n r}{\rho_n}\right) + \frac{c_2t}{\rho_n^d \psi(\rho_n)} \\
&= c_3 \Phi^{-1}(t)^{-d} \exp\left(-C_2 \frac{2^{n(1-\alpha)}r^{1/N}} {\Phi^{-1}(t)^{1/N}} \right) + \frac{c_2t}{\rho_n^d \psi(\rho_n)}.
\end{align*}
Using the above estimate we get that
\begin{align*}
&\int_{B(x,r)^c}p(t,x,y)dy
=\sum_{n=0}^{\infty}\int_{B(x,2^{n+1}r)\setminus B(x,2^nr)}p(t,x,y)dy\\
&\leq c_4\sum_{n=0}^{\infty} (2^n r)^d \Phi^{-1}(t)^{-d} \exp\left(-C_2 \frac{2^{n(1-\alpha)}r^{1/N}} {\Phi^{-1}(t)^{1/N}} \right) + c_4\sum_{n=0}^\infty (2^n r)^d \frac{t}{\rho_n^d \psi(\rho_n)}=:I_1+I_2.
\end{align*}
We first estimate $I_1$. Using $\Phi^{-1}(t) < r$,  \eqref{e:3.4.1}, and the fact that $\sup_{s \ge 1} s^d \exp(-\frac{C_2}{4} s^{1/N}) < \infty$,
\begin{align}
I_1
&= c_4\sum_{n=0}^{\infty}\left(\frac{r}{\Phi^{-1}(t)}\right)^{d}2^{nd}\exp\left(-C_2\frac{2^{n(1-\alpha)}r^{1/N}}{\Phi^{-1}(t)^{1/N}}\right)\nn\\
&\leq c_4e^{c_1d}\left(\frac{r}{\Phi^{-1}(t)}\right)^{d}\exp\left(-\frac{C_2}{2}\frac{r^{1/N}}{\Phi^{-1}(t)^{1/N}}\right)\sum_{n=0}^{\infty}2^{-nd}\nn\\
&\leq c_5\exp\left(-\frac{C_2 r^{1/N}}{4\Phi^{-1}(t)^{1/N}}\right). \label{e:wq1}
\end{align}
We next estimate $I_2$. By \eqref{comp1},  $t<\Phi(\rho_n)<\psi(\rho_n)$ and $\Phi^{-1}>\psi^{-1}$, $L(\beta_1, C_L)$ on $\psi$ and $\alpha({d+\beta_1})>{d}$,  we have
\begin{align*}
I_2 & =c_4\sum_{n=0}^\infty\frac{(2^nr)^d}{\rho_n^d}\frac{\psi(\psi^{-1}(t))}{\psi(\rho_n)}
\leq c_4C_L^{-1}\sum_{n=0}^\infty \left(\frac{2^nr}{\rho_n}\right)^{d}\left(\frac{\psi^{-1}(t)}{\rho_n}\right)^{\beta_1}
\\&= c_4C_L^{-1}\left(\frac{\Phi^{-1}(t)}{r}\right)^{-\frac{d+\beta_1}{N}}\left(\frac{\psi^{-1}(t)}{r}\right)^{\beta_1}\sum_{n=0}^\infty 2^{n(d-\alpha(d + \beta_1))}\\
&= c_6\left(\frac{\Phi^{-1}(t)}{r}\right)^{-\frac{d+\beta_1}{N}} \left(\frac{\psi^{-1}(t)}{r}\right)^{\beta_1} \leq c_6\left(\frac{\psi^{-1}(t)}{r}\right)^{\beta_1-\frac{d+\beta_1}{N}}
=c_6\left(\frac{\psi^{-1}(t)}{r}\right)^{\eta}.
\end{align*}
Thus, by above estimates of $I_1$ and $I_2$, we obtain \eqref{e:tailN}.

By $\eta<\beta_1$ and \eqref{e:tailN}, assumptions in Lemma \ref{l:tail} hold with $f(r,t):=\big(r/\Phi^{-1}(t)\big)^{1/N}$. 
Now \eqref{e:expN} follows from Lemma \ref{l:tail}.
\qed
 
By using Proposition \ref{p:exp}, we obtain the upper bound in Theorem \ref{t:exp-1}.

\vspace{3mm}

\noindent{\bf Proof of \eqref{e:exp-1}.} 
By Proposition \ref{p:exp}, there are $a_0, c_0>0$ and $N\in\bN$ such that
\begin{equation*}
	p(t,x,y) 
	\le \frac{c_0\, t}{|x-y|^d\psi(|x-y|)} + c_0\,\Phi^{-1}(t)^{-d} \exp{\left(-\frac{a_0|x-y|^{1/N}}{\Phi^{-1}(t)^{1/N}}\right)}, 
	\end{equation*}
for all $t>0$ and $x,y\in\bR^d$.
Similar to the proof of Proposition \ref{p:exp}, we will show that there exist $a_1>0$ and $c_1>0$  such that for any $t>0$ and $r>0$,
\begin{equation}\label{tail}
\int_{B(x,r)^c} p(t, x,y)\,dy \le c_1\bigg(\frac{\psi^{-1}(t)}{r}\bigg)^{\beta_1/2}+c_1\exp{\left(-\frac{a_1r^2}{\Phi^{-1}(t)^2}\right)}.
\end{equation}
Let $\theta:=\frac{\beta_1}{4d+3\beta_1} \in (0, 1)$ and $C_0=\frac{2C_1}{C_2}$, where $C_1$ and $C_2$ are the constants in Lemma \ref{l:5.2}. Without loss of generality, we may and do assume that $C_0\ge1$. Firstly, when $r \le C_0\Phi^{-1}(t)$ we have
\begin{equation}
\label{tail1}
\int_{B(x,r)^c} p(t,x,y)dy \le 1 \le  e^{a_0C_0^2} \exp{\left(-\frac{a_0r^2}{\Phi^{-1}(t)^2}\right)}. 
\end{equation}
Secondly, we consider the case $r > C_0\frac{\Phi^{-1}(t)^{1+\theta}}{\psi^{-1}(t)^\theta}$. For $|x-y| > r$, there is a $\theta_0 \in (\theta,\infty)$ such that $|x-y|=C_0\Phi^{-1}(t)^{1+\theta_0}/\psi^{-1}(t)^{\theta_0}$. Note that there exists a positive constant $c_2=c_2(\theta)$ such that $s^{-d - \beta_2 - \beta_2/\theta} \ge c_2\exp(-a_0 s^{1/N})$ for $s  \ge1$. Thus, $$\left(\left(\frac{\Phi^{-1}(t)}{\psi^{-1}(t)}\right)^{\theta_0}\right)^{-d - \beta_2 - \beta_2/\theta} \ge c_2\exp\left(-a_0 \left(\frac{\Phi^{-1}(t)}{\psi^{-1}(t)}\right)^{{\theta_0}/N}\right)$$ holds for all $\theta_0 \in (\theta,\infty)$. Using this and $U(\beta_2, C_U)$ condition on $\psi$,
\begin{align}\label{jgeqexp}
\frac{ t}{|x-y|^d\psi(|x-y|)}
&= C_0^{-d}\Phi^{-1}(t)^{-d} \left(\frac{\psi^{-1}(t)}{\Phi^{-1}(t)}\right)^{d\theta_0} \frac{\psi(C_0\psi^{-1}(t)^{1+\theta_0} / \psi^{-1}(t)^{\theta_0})}{\psi(C_0\Phi^{-1}(t)^{1+\theta_0} / \psi^{-1}(t)^{\theta_0})} \frac{\psi(\psi^{-1}(t))}{\psi(C_0 \psi^{-1}(t))}\nn \\
&\ge C_0^{-d-\beta_2}C_U^{-2}\Phi^{-1}(t)^{-d}\left(\frac{\psi^{-1}(t)}{\Phi^{-1}(t)}\right)^{d\theta_0} \left(\frac{\psi^{-1}(t)}{\Phi^{-1}(t)}\right)^{(1+\theta_0)\beta_2}\nn \\
&= C_0^{-d-\beta_2}C_U^{-2}\Phi^{-1}(t)^{-d}\left(\left(\frac{\Phi^{-1}(t)}{\psi^{-1}(t)}\right)^{\theta_0}\right)^{-d-\beta_2-\beta_2/\theta_0} \nn\\
&\ge c_2C_0^{-d-\beta_2}C_U^{-2}\Phi^{-1}(t)^{-d}\exp{\left(-\frac{a_0\Phi^{-1}(t)^{\theta_0/N}}{\psi^{-1}(t)^{\theta_0/N}}\right)}\nn\\
&= c_2C_0^{-d-\beta_2}C_U^{-2}\Phi^{-1}(t)^{-d} \exp{\left(-\frac{a_0|x-y|^{1/N}}{\Phi^{-1}(t)^{1/N}}\right)}.
\end{align}
Thus, for $|x-y|>r$,
\begin{align*}
p(t,x,y) &\le \frac{c_0\, t}{|x-y|^d\psi(|x-y|)} + c_0\Phi^{-1}(t)^{-d} \exp{\left(-\frac{a_0|x-y|^{1/N}}{\Phi^{-1}(t)^{1/N}}\right)} \le \frac{c_3\, t}{|x-y|^d\psi(|x-y|)}.
\end{align*}
Using this, \eqref{int_outball}, $L(\beta_1,C_L)$ condition on $\psi$, and the fact that  $r>C_0\psi^{-1}(t)$ which follows from \eqref{comp1},  
\begin{align}
\int_{B(x,r)^c}p(t,x,y)dy &\le c_3 \int_{B(x,r)^c} \frac{t}{|x-y|^d \psi(|x-y|)} dy 
\le c_4 \frac{t}{\psi(r)}\nn\\
&\le c_4C_L^{-1}\left(\frac{\psi^{-1}(t)}{r}\right)^{\beta_1} \le c_4C_L^{-1}\left(\frac{\psi^{-1}(t)}{r}\right)^{\beta_1/2}. \label{tail2}
\end{align}

Now consider the case $C_0\Phi^{-1}(t) < r \le C_0\Phi^{-1}(t)^{1+\theta} / \psi^{-1}(t)^{\theta}$. In this case, there exists $\theta_0 \in (0,\theta]$ such that $r= C_0\Phi^{-1}(t)^{1+\theta_0}/\psi^{-1}(t)^{\theta_0}$. Define $\rho_n = C_02^{n\alpha} \Phi^{-1}(t)^2/r$, where $\alpha \in (d/(d +\beta_1),1)$. Since $C_0 = \frac{2C_1}{C_2}$, using \eqref{e:Phi} 
\begin{align*}
 \frac{C_1 t}{\Phi(\rho_n)} -  \frac{C_2 2^n r}{\rho_n} 
&\le 
  \frac{C_1\Phi(\Phi^{-1}(t))}{\Phi\left(\Phi^{-1}(t) \frac{C_0\Phi^{-1}(t)}{r}\right)} -  \frac{C_2 2^{n(1-\alpha)}r^2}{C_0\Phi^{-1}(t)^2} \\
&\le C_1 \left(\frac{\Phi^{-1}(t)}{\Phi^{-1}(t)\frac{C_0\Phi^{-1}(t)}{r}}\right)^2 - C_2 2^{n(1-\alpha)} \frac{r^2}{C_0\Phi^{-1}(t)^2} 
\le \left(\frac{C_1}{C_0^2} - \frac{C_2}{C_0} 2^{n(1-\alpha)}\right) \frac{r^2}{\Phi^{-1}(t)^2}\\
& \le -\frac{C_2}{2C_0}2^{n(1-\alpha)} \frac{r^2}{\Phi^{-1}(t)^2}
=: -a_2 2^{n(1-\alpha)} \frac{r^2}{\Phi^{-1}(t)^2}.
\end{align*}
By the above inequality, Lemma \ref{l:4.24}, and Lemma \ref{l:5.2},  for $2^n r \le |x-y| < 2^{n+1}r$,
\begin{align*}
p(t,x,y) &\le \frac{c_5\,t}{\rho_n^d \psi(\rho_n)} + c_5\Phi^{-1}(t)^{-d} \exp\left(C_1 \frac{t}{\Phi(\rho_n)} -C_2 \frac{2^n r}{\rho_n}\right) \nn \\
 &\le \frac{c_5\,t}{\rho_n^d \psi(\rho_n)} + c_5\Phi^{-1}(t)^{-d} \exp \left( -a_22^{n(1-\alpha)} \frac{r^2}{\Phi^{-1}(t)^2} \right).
\end{align*}
Using this, we have
\begin{align*}
&\int_{B(x,r)^c} p(t,x,y)dy
= \sum_{n=0}^\infty \int_{B(x,2^{n+1}r) \backslash B(x,2^n r)} p(t,x,y) dy \\
&\le c_6\sum_{n=0}^\infty (2^nr)^d \Phi^{-1}(t)^{-d} \exp\left(-a_2 \frac{2^{n(1-\alpha)} r^2}{\Phi^{-1}(t)^2}\right)
+c_6\sum_{n=0}^\infty (2^nr)^d \left(\frac{t}{\rho_n^d \psi(\rho_n)}\right) \\
&=: c_6(I_1+I_2).
\end{align*}
 Using \eqref{e:3.4.1} and $r>C_0\Phi^{-1}(t)$, 
the proof of the upper bound of $I_1$  is the same as the one in \eqref{e:wq1}. Thus,  
we have
\begin{align*}
I_1 &\le c_7\left( \frac{ r}{\Phi^{-1}(t)} \right)^d\exp{\left(-\frac{a_2r^2}{2\Phi^{-1}(t)^2}\right)} \sum_{n=0}^\infty 2^{-nd} \le c_8\exp{\left(-\frac{a_2r^2}{4\Phi^{-1}(t)^2}\right)}.
\end{align*}
We next estimate $I_2$. Note that $\rho_n \ge \rho_0= C_0 \Phi^{-1}(t)^{1+\theta_0} \psi^{-1}(t)^{\theta_0} \ge C_0 \psi^{-1}(t)$. Thus,  we have 
\begin{align*}
I_2 & =\sum_{n=0}^\infty \left(\frac{2^nr}{\rho_n}\right)^d\frac{t}{\psi(\rho_n)} = \sum_{n=0}^\infty \left(\frac{2^{n(1-\alpha)} r^2}{C_0\Phi^{-1}(t)^2}\right)^d \frac{\psi(C_0\psi^{-1}(t))}{\psi(\rho_n)} \frac{\psi(\psi^{-1}(t))}{\psi(C_0\psi^{-1}(t))}\\
&\le c_9\sum_{n=0}^\infty 2^{n(d-\alpha(d+\beta_1))} \left( \frac{r}{\Phi^{-1}(t)} \right)^{2d} \left( \frac{r^2}{\Phi^{-1}(t)^2} \cdot \frac{\psi^{-1}(t)}{r} \right)^{\beta_1}  \\
&= c_{10} \left( \frac{\Phi^{-1}(t)}{r} \right)^{-2(d+\beta_1)} \left(\frac{\psi^{-1}(t)}{r}\right)^{\beta_1}.
\end{align*}
We now bound $\frac{\Phi^{-1}(t)}{r}$ in terms of $\frac{\psi^{-1}(t)}{r}$ carefully. 
Since $r= C_0\Phi^{-1}(t)^{1+\theta_0}/\psi^{-1}(t)^{\theta_0}$, ${C_0\psi^{-1}(t)}<C_0\Phi^{-1}(t) <r$,  and $\theta_0 \le \theta$, 
$$\frac{C_0\Phi^{-1}(t)}{r} = \left( \frac{\psi^{-1}(t)}{\Phi^{-1}(t)} \right)^{\theta_0} = \left( C_0\frac{\psi^{-1}(t)}{r} \right)^{\theta_0/(1+\theta_0)} \ge \left(C_0 \frac{\psi^{-1}(t)}{r} \right)^{\theta/(1+\theta)}.$$
By using $\theta=\frac{\beta_1}{4d+3\beta_1}$, we have
$$ \left( \frac{\Phi^{-1}(t)}{r} \right)^{-2(d+\beta_1)} \le c_{11}\left( \frac{\psi^{-1}(t)}{r} \right)^{-\frac{2(d+\beta_1)\theta}{1+\theta}} = c_{11}\left( \frac{\psi^{-1}(t)}{r} \right)^{-\beta_1/2}.$$
Thus, 
$$I_2 \le c_{10} \left( \frac{\Phi^{-1}(t)}{r} \right)^{-2(d+\beta_1)} \left(\frac{\psi^{-1}(t)}{r}\right)^{\beta_1} \le c_{10}c_{11}\left(\frac{\psi^{-1}(t)}{r}\right)^{\beta_1/2}.$$
Using estimates of $I_1$ and $I_2$, we arrive
\begin{equation}\label{e:p3}
\int_{B(x,r)^c} p(t,x,y)dy \le c_6(I_1 + I_2) \le c_{12}\bigg(\frac{\psi^{-1}(t)}{r}\bigg)^{\beta_1/2}+c_{12}\exp{\left(-\frac{a_2r^2}{4\Phi^{-1}(t)^2}\right)}. 
\end{equation}
Combining \eqref{tail1}, \eqref{tail2} and \eqref{e:p3}, we obtain \eqref{tail}.

Now, \eqref{e:exp-1}
 follows from \eqref{tail} and Lemma \ref{l:tail} with $f(r,t):=\big(r/\Phi^{-1}(t)\big)^{2}$.
\qed

Recall that, without loss of generality, whenever
$\Phi$ satisfies the weak lower scaling property at infinity with index $\delta>1$, we have {assumed} that $\Phi$ satisfies $L^1(\delta, \wt C_L)$ instead of $L^a(\delta, \wt C_L)$. We also recall our notations: 
\begin{align*}
\scK(s)=\sup_{b\le s}\frac{\Phi(b)}{b},\quad\wt \Phi(t)=\Phi(1)t^2\1_{\{0<t< 1\}}+\Phi(t)\1_{\{t\ge 1\}}, \quad  \scK_\infty(s)=\sup_{b\le s}\frac{\wt\Phi(b)}{b}.
\end{align*} 

We are now ready to prove the sharp upper bound of $p(t,x,y)$, which is the most delicate part of this paper. 
\begin{thm}\label{t:UHK-exp} 
\begin{enumerate}
\item[(1)] Assume that $\Phi$ satisfies $L_a(\delta, \wt C_L)$ with $\delta >1$. Then for any $T>0$, there exist constants $a_U>0$ and $c>0$ such that for every $x,y \in \R^d$ and $t< T$,
	\begin{equation}
	\label{e:UHK-exp1}
	p(t,x,y) 
	\leq \frac{c\, t}{|x-y|^d\psi(|x-y|)} + c\,\Phi^{-1}(t)^{-d} \exp{\left(-\frac{a_U|x-y|}{\scK^{-1}(t/|x-y|)}\right)}. 
	\end{equation}
Moreover, if $\Phi$ satisfies $L(\delta, \wt C_L)$, then \eqref{e:UHK-exp1} holds for all $t<\infty$.
\item[(2)] Assume that $\Phi$ satisfies $L^1(\delta, \wt C_L)$ with $\delta >1$. Then for any $T>0$, there exist constants $a_U'>0$ and $c'>0$ such that for every $x,y \in \R^d$ and $t \geq T$,
	\begin{equation*}
	p(t,x,y) 
	\leq \frac{c'\, t}{|x-y|^d\psi(|x-y|)} + c'\,\Phi^{-1}(t)^{-d} \exp{\left(-\frac{a_U'|x-y|}{\scK_{\infty}^{-1}(t/|x-y|)}\right)}. 
	\end{equation*}
\end{enumerate}
\end{thm}
\pf Take $\theta = \frac{\beta_1(\delta-1)}{2\delta d+\delta \beta_1+\beta_1}$ and $\wt C_0=\left(\frac{2C_1}{C_2\wt C_L^2}\right)^{1/(\delta-1)}$, where $C_1$ and $C_2$ are constants in Lemma \ref{l:5.2}. Without loss of generality, we may and  do assume that $\wt C_0\ge1$. Note that $\theta $ satisfies $\frac{\delta(d+\beta_1)}{\delta-1}\frac{ \theta}{1+\theta} = \frac{\beta_1}{2}$ and $\theta< \delta-1$. Let $\alpha \in (d/(d +\beta_1),1)$.

(1) Again we will show that there exist $a_1>0$ and $c_1>0$   such that for any $t\leq T$ and $r>0$,
\begin{equation}\label{e:tail-F}
\int_{B(x,r)^c} p(t, x,y)\,dy\le c_1\bigg(\frac{\psi^{-1}(t)}{r}\bigg)^{\beta_1/2}+c_1\exp{\left(-\frac{a_1r}{\scK^{-1}(t/r)}\right)}.
\end{equation}
When $r \le \wt C_0\Phi^{-1}(t)$ using \eqref{e:F3} we have for $t\leq T$
\begin{equation}
\label{e:tail-F1}
\int_{B(x,r)^c} p(t,x,y)dy \le 1\le  e^{c_2} \exp{\left(-\frac{\wt C_0\Phi^{-1}(t)}{\scK^{-1}(t/(\wt C_0\Phi^{-1}(t)))}\right)} \le  e^{c_2} \exp{\left(-\frac{r}{\scK^{-1}(t/r)}\right)}. 
\end{equation}
The proof of case $r > \wt C_0\frac{\Phi^{-1}(t)^{1+\theta}}{\psi^{-1}(t)^\theta}$ is exactly same as the corresponding part in the proof of \eqref{e:exp-1} in Theorem \ref{t:exp-1}. 

Now consider the case $\wt C_0\Phi^{-1}(t) < r \le \wt C_0\Phi^{-1}(t)^{1+\theta} / \psi^{-1}(t)^{\theta}$. In this case, there exists $\theta_0 \in (0,\theta]$ such that $r=\wt C_0 \Phi^{-1}(t)^{1+\theta_0}/\psi^{-1}(t)^{\theta_0}$.  Define $\rho=\scK^{-1}(t/r)$ and $\rho_n = \wt C_02^{n\alpha} \rho$ for integer $n \ge 0$. 
Note that for $t\leq T$ and $\wt C_0\Phi^{-1}(t)<r$, we have $t\le T\wedge \Phi(r)$. Thus, by \eqref{e:F4}
\begin{align}\label{rhobound}
\rho \leq \rho_0 =\wt C_0 \rho\le \wt C_0\frac{\Phi^{-1}(t)^2}{r}\leq\wt C_0\Phi^{-1}(T)\frac{\Phi^{-1}(t)}{r}\leq \Phi^{-1}(T).
\end{align}
By Remark \ref{mwsc}, we may assume that $\Phi^{-1}(T) <a$.
Thus, by \eqref{rhobound}, \eqref{e:F1}, the condition $L_a(\delta, \wt C_L)$ on $\Phi$, and the definition of $\wt C_0$,
we have
\begin{align}\label{e:exp232}
&C_1 \frac{t}{\Phi(\rho_n)} - C_2 \frac{2^n r}{\rho_n} \le C_1 \frac{\Phi(\rho)}{\Phi(\rho_0)}\frac{t}{\Phi(\rho)} - \frac{C_2}{\wt C_0} \frac{2^{n(1-\alpha)} r}{\rho} 
\nn\\& = \frac{C_2}{\wt C_0\rho} \left( \frac{\wt C_0C_1}{C_2}\frac{\Phi(\rho)}{\Phi(\rho_0)}\frac{\rho}{\Phi(\rho)}t -  2^{n(1-\alpha)} r \right)\le \frac{C_2}{\wt C_0\rho} \left( \frac{\wt C_0C_1}{C_2\wt C_L}\frac{\Phi(\rho)}{\Phi(\rho_0)}\frac{t}{\scK(\rho)} -  2^{n(1-\alpha)} r \right)
\nn\\
&\leq  \frac{C_2r}{\wt C_0\rho} \left( \frac{\wt C_0C_1}{C_2\wt C_L}\frac{\Phi(\rho)}{\Phi(\rho_0)} - 2^{n(1-\alpha)} \right) \leq  \frac{C_2r}{\wt C_0\rho} \left( \frac{\wt C_0^{1-\delta}C_1}{C_2\wt C_L^2} - 2^{n(1-\alpha)} \right) \nn\\
 &=  \frac{C_2r}{\wt C_0\rho} \left(\frac{1}{2} - 2^{n(1-\alpha)} \right)
\le -c_5 2^{n(1-\alpha)} \frac{r}{\rho}.
\end{align}
Combining \eqref{e:exp232}, Lemmas \ref{l:4.24} and \ref{l:5.2}, we have that 
\begin{align*}
p(t,x,y) &\le \frac{c_6t}{\rho_n^d \psi(\rho_n)} + c_6\Phi^{-1}(t)^{-d} \exp\left(C_1 \frac{t}{\Phi(\rho_n)} -C_2 \frac{2^n r}{\rho_n}\right) \nn \\ &\le \frac{c_6t}{\rho_n^d \psi(\rho_n)} + c_6\Phi^{-1}(t)^{-d} \exp \left( -c_52^{n(1-\alpha)} \frac{r}{\rho} \right). 
\end{align*}
With above estimate, we get that
\begin{align*}
\int_{B(x,r)^c} p(t,x,y)dy &\le \sum_{n=0}^\infty \int_{B(x,2^{n+1}r) \backslash B(x,2^n r)} p(t,x,y) dy \\
&\le  c_7\sum_{n=0}^{\infty}\left( \frac{2^n r}{\Phi^{-1}(t)} \right)^d \exp \left(-c_5 \frac{2^{n(1-\alpha)} r}{\rho}\right)
+ c_7\sum_{n=0}^\infty  \frac{t(2^nr)^d}{\rho_n^d \psi(\rho_n)}\\
&:= I_1+I_2.
\end{align*}
We first estimate $I_1$. Note that by \eqref{e:F4}, $r/\rho \geq    (r/\Phi^{-1}(t))^2 \ge \wt C_0^2$. Using this, \eqref{e:F4}, and \eqref{e:3.4.1} we have
\begin{align*}
I_1 &
\le c_8\sum_{n=0}^{\infty}\left( \frac{r}{\rho} \right)^{d/2} 2^{nd}\exp \left(-c_5 \frac{2^{n(1-\alpha)} r}{\rho}\right) \leq c_9 \left( \frac{r}{\rho} \right)^{d/2}\exp{\left(-\frac{c_5r}{2\rho}\right)} \sum_{n=0}^\infty 2^{-nd} \nn\\
&\leq c_{10} \exp{\left(-\frac{c_5r}{4\rho}\right)} \sum_{n=0}^\infty 2^{-nd}\le c_{11}\exp{\left(-\frac{c_5r}{4\rho}\right)}.
\end{align*}
We next estimate $I_2$. 
By using \eqref{e:F4}, $r = \wt C_0\Phi^{-1}(t)^{1+\theta_0} / \psi^{-1}(t)^{\theta_0}$, $\psi^{-1}(t)\le\Phi^{-1}(t)$, and  $\theta_0\le\theta<\delta-1$, we have
\begin{align*}
\frac{\Phi^{-1}(t)}{\rho}\le C_3\left(\frac{r}{\Phi^{-1}(t)}\right)^{1/(\delta-1)}=C_3\wt C_0^{1/(\delta-1)}\left(\frac{\Phi^{-1}(t)}{\psi^{-1}(t)}\right)^{\theta_0/(\delta-1)}\le C_3\wt C_0^{1/(\delta-1)}\frac{\Phi^{-1}(t)}{\psi^{-1}(t)}. 
\end{align*}
Thus, we have $\rho_n>\rho \ge C_3^{-1}\wt C_0^{-1/(\delta-1)}\psi^{-1}(t)$. 
Using this, $L(\beta_1, C_L)$ condition on $\psi$, and \eqref{e:F4},
\begin{align*}
I_2 & =c_7\sum_{n=0}^\infty \left(\frac{2^nr}{\rho_n}\right)^d \frac{t}{\psi(\rho_n)}
\leq c_{11}\sum_{n=0}^\infty \left(\frac{2^nr}{\rho_n}\right)^{d}\left(\frac{\psi^{-1}(t)}{\rho_n}\right)^{\beta_1}\\
& = c_{12}\sum_{n=0}^\infty 2^{n(d-\alpha(d+\beta_1))} \left( \frac{r}{\rho} \right)^{d} \left( \frac{r}{\rho} \cdot \frac{\psi^{-1}(t)}{r} \right)^{\beta_1} 
=c_{13} \left(\frac{r}{\rho}\right)^{d+\beta_1} \left(\frac{\psi^{-1}(t)}{r}\right)^{\beta_1} \\
& \le c_{14} \left( \frac{\Phi^{-1}(t)}{r} \right)^{-\frac{\delta}{\delta-1}(d+\beta_1)} \left(\frac{\psi^{-1}(t)}{r}\right)^{\beta_1}.
\end{align*}
Since $r=\wt C_0 \Phi^{-1}(t)^{1+\theta_0}/\psi^{-1}(t)^{\theta_0}$, $\psi^{-1}(t)\le\Phi^{-1}(t)<r$,  and $\theta_0 \le \theta$,
\begin{align*}
\frac{\Phi^{-1}(t)}{r} =\wt C_0^{-1} \left( \frac{\psi^{-1}(t)}{\Phi^{-1}(t)} \right)^{\theta_0} = \wt C_0^{-1/(1+\theta_0)}\left( \frac{\psi^{-1}(t)}{r} \right)^{\theta_0/(1+\theta_0)} \ge \wt C_0^{-1}\left( \frac{\psi^{-1}(t)}{r} \right)^{\theta/(1+\theta)}.
\end{align*}
Recall that $\theta$ satisfies
$$\frac{\delta(d+\beta_1)}{\delta-1}\frac{ \theta}{1+\theta} = \frac{\beta_1}{2}. $$
Thus, we have
$$ \left( \frac{\Phi^{-1}(t)}{r} \right)^{-\frac{\delta}{\delta-1}(d+\beta_1)} \le \wt C_0^{\delta(d+\beta_1)/(\delta-1)}\left( \frac{\psi^{-1}(t)}{r} \right)^{-\frac{\delta(d+\beta_1) }{\delta-1}\frac{\theta}{1+\theta}} = \wt C_0^{\delta(d+\beta_1)/(\delta-1)}\left( \frac{\psi^{-1}(t)}{r} \right)^{-\beta_1/2},$$
which implies
\begin{align*}
I_2 \le c_{14} \left( \frac{\Phi^{-1}(t)}{r} \right)^{-\frac{\delta}{\delta-1}(d+\beta_1)} \left(\frac{\psi^{-1}(t)}{r}\right)^{\beta_1} \le c_{15}\left(\frac{\psi^{-1}(t)}{r}\right)^{\beta_1/2}.
\end{align*}
Using estimates of $I_1$ and $I_2$, we obtain
\begin{equation}\label{e:tail-F3}
\int_{B(x,r)^c} p(t,x,y)dy \le I_1 + I_2 \le c_{15}\bigg(\frac{\psi^{-1}(t)}{r}\bigg)^{\beta_1/2}+c_{15}\exp{\left(-\frac{a_1r}{\scK^{-1}(t/r)}\right)}.
\end{equation} Combining \eqref{e:tail-F1}, \eqref{tail2} and \eqref{e:tail-F3} we obtain \eqref{e:tail-F}. 

Let $f(r,t):=\frac{r}{\scK^{-1}(t/r)}$. Then, by \eqref{rhobound} and  Lemma \ref{l:F}, we see that $f(r,t)$ satisfies the condition in Lemma \ref{l:tail}. Thus, by Lemma \ref{l:tail}, we obtain 
$$p(t,x,y) 
	\leq \frac{c_{16}\, t}{|x-y|^d\psi(|x-y|)} + c_{16}\,\Phi^{-1}(t)^{-d} \exp{\left(-\frac{c_{17}|x-y|}{\scK^{-1}(c_{18}t/|x-y|)}\right)}. $$
Since $ |x-y| \ge \wt C_0 \Phi^{-1} (t) \ge c_{19} t^{1/\delta} \ge t$, we can apply 
\eqref{e:wscF} and get $\scK^{-1}(c_{18}t/|x-y|) \le c_{20}\scK^{-1}(t/|x-y|)$. We have proved the first claim of the theorem. 

(2) The proof of the second claim is similar to the proof of the first claim. Let $f(r,t):=\frac{r}{\scK_{\infty}^{-1}(t/r)}$. Then, by \eqref{e:F_13} and \eqref{e:F_11}, we see that $f$ satisfies $f(c_0 r,t) \asymp f( r,t)$ and the assumption $(i)$  in Lemma \ref{l:tail} holds for $t>0$ and $r>0$.   Thus, it suffices to show that there exist $a_1, c_1>0$ such that for any $t\geq T$ and $r>0$,
\begin{equation}\label{e:tail_Finf}
\int_{B(x,r)^c} p(t, x,y)\,dy\le c_1\bigg(\frac{\psi^{-1}(t)}{r}\bigg)^{\beta_1/2}+c_1\exp{\left(-\frac{a_1r}{\scK_{\infty}^{-1}(t/r)}\right)}.
\end{equation}
 Since the proof of \eqref{e:tail_Finf} for the cases $r\le \wt C_0\Phi^{-1}(t)$ and  $r> \wt C_0\Phi^{-1}(t)^{1+\theta}/\psi^{-1}(t)^{\theta}$ are the same as that for (1), we only prove that \eqref{e:tail_Finf} holds for $\wt C_0\Phi^{-1}(t)<r\leq \wt C_0\Phi^{-1}(t)^{1+\theta}/\psi^{-1}(t)^{\theta}$. As (1), take $\theta_0 \in (0,\theta]$ such that $r=\wt C_0 \Phi^{-1}(t)^{1+\theta_0}/\psi^{-1}(t)^{\theta_0}$ and define $\rho=\scK_{\infty}^{-1}(t/r)$ and $\rho_n = \wt C_02^{n\alpha} \scK_{\infty}^{-1}(t/r)$ for integer $n \ge 0$.  
By Lemma \ref{PhiwtPhi}, we see that $\wt \Phi\leq \Phi$ and $\wt \Phi$ satisfies $L(\delta, \wt C_L)$. Using these, \eqref{e:F_11}, and the definition of $\wt C_0$,
we follow the argument in \eqref{e:exp232} and get
\begin{align}\label{e:exp2323}
C_1 \frac{t}{\Phi(\rho_n)} - C_2 \frac{2^n r}{\rho_n} \le C_1 \frac{t}{\wt \Phi(\rho_n)} - C_2 \frac{2^n r}{\rho_n}  
\le -c_2 2^{n(1-\alpha)} \frac{r}{\rho}.
\end{align}
Combining \eqref{e:exp2323} and Lemmas \ref{l:4.24} and \ref{l:5.2}, we get that for $\wt C_0\Phi^{-1}(t)<r\leq \wt C_0\Phi^{-1}(t)^{1+\theta}/\psi^{-1}(t)^{\theta}$,
\begin{align*}
&\int_{B(x,r)^c} p(t,x,y)dy \le \sum_{n=0}^\infty \int_{B(x,2^{n+1}r) \backslash B(x,2^n r)} p(t,x,y) dy \nn\\
&\le c_4\sum_{n=0}^\infty (2^nr)^d \Phi^{-1}(t)^{-d} \exp \left(-c_2 2^{n(1-\alpha)} \frac{r}{\rho}\right) + c_4\sum_{n=0}^\infty (2^nr)^d   \frac{t}{\rho_n^d \psi(\rho_n)}=: I_1+I_2.
\end{align*}
Since $r/\rho \geq  C_4^{-1}  (r/\Phi^{-1}(t))^2 \ge C_4^{-1}\wt C_0^2$ by \eqref{e:F_14}, the estimates of $I_1$ and $I_2$ are  similar to the arguments in (1). In fact, using \eqref{e:F_14} and \eqref{e:3.4.1}, we have
\begin{align*}
I_1 &\leq c_5\left(\frac{r}{\rho}\right)^{d/2} \exp{\left(-\frac{c_2r}{2\rho}\right)} \sum_{n=0}^\infty 2^{-nd}\leq c_6 \exp{\left(-\frac{c_2r}{4\rho}\right)} \sum_{n=0}^\infty 2^{-nd}
\le c_{7}\exp{\left(-\frac{c_5r}{4\rho}\right)}.
\end{align*}
Moreover, 
\begin{align*}
\frac{\Phi^{-1}(t)}{\rho}\le C_4\left(\frac{r}{\Phi^{-1}(t)}\right)^{1/(\delta-1)}=C_4\wt C_0^{1/(\delta-1)}\left(\frac{\Phi^{-1}(t)}{\psi^{-1}(t)}\right)^{\theta_0/(\delta-1)}\le C_4\wt C_0^{1/(\delta-1)}\frac{\Phi^{-1}(t)}{\psi^{-1}(t)},
\end{align*}
which implies $\rho_n>\rho \ge C_4^{-1}\wt C_0^{-1/(\delta-1)}\psi^{-1}(t)$. 
Using this, $L(\beta_1, C_L)$ condition on $\psi$, and \eqref{e:F_14}, by the same argument as in the proof of (1), we have
\begin{align*}
I_2 \le c_{8} \left( \frac{\Phi^{-1}(t)}{r} \right)^{-\frac{\delta}{\delta-1}(d+\beta_1)} \left(\frac{\psi^{-1}(t)}{r}\right)^{\beta_1} \le c_{9}\left(\frac{\psi^{-1}(t)}{r}\right)^{\beta_1/2}.
\end{align*}
Using estimates of $I_1$, $I_2$, and $\rho=\scK_{\infty}^{-1}(t/r)$,  we obtain \eqref{e:tail_Finf}.
\qed

Combining Theorems \ref{t:uhkd} and \ref{t:UHK-exp} and Lemma \ref{l:F}, we get the desired upper bounds of $p(t,x,y)$.
\begin{thm}
(1) Assume that $\Phi$ satisfies $L_a(\delta, \wt C_L)$ with $\delta >1$ and $T>0$. Then there exist constants $a_U>0$ and $c>0$ such that for every $x,y \in \R^d$ and $t\leq T$,
	\begin{equation}
	\label{e:UHK0}
	p(t,x,y) 
	\leq c\Phi^{-1}(t)^{-d}\wedge \left(\frac{c\, t}{|x-y|^d\psi(|x-y|)} + c\Phi^{-1}(t)^{-d} \exp{\left(-\frac{a_U|x-y|}{\scK^{-1}(t/|x-y|)}\right)}\right). 
	\end{equation}
Moreover, if $a=\infty$, then \eqref{e:UHK0} holds for all $t<\infty$.

(2) Assume that $\Phi$ satisfies $L^1(\delta, \wt C_L)$ with $\delta >1$ and $T>0$. Then there exist constants $a_U'$ and $c>0$ such that for every $x,y \in \R^d$ and $t \geq T$,
	\begin{equation*}
	p(t,x,y) 
	\leq c\Phi^{-1}(t)^{-d}\wedge \left(\frac{c\, t}{|x-y|^d\psi(|x-y|)} + c\Phi^{-1}(t)^{-d} \exp{\left(-\frac{a_U'|x-y|}{\scK_{\infty}^{-1}(t/|x-y|)}\right)}\right). 
	\end{equation*}
\end{thm}
\subsection{Off diagonal lower bound estimates}\label{s:ODLE}
\begin{prop}\label{l:LHK-psi}
	There exist constants $\delta_1 \in (0,1/2)$ and $C_5>0$ such that
\begin{align}\label{e:LHK-psi}
p(t,x,y) \ge  C_5\begin{cases} \Phi^{-1}(t)^{-d} &\text{ if } |x-y| \le \delta_1\Phi^{-1}(t) \\
 \displaystyle \frac{t}{|x-y|^d \psi(|x-y|)} &\text{ if } |x-y| \ge \delta_1\Phi^{-1}(t) .
\end{cases}
\end{align}
\end{prop}
\pf The proof of the proposition is standard. For example, see \cite[Proposition 5.4]{HKE}. 

Let $\delta_1=\eps/2<1/2$ where $\eps$ is the constant in Theorem  \ref{t:NDL}. Then 
by Theorem  \ref{t:NDL},
\begin{align}
\label{e:nbe1}
p(t,x,y) \ge  p^{B(x,\Phi^{-1}(t)/\eps)}(t,x,y) \ge c_0\Phi^{-1}(t)^{-d} \quad \text{ for all } |x-y| \le \delta_1\Phi^{-1}(t). 
\end{align}
Thus, we have  \eqref{e:LHK-psi} when $|x-y| \le \delta_1\Phi^{-1}(t)$.

  By Lemma \ref{l:4.2} we have
$$ \P^x(\tau_{B(x,r)} \le t) \le \frac{c_1 t}{\Phi(r)} $$
for any $r>0$ and $x \in \R^d$. 
Let $\delta_2:= (C_L/2)^{1/\beta_1}\delta_1\in (0,\delta_1)$ so that $\delta_1 \Phi^{-1}((1-b)t) \ge \delta_2 \Phi^{-1}(t) $ holds for all $b \in(0, 1/2]$.
Then choose 
$\lambda \le c_1^{-1}C_U^{-1} (2\delta_2/3)^{\beta_2}/2 <1/2$
small enough so that
 $c_1 \lambda t /\Phi(2\delta_2 \Phi^{-1}(t)/3)  \le \lambda c_1C_U (2\delta_2/3)^{- \beta_2}\le 1/2$. 
Thus we have  $\lambda \in (0,1/2)$ and $\delta_2 \in (0,\delta_1)$ (independent of $t$) such that
\begin{align}
\label{e:nbe3}
\delta_1 \Phi^{-1}((1-\lambda)t) \ge \delta_2 \Phi^{-1}(t),  \quad \text{for all } t>0
\end{align}
and
\begin{align}
\label{e:nbe4}
 \P^x(\tau_{B(x,2\delta_2 \Phi^{-1}(t)/3)} \le \lambda t) \le 1/2, \quad \text{for all } t>0 \text{ and } x \in \R^d.\end{align}
  For the remainder of the proof we assume that $|x-y| \ge \delta_1\Phi^{-1}(t)$.
Since, using \eqref{e:nbe1} and \eqref{e:nbe3}, 
\begin{align*}
p(t,x,y) &\ge \int_{B(y,\delta_1 \Phi^{-1}((1-\lambda)t))} p(\lambda t,x,z) p((1-\lambda)t,z,y)dy \\
&\ge \inf_{z \in B(y,\delta_1 \Phi^{-1}((1-\lambda)t))} p((1-\lambda)t,z,y) \int_{B(y,\delta_1 \Phi^{-1}((1-\lambda)t))} p(\lambda t,x,z)dz \\
&\ge c_0\Phi^{-1}(t)^{-d} \P^x ( X_{\lambda t} \in B(y,\delta_2 \Phi^{-1}(t))), 
\end{align*}
it suffices to prove 
\begin{align}
\label{e:nbe2}
\P^x ( X_{\lambda t} \in B(y,\delta_2 \Phi^{-1}(t)))  \ge c_2\frac{t\Phi^{-1}(t)^d}{|x-y|^d\psi(|x-y|)}.
\end{align}
For $A \subset \R^d$, let $\sigma_A := \inf\{t>0:X_t \in A \}$.
Using  \eqref{e:nbe4} and the strong Markov property we have 
\begin{align*}
& \P^x(X_{\lambda t} \in B(y,\delta_2 \Phi^{-1}(t))) \\
&\ge \P^x \left( \sigma_{B(y,\delta_2\Phi^{-1}(t)/3)} \le \lambda t; \sup_{s \in [\tau_{B(y,\delta_2\Phi^{-1}(t)/3)},\lambda t]} |X_s - X_{\tau_{B(y,\delta_2\Phi^{-1}(t)/3)}} | \le 2\delta_2\Phi^{-1}(t)/3 \right) \\
&\ge \P^x (\sigma_{B(y,\delta_2\Phi^{-1}(t)/3)} \le \lambda t ) \inf_{z \in B(y,\delta_2 \Phi^{-1}(t)/3)} \P^z (\tau_{B(z,2\delta_2 \Phi^{-1}(t)/3)} >\lambda t) \\
&\ge \frac{1}{2}\P^x (\sigma_{B(y,\delta_2\Phi^{-1}(t)/3)} \le \lambda t) \\
&\ge \frac{1}{2}\P^x \Big(X_{(\lambda t) \land \tau_{B(x,2\delta_2 \Phi^{-1}(t)/3)}} \in B(y,\delta_2\Phi^{-1}(t)/3) \Big).
\end{align*} 
Since $|x-y| \ge \delta_1 \Phi^{-1}(t) > \delta_2 \Phi^{-1}(t)$, clearly $B(y, \delta_2 \Phi^{-1}(t)/3) \subset \overline{B(x,2\delta_2\Phi^{-1}(t)/3)}^c$. Thus by \eqref{e:J_psi}, L\'evy system and \eqref{e:nbe4}, we have 
\begin{align*}
&\P^x \Big(X_{(\lambda t) \land \tau_{B(x,2\delta_2 \Phi^{-1}(t)/3)}} \in B(y,\delta_2\Phi^{-1}(t)/3) \Big) \\
& \quad = \E^x \left[  \sum_{s \le (\lambda t) \land \tau_{B(x,2\delta_2 \Phi^{-1}(t)/3)}} \1_{ \{ X_s \in B(y,\delta_2\Phi^{-1}(t)/3)  \}  }   \right]\\
& \quad \ge \E^x \Bigg[ \int_0^{(\lambda t) \land \tau_{B(x,2\delta_2 \Phi^{-1}(t)/3)}} ds \int_{B(y, \delta_2 \Phi^{-1}(t)/3)} J(X_s,u)du   \Bigg] \\
& \quad \ge c_3\E^x \Bigg[ \int_0^{(\lambda t) \land \tau_{B(x,2\delta_2 \Phi^{-1}(t)/3)}} ds \int_{B(y, \delta_2 \Phi^{-1}(t)/3)} \frac{1}{|X_s-u|^d\psi(|X_s-u|)}du   \Bigg] \\
& \quad \ge c_4\E^x [(\lambda t) \land \tau_{B(x,2\delta_2 \Phi^{-1}(t)/3)} ] \big(\delta_2 \Phi^{-1}(t)/3 \big)^d \frac{1}{|x-y|^d \psi(|x-y|)} \\
& \quad \ge c_4 (\lambda t) \P^x( \tau_{B(x,2\delta_2 \Phi^{-1}(t)/3)} \ge \lambda t )
\big(\delta_2 /3 \big)^d \frac{\Phi^{-1}(t)^d}{|x-y|^d \psi(|x-y|)}\\
& \quad \ge c_4 2^{-1} \lambda (\delta_2 /3)^d\frac{t\Phi^{-1}(t)^d}{|x-y|^d\psi(|x-y|)},
\end{align*}
where in the third inequality we used the fact that
$$|X_s -u | \le |X_s -x| + |x-y| + |y-u| \le |x-y| + \delta_2 \Phi^{-1}(t) \le 2|x-y|. $$
Thus, combining the above two inequality, we have proved \eqref{e:nbe2}.
\qed

We now give the two-sided sharp estimate for Green function.

\vspace{3mm}

\noindent
{\bf Proof of Corollary \ref{c:Green}.}
Let $\wt \delta :=\beta_2 \wedge 2<d$ and $r=|x-y|$.
By  Lemma \ref{l:rel1}, $\Phi$ satisfies $L(\beta_1,  C_L)$ and $U(\wt \delta,  C_U)$. 

For the lower bound, we use 
Proposition \ref{l:LHK-psi} and Lemma \ref{lem:inverse} and get
\begin{align*}
G(x,y)&\geq \int_{\Phi(r/\delta_1)}^\infty p(t,x,y)dt \geq c_1\int_{\Phi(r/\delta_1)}^{2\Phi(r/\delta_1)} \Phi^{-1}(t)^{-d}dt\geq c_2r^{-d}\Phi(r).
\end{align*}
For the upper bound, using the change of variables, the integration by parts and \eqref{e:Phi}, 
\begin{align*}
\int_{\Phi(r)}^\infty \Phi^{-1}(t)^{-d}dt &=\int_0^{r^{-1}}s^d d(-\Phi(s^{-1}))=-\big[s^d\Phi(s^{-1})\big]^{r^{-1}}_0 +d\int_0^{r^{-1}}s^{d-1}\Phi(s^{-1})ds\\
&\leq -\big[s^d\Phi(s^{-1})\big]^{r^{-1}}_0 + c_2\Phi(r)\int_0^{r^{-1}} s^{d-1}\left(\frac{s^{-1}}{r}\right)^{\wt \delta}ds.
\end{align*}
By using the condition $\wt \delta<d$, we get that
$
\int_{\Phi(r)}^\infty \Phi^{-1}(t)^{-d}dt \le c_3 \Phi(r)r^{-d}.
$
Using this inequality and  Theorem \ref{t:UHK-Phi}, we conclude that
\begin{align*}
G(x,y)&=\int_0^{\Phi(r)}p(t,x,y)dt +\int_{\Phi(r)}^\infty p(t,x,y)dt\\
&\leq \frac{c_3}{\Phi(r)r^d}
\int_0^{\Phi(r)} tdt+c_3\int_{\Phi(r)}^\infty \Phi^{-1}(t)^{-d}dt \le c_4r^{-d}\Phi(r).
\end{align*}
\qed

By using $\scK$ and $\scK_\infty$, we give the lower bound of $p(t,x,y)$ under $L_a(\delta,\wt C_L)$ or $L^1(\delta,\wt C_L)$ on $\Phi$ with $\delta>1$.  See \cite[Lemma 3.1--3.2]{Sz17} for similar bound for L\'evy processes.
\begin{prop}\label{l:LHK-exp1}
Suppose $\Phi$ satisfies $L_a(\delta, \wt C_L)$ with $\delta>1$ and for some $a>0$. For $T>0$ there exist  $C>0$ and $a_L>0$ such that for any $t\leq T$ and $x,y \in \R^d$,
\begin{align}\label{proplow}
p(t,x,y)\geq C\Phi^{-1}(t)^{-d}\exp \left(-a_L\frac{|x-y|}{\scK^{-1}(t/|x-y|)} \right).
\end{align}
Moreover, if $a=\infty$, then \eqref{proplow} holds for all $t<\infty$.
\end{prop}

\pf Let $r=|x-y|$.
By Proposition \ref{l:LHK-psi} and Remark \ref{mwsc}, without loss of generality, we assume that $\delta_1\Phi^{-1}(t) \leq r$ and $a\geq \delta_1\Phi^{-1}(T)$ where $\delta_1$ is the constants in Proposition \ref{l:LHK-psi}.
 Let $k=\ceil{{3r\delta_1^{-1}}/{\scK^{-1}(\frac{\delta_1t}{3r})}}$. Note that by \eqref{e:F4}, $$\scK^{-1}\left(\frac{t}{r}\right) \leq \frac{\Phi^{-1}(t)^2}{r}\leq \delta_1\Phi^{-1}(t)\leq \delta_1\Phi^{-1}(T)\leq a.$$ Thus by \eqref{e:wscF} we have
$
\scK^{-1}\left(\frac{t}{r}\right) \leq \wt C_L^{-1}(\frac{3}{\delta_1})\scK^{-1}(\frac{\delta_1t}{3r}).
$
Since
$
\frac{\delta_1t}{3r} \leq \frac{\delta_1\Phi(r/\delta_1)}{3r}\leq \frac{1}{3}\scK(\frac{r}{\delta_1})$, we see that $\scK^{-1}(\frac{\delta_1t}{3r})\leq \frac{r}{\delta_1} $, hence
\begin{align}\label{ku}
3\leq k\leq \frac{4r}{\delta_1\scK^{-1}(\frac{\delta_1t}{3r})}\leq \frac{12 \wt C_L^{-1}r}{\delta_1^2\scK^{-1}(\frac{t}{r})}.
\end{align}
By \eqref{e:F1} and our choice of $k$ 
$$
\Phi\left(\frac{3r}{\delta_1k}\right)\frac{\delta_1k}{r}\leq3\scK\left(\frac{3r}{\delta_1k}\right)\leq \frac{\delta_1t}{r}.
$$
Thus, we have 
$
\frac{r}{k}\leq \frac{\delta_1}{3}\Phi^{-1}(t/k).
$
Let $z_l=x+\frac{l}{k}(y-x)$, $l=0,1,\cdots, k-1$. For $\xi_l \in B(z_l, \frac{\delta_1}{3}\Phi^{-1}(\frac{t}{k}))$ and $\xi_{l-1} \in B(z_{l-1}, \frac{\delta_1}{3}\Phi^{-1}(\frac{t}{k}))$, 
$|\xi_l-\xi_{l-1}|\leq |\xi_l-z_l| + |z_l-z_{l-1}|+|z_{l-1}-\xi_{l-1}|\leq \delta_1\Phi^{-1}(t/k).$
Thus by Proposition \ref{l:LHK-psi}, $p(\frac{t}{k},\xi_{l-1},\xi_l) \geq C_5\Phi^{-1}(t/k)^{-d}$. Using the semigroup property   and \eqref{ku}, we get
\begin{align}
&p(t,x,y) \geq \int_{B(z_{k-1},\frac{\delta_1}{3}\Phi^{-1}(t/k))} \cdots\int_{B(z_{1},\frac{\delta_1}{3}\Phi^{-1}(t/k))}p(\tfrac{t}{k},x,\xi_1)\cdots p(\tfrac{t}{k},\xi_{k-1},y)d\xi_1\cdots d\xi_{k-1}\nn \\
&\geq C_5^k \Phi^{-1}(t/k)^{-dk} \prod_{l=1}^{k-1}\Big|B(z_l,\tfrac{\delta_1}{3}\Phi^{-1}(t/k))\Big| = c_2c_3^k\Phi^{-1}(t/k)^{-dk}\left(\frac{\delta_1}{3}\Phi^{-1}(t/k)\right)^{d(k-1)}\nn \\
&\geq c_2\left(\frac{c_3\delta_1^d}{3^d}\right)^k\Phi^{-1}(t)^{-d}
\geq c_2\Phi^{-1}(t)^{-d}\exp\left({-C_6 k}\right)\label{3} \\
&\geq c_2\Phi^{-1}(t)^{-d}\exp\left({-c_4\frac{r}{\scK^{-1}(t/r)}}\right). \nn
\end{align}
This finishes the proof. Here we record that the constant $C_6$ in \eqref{3} depends only on $d$ and constants $\delta_1, C_5$ in \eqref{e:LHK-psi}. \qed

Recall that $\beta_2$ is the upper scaling index of $\psi$.
\begin{prop}\label{l:LHK-exp1n}
Suppose $\Phi$ satisfies $L^1(\delta,\wt C_L)$ with $\delta>1$.  For any $T>0$ and $\theta>0$ satisfying $\frac{1}{\delta}+\theta(\frac{1}{\delta}-\frac{1}{\beta_2})\leq 1$, there exist  $c_1$, $c_2>0$ and $a_L'>0$ such that for 
$(t,x,y)\in[T,\infty)\times\bR^d\times\bR^d$ satisfying
$\delta_1\Phi^{-1}(t) < |x-y| \le c_1\Phi^{-1}(t)^{1+\theta} / \psi^{-1}(t)^{\theta}$,
$$
p(t,x,y)\geq c_2\Phi^{-1}(t)^{-d}\exp \left(-a_L'\frac{|x-y|}{\scK_{\infty}^{-1}(t/|x-y|)} \right),
$$
where $\delta_1$ is the constant in Proposition \ref{l:LHK-psi}.
\end{prop}
\pf
Without loss of generality, we assume that $T\ge\Phi(1)$. 
Take $c_1>0$ small so that
$$c_1(\Phi^{-1}(T)\wt C_L^{-1/\delta}T^{-{1}/{\delta}})^{1+\theta}(\psi^{-1}(T)^{-1}C_U^{1/\beta_2}T^{1/\beta_2})^{\theta}(1\vee T^{-1})\le \frac{\delta_1}{3\scK_{\infty}(2)}.$$
Since  $\psi$ satisfies $U(\beta_2, C_U)$ and $\Phi$ satisfies $L^1(\delta, \wt C_L)$, we see that  for $t\geq T\geq \Phi(1)$, $\psi^{-1}(t)\geq \psi^{-1}(T)C_U^{-1/\beta_2}T^{-1/\beta_2}t^{1/\beta_2}$ and $\Phi^{-1}(t)\leq \Phi^{-1}(T)\wt C_L^{-1/\delta}T^{-{1}/{\delta}}t^{1/\delta}$ by Lemma \ref{lem:inverse}. Thus, we have
\begin{align}\label{lowerbound} 
|x-y| &\le c_1\Phi^{-1}(t)^{1+\theta} / \psi^{-1}(t)^{\theta}
\leq c_1(\Phi^{-1}(T)\wt C_L^{-1/\delta}T^{-{1}/{\delta}})^{1+\theta}(\psi^{-1}(T)^{-1}C_U^{1/\beta_2}T^{1/\beta_2})^{\theta}t^{\frac{1}{\delta}+\frac{\theta}{\delta}-\frac{\theta}{\beta_2}}\nn\\
&\leq c_1(\Phi^{-1}(T)\wt C_L^{-1/\delta}T^{-{1}/{\delta}})^{1+\theta}(\psi^{-1}(T)^{-1}C_U^{1/\beta_2}T^{1/\beta_2})^{\theta}(1\vee T^{-1})t
\le \frac{\delta_1 t}{3\scK_{\infty}(2)},
 \end{align}
where the third inequality follows from $\frac{1}{\delta}+\theta(\frac{1}{\delta}-\frac{1}{\beta_2})\leq 1$.
Let $r =|x-y|$ and $k=\ceil{{3r}\delta_1^{-1}/{\scK_{\infty}^{-1}(\frac{\delta_1t}{3r})}}$. 
Since $r\ge \delta_1\Phi^{-1}(t)\ge\delta_1$, we have by \eqref{e:F_11} that
${\delta_1t}/{r} \leq {\delta_1\Phi(r/\delta_1)}/{r}={\delta_1\wt\Phi(r/\delta_1)}/{r}\leq  \scK_{\infty}(r/\delta_1).$ 
Thus, 
$
\scK_{\infty}^{-1}({\delta_1t}/{3r})
\le \scK_{\infty}^{-1}(\frac13 \scK_{\infty}({r}/{\delta_1}))
\le \frac{r}{\delta_1},
$
which implies that 
\begin{align}\label{ku2}
3\le k\le \frac{4r}{\delta_1\scK_{\infty}^{-1}(\frac{\delta_1t}{3r})}
\le \frac{12\wt C_L^{-1}r}{\delta_1^2\scK_{\infty}^{-1}(\frac{t}{r})}.
\end{align}
On the other hand, since
$\scK_{\infty}^{-1}(\frac{\delta_1t}{3r})\ge \scK_{\infty}^{-1}(\scK_\infty(2))=2$ and $\frac{3r}{\delta_1}\ge 3\Phi^{-1}(T)\ge 3$,
we see that
$$\frac{3r}{\delta_1\scK_{\infty}^{-1}(\frac{\delta_1t}{3r})}\le k<\frac{3r}{\delta_1}.$$
Thus, by the above inequality and \eqref{e:F_11}, we get
$$
\Phi\left(\frac{3r}{\delta_1k}\right)\frac{\delta_1k}{3r}=\wt \Phi\left(\frac{3r}{\delta_1k}\right)\frac{\delta_1k}{3r} \leq \scK_{\infty}\left(\frac{3r}{\delta_1k}\right)\leq \frac{\delta_1t}{3r},
$$
which yields
$
\frac{r}{k}\leq \frac{\delta_1}{3}\Phi^{-1}(t/k).
$
Using this, Proposition \ref{l:LHK-psi}, the semigroup property and  \eqref{ku2}, the remaining part of the proof is same as the one in the proof of Proposition \ref{l:LHK-exp1}.
 \qed

\begin{thm}\label{c:UHK-scinf}
(1) Suppose $\Phi$ satisfies $L_a(\delta, \wt C_L)$ with $\delta>1$ and for some $a>0$. For any $T>0$, there exist  $C>0$ and $a_L>0$ such that for any $(t, x, y)\in(0, T)\times\bR^d\times\bR^d$,
\begin{align}\label{lint}
p(t,x,y)\geq C\left(\Phi^{-1}(t)^{-d}\wedge\left( \frac{t}{|x-y|^d \psi(|x-y|)} +\Phi^{-1}(t)^{-d}\exp \left(-\tfrac{a_L|x-y|}{\scK^{-1}(t/|x-y|)} \right)\right)\right).
\end{align}
Moreover, if $a=\infty$, then \eqref{lint} holds for all $t<\infty$.

(2) Suppose $\Phi$ satisfies $L^1(\delta, \wt C_L)$ with $\delta>1$. For any $T>0$, there exist constants $C>0$ and $a_L'$ such that for $(t, x, y)\in[T, \infty)\times\bR^d\times\bR^d$,
$$
	p(t,x,y) 
	\geq  C\left(\Phi^{-1}(t)^{-d}\wedge\left( \frac{ t}{|x-y|^d\psi(|x-y|)} + \Phi^{-1}(t)^{-d} \exp{\left(-\tfrac{a_L'|x-y|}{\scK_{\infty}^{-1}(t/|x-y|)}\right)}\right)\right). 
	$$
In particular, if $\delta=2$, then  $ \scK_{\infty}^{-1}(s)\asymp s$.
\end{thm}
\pf
(1) is a direct consequence of Proposition \ref{l:LHK-psi} and  Proposition \ref{l:LHK-exp1}.

\noindent(2) Let $r=|x-y|$. 
By Proposition \ref{l:LHK-psi}, it is enough to show that for $t\ge T$ and $r>\delta_1\Phi^{-1}(t)$,
\begin{align}\label{pf:lower0}
p(t,x,y)\geq C\Phi^{-1}(t)^{-d}\exp \left(-a_L'\frac{r}{\scK_{\infty}^{-1}(t/r)} \right).
\end{align}
By  Proposition \ref{l:LHK-exp1n}, it suffices to show \eqref{pf:lower0} when $(c_1\Phi^{-1}(t)^{1+\theta} / \psi^{-1}(t)^{\theta})\vee \delta_1\Phi^{-1}(t)<r$, where $c_1$ and $\theta$ are the constants in Proposition \ref{l:LHK-exp1n}.

By Proposition \ref{l:LHK-psi}, we have
$p(t,x,y)\geq \frac{c_2t}{r^d \psi(r)}.$ 
By the same argument as in \eqref{jgeqexp}, we have that for some $N\ge1$,
$$
 \frac{t}{r^d \psi(r)} \geq c_3\Phi^{-1}(t)^{-d} \exp{\left(-\frac{a_2r^{1/N}}{\Phi^{-1}(t)^{1/N}}\right)}.
$$
Since $r\geq \delta_1\Phi^{-1}(t)$ and $T\le t$, by \eqref{e:F_14}, we get
\begin{align*}
\left(\frac{\delta_1^{-1}r}{\Phi^{-1}(t)}\right)^{1/N}
\leq \left(\frac{\delta_1^{-1}r}{\Phi^{-1}(t)}\right)^{2}
\leq  C_4\frac{\delta_1^{-1}r}{\scK_{\infty}^{-1}(\delta_1t/r)}.
\end{align*}
By \eqref{e:wscF_1} and Lemma \ref{lem:inverse},  $\scK_{\infty}$ satisfies $L(\delta-1, \wt C_L^{-2/(\delta-1)})$, which yields
\begin{align*}
\frac{\delta_1^{-1}r}{\scK_{\infty}^{-1}(\delta_1t/r)}
\leq \frac{\wt C_L^{-2/(\delta-1)}\delta_1^{-\delta/(\delta-1)}r}{\scK_{\infty}^{-1}(t/r)}.
\end{align*}
Thus,  
$
p(t,x,y)\geq c_2c_3\Phi^{-1}(t)^{-d} \exp{\left(-\frac{a_3r}{\scK_{\infty}^{-1}(t/r)}\right)}.
$
\qed

\begin{remark}\label{r:replace}
\rm
Suppose $\Phi$ satisfies $L^1(\delta, \wt C_L)$. Let $\theta=\frac{\beta_1(\delta-1)}{2\delta d+\delta \beta_1+\beta_1}$ which is defined in the proof of Theorem \ref{t:UHK-exp}. Then, $\theta$ satisfies $\frac{1}{\delta}+\theta(\frac{1}{\delta}-\frac{1}{\beta_2})\leq 1$. Thus, for any $T>0$ and $C>0$,  if $(t, x, y)\in[T,\infty)\times\bR^d\times\bR^d$ satisfies $C\Phi^{-1}(t)<|x-y|\le C\Phi^{-1}(t)\left(\frac{\Phi^{-1}(t)}{\psi^{-1}(t)}\right)^{\theta}$, then 
\begin{align*}
c_0^{-1}\Phi^{-1}(t)^{-d}\exp \left(-\frac{a_L'|x-y|}{\scK_{\infty}^{-1}(t/|x-y|)} \right)
\le p(t,x,y)\le c_0\Phi^{-1}(t)^{-d}\exp \left(-\frac{a_U'|x-y|}{\scK_{\infty}^{-1}(t/|x-y|)} \right).
\end{align*}
In this case, as one can see from the estimates  in \eqref{lowerbound}, there exists $c=c(C, T, \delta, \beta_1, \beta_2, \wt C_L, C_U)$ such that $\frac{t}{|x-y|}\ge c$. Thus, we only need $\scK_{\infty}(s)$ for $s\ge c_0:=\scK_{\infty}^{-1}(c)$ to estimate $p(t ,x,y)$. 

On the other hand, by \eqref{e:FcompFinf}, there exists $c_1=c_1(c_0)>0$ such that for $s\ge c_0$,
$$c_1^{-1}\sup_{c_0\le b\le s}\frac{\Phi(b)}{b}\le \scK_{\infty}(s)\le c_1\sup_{c_0\le b\le s}\frac{\Phi(b)}{b}.$$
Thus, in case of $t\ge T$ and $C\Phi^{-1}(t)<|x-y|\le C\Phi^{-1}(t)\left(\frac{\Phi^{-1}(t)}{\psi^{-1}(t)}\right)^{\theta}$, we may replace $\scK_{\infty}^{-1}$ with the generalized inverse of the function $f(s)=\sup_{c_0\le b\le s}\tfrac{\Phi(b)}{b}$.
\end{remark}

\noindent{\bf Proof of Corollary \ref{c:main2}.}
Since the upper bound is a direct consequence of Theorem \ref{t:exp-1}, we show the lower bound in \eqref{sbm}. 
Let $r=|x-y|$ and $\phi(s):=\Phi(s^{-1/2})^{-1}$. Since $\Phi$ satisfies $L_a(\delta, \wt C_L)$, $\phi$ satisfies $L^{1/a^{2}}(\delta/2, \wt C_L)$. Let $Z$ be a subordinate Brownian motion whose Laplace exponent is $\phi$. Then, by \cite[Proposition 3.5]{M}, for any $T>0$, there exist positive constants $\wt a_L$ and $c_1$ such that for all $(t,x,y)\in(0, T)\times\bR^d\times\bR^d$, 
\begin{align*}
&c_1\,\Phi^{-1}(t)^{-d/2}\exp\left(-\frac{\wt a_Lr^2}{\Phi^{-1}(t)^2}\right)
\le p^{Z}(t,x,y).
\end{align*}
On the other hand, by Theorem \ref{t:main}, there exist positive constants $a_U$ and $c_2$ such that
\begin{align*}
p^{Z}(t,x,y)
\le c_2\,\Phi^{-1}(t)^{-d/2}\exp\left(-\frac{a_Ur}{\scK^{-1}(t/r)}\right)+ \frac{c_2t}{r^d\psi(r)}.
\end{align*}
Let $a_L\ge a_U$ be a constant in Theorem \ref{t:main} and $A:=a_L/a_U\ge1$. Then,
$
{\scK^{-1}(t/Ar)}\le {\scK^{-1}(t/r)},
$
which implies that for all $(t,x,y)\in(0, T)\times\bR^d\times\bR^d$,
\begin{align*}
c_1\,\Phi^{-1}(t)^{-d/2}\exp\left(-\frac{\wt a_L A^2r^2}{\Phi^{-1}(t)^2}\right)
&\le c_2\,\Phi^{-1}(t)^{-d/2}\exp\left(-\frac{a_U Ar}{\scK^{-1}(t/Ar)}\right)+ \frac{c_2t}{(Ar)^d\psi(Ar)}\nn\\
&\le c_2\,\Phi^{-1}(t)^{-d/2}\exp\left(-\frac{a_Lr}{\scK^{-1}(t/r)}\right)+ \frac{c_3t}{r^d\psi(r)}.
\end{align*}
Thus,  by Theorem \ref{t:main}, we obtain the desired results.
\qed

\section{Application to the Khintchine-type law of iterated logarithm}\label{s:App}

In this section, we apply our main results in previous sections and 
 show that,  if our symmetric jump process has  the  finite second moment,  the Khintchine-type law of iterated logarithm at the infinity holds. Furthermore, we will also prove the converse.

We first establish the zero-one law for tail events.

\begin{thm}\label{t:01tail}
Let $A$ be a tail event. Then, either $\P^x(A)=0$ for all $x$ or else $\P^x(A)=1$ for all $x \in \R^d$.
\end{thm}
\pf
Fix $t_0,\eps>0$ and $x_0 \in \R^d$. 
Note that, by Lemma \ref{l:4.2}, there exists 
$c_1>0$ such that 
\begin{equation}\label{eq:enwq1}
\P^{x_0}(\sup_{s\leq t_0}|X_s-x_0|>
c_1\Phi^{-1}(t_0))<\eps.
\end{equation}
While, using Theorem \ref{t:PHR} to $P_{t}f$, 
we can  choose $t_1>0$
large so that for all $f \in L^\infty (\R^d)$ and $x\in \R^d$ with  $|x-x_0|\leq
c_1\Phi^{-1}(t_0)$, 
\begin{align}\label{eq:copt}
|P_{t_1}f(x)-P_{t_1}f(x_0)| &\le 
 c_2 \left(\frac{|x_0-x|}{\Phi^{-1}(t_1)}\right)^\theta \sup_{t>0}
\Vert P_{t}f\Vert_\infty \nn\\
&\le c_2 \left(\frac{c_1\Phi^{-1}(t_0)}{\Phi^{-1}(t_1)}\right)^\theta \sup_{t>0}
\Vert f\Vert_\infty<\eps\Vert f\Vert_\infty.
\end{align}
Note that \eqref{eq:enwq1} and \eqref{eq:copt} are same as \cite[(A.6) and (A.7)]{KKW}, and the proof of the theorem is exactly same as that of \cite[Theorem 2.10]
{KKW}. 
\qed

From \eqref{e:intcon} and \eqref{e:J_psi}, we see that the following three conditions are equivalent: 
\begin{align}
\label{e:fsecond1}
\sup_{x \in \R^d} \Big(\text{or } \inf_{x \in \R^d} \Big)\int_{\R^d} J(x,y) |x-y|^2dy < \infty ;
\end{align}
\begin{align}
\label{e:fsecond2}
c^{-1} r^2 \le \Phi(r) \le c r^2, \quad r>1;
\end{align}
\begin{align}
\label{e:fsecondp}
\int_0^\infty \frac{sds}{\psi(s)}< \infty.
\end{align} 
Using 
Theorem \ref{t:exp-1}, we see that under the assumption  \eqref{e:intcon},
the above conditions  \eqref{e:fsecond1}--\eqref{e:fsecondp}
are also equivalent to  the following weak and strong finite second moment conditions:
\begin{align}
\label{e:fsecond4}
\sup_{x \in \R^d}  \E^x[|X_t-x|^2]  < \infty 
\quad \text{ for all } \,\, t >0;
\end{align}
\begin{align}
\label{e:fsecond41}
 \inf_{x \in \R^d} \E^x[|X_t-x|^2]  < \infty 
\quad \text{ for some } \,\, t >0.
\end{align}
Here is the main result of this section. 
\begin{thm}\label{t:litrtd} Suppose $X$ is symmetric pure-jump process whose jumping density $J$ satisfies 
\eqref{e:J_psi}.
\noindent
(1) If   \eqref{e:fsecond2} holds, then 
there exists a constant $c\in(0,\infty)$ such that  for
all $x\in \R^d$, 
\begin{align*}
\limsup_{t\to \infty} \frac{ |X_t-x|}{(t\log\log t)^{1/2}}=c \quad \mbox{for} \quad \P^x\mbox{-a.e.}
\end{align*}

\noindent
(2) Suppose that  \eqref{e:intcon} holds but  \eqref{e:fsecond2} does not hold, i.e., \begin{align}
\label{e:fsecond7}
\int_0^\infty \frac{sds}{\psi(s)}=\infty.
\end{align} 
Then
for
all $x\in \R^d$,
\begin{align*}
\limsup_{t\to \infty} \frac{ |X_t-x|}{(t \log\log t)^{1/2}}=\infty \quad \mbox{for} \quad \P^x\mbox{-a.e.}
\end{align*}

\end{thm}

\pf  Without loss of generality we assume that $\Phi(1)=1$. Let $h(t)=t^{1/2}(\log\log t)^{1/2}$. We first observe that, by applying the change of variable $s=h(t)$,  
\begin{align}\label{e:cvi}
\int_{2 (\log\log 4)^{1/2}}^\infty \frac{sds}{\psi(s)} =\frac{1}{2}\int_4   ^\infty \frac{(\log\log t)+ (\log t)^{-1}}{\psi(h(t))}dt \asymp 
\int_4   ^\infty \frac{\log\log t}{\psi(h(t))}dt,
\end{align}
and 
\begin{align}\label{e:cvi1}
\int_{2 (\log\log 4)^{1/2}}^\infty \frac{s ds}{\psi(s) \log\log s} =\frac{1}{2}\int_4   ^\infty \frac{(\log\log t)+ (\log t)^{-1}}{\psi(h(t))
\log\log   [ t^{1/2}(\log\log t)^{1/2} ]}dt 
\asymp
\int_4   ^\infty \frac{dt}{\psi(h(t))}.
\end{align}

(1) 
By \eqref{e:fsecondp} and \eqref{e:cvi},
\begin{align}\label{e:uplog}
&\sum_{k=3}^\infty 
\frac{2^k}{\psi(2^{k/2}(\log\log 2^k)^{1/2})} \le c_0 \sum_{k=2}^\infty \int_{2^k}^{2^{k+1}} \frac{d t}{\psi(h(t))}= \int_4^\infty \frac{d t}{\psi(h(t))} < \infty.
\end{align}
Since we have $L^1(2, C_L)$ under \eqref{e:fsecond2},  by Theorem \ref{t:main}(2)  we have that for all $C>0$, $t>4$ and 
$t\le u \le 4t$, 
\begin{align*}
 &\P^x(|X_u-x|>C h(t)) =
\int_{|x-y| >C h(t)} p(u,x,y)dy \\
&\le c_1 \left(t\int_{|x-y| >C h(t)}
 J(x,y) dy+ t^{-d/2}\int_{|x-y| >C h(t)} \exp\left(-c_2\tfrac{|x-y|^2}{t} \right) dy \right)\\
 &\le c_3 \left(t\int_{C h(t)}^{\infty} \frac1{s\psi(s)}  ds+ t^{-d/2}\int_{C h(t)}^{\infty} \exp\left(-c_2\tfrac{s^2}{t} \right) s^{d-1} ds \right)=:c_3(I+II).
 \end{align*}
 First, letting $s_1 = \frac{s^2}{t}$ and we obtain
 \begin{align*}
II \le c_4\int_{C^2 \log\log t}^\infty e^{-c_2s_1}s_1^{d/2-1}ds_1 \le c_5 \int_{C^2 \log \log t}^\infty e^{-c_2s/2}ds \le \frac{c_5}{c_2} (\log t)^{-C^2c_2/2}.
  \end{align*}
  Let $C:=1\vee 2 c_2^{-1/2}$ so that  $II \le c_5  (\log t)^{-2}.$
While, by Lemma \ref{l:int_outball}, 
    \begin{align*}
I  \le c_6 \frac{t}{\psi(Ch(t))} \le c_7 \frac{t}{\psi(h(t))}. 
    \end{align*}
    Thus, for every $t>4$ and 
$t\le u \le 4t$, 
 \begin{align*}   
  \P^x(|X_u-x|>C h(t)) \le   c_8 \left( \frac1{(\log t)^2} + \frac{t}{\psi(h(t))}  \right).
    \end{align*}
    Using this and the strong Markov property,
    with $t_k=2^k$, $k \ge 3$ we get 
     \begin{align*}   
 & \P^x(|X_s-x|>2C h(s) \text{ for some } s \in [t_{k-1}, t_{k}])
  \le  \P^x(
  \tau_{B(x, C h(t_{k-1}))} \le t_k)
\\
  & 
 \le 2  \sup_{s \le t_k, z \in \R^d}  \P^z(|X_{t_{k+1}-s}-z|>C h(t_{k-1}))
\le   c_9 \left( \frac1{k^2} + \frac{2^k}{\psi(2^{k/2}(\log\log 2^k)^{1/2})} \right),
    \end{align*}
   where we followed the calculations in \eqref{cac1}. Therefore, by \eqref{e:uplog}
  \begin{align*}   
 &\sum_{k=3}^\infty \P^x(|X_s-x|>2C h(s) \text{ for some } s \in [t_k, t_{k+1}])<\infty. 
 \end{align*} 
By the Borel-Cantelli lemma, the above implies that 
 $$
  \P^x(|X_t-x| \le 2C h(t)\text{ for all sufficiently large } t)=1.
 $$
 Thus,
 $$
 \limsup_{t \to \infty}\frac{|X_t-x|}{h(t)} \le 2C.
 $$
 Since $\psi(r) \ge  \Phi(r) $,  
 by \eqref{e:fsecond2}, $J(x,y) \le c_9 |x-y|^{-d-2}$ for $|x-y|>1$. Thus, by
 \cite[Theorem 1.2(2)]{SW} (which can be proved using Theorem \ref{t:main}(2) and the second Borel-Cantelli lemma), we have that there exists $c_{10}>0$ such that 
  \begin{align*} 
&  \P^x(|X_t-x| > c_{10} h(t) \text{ for infinitely many } t)\\&=1-  \P^x(|X_t-x| \le c_{10} h(t) \text{ for all sufficiently large }t)=1.
 \end{align*} 
 Therefore,
$$
\P^x(c_{10} \le \limsup_{t \to \infty}\frac{|X_t-x|}{h(t)} \le 2C)=1.
 $$
Now using the zero-one law in  Theorem \ref{t:01tail}, we conclude that  there exists $c \in [c_{10}, 2C]$ such that
$$\P^x( \limsup_{t \to \infty}\frac{|X_t-x|}{h(t)} = c_{11} )=1, \quad \text{for all } x \in \R^d.
$$

      (2)  
 Using Theorem \ref{t:UHK-Phi}, there is $\lambda \in (0,1)$ such that 
   $$
 \sup_{t \ge 1}  \sup_{y \in \R^d} \int_{|z-y| <\lambda \Phi^{-1}(t)} p(t, y,z)dz
 \le c_0  \sup_{t \ge 1} |\lambda \Phi^{-1}(t)|^d (\Phi^{-1}(t))^{-d}=c_0\lambda^d
  <1/2.
   $$
   Let $t_k=2^k$. By the strong Markov property, we have that for all $C>0$ 
   \begin{align*}  
  \P^x( |X_{t_{k+1}}-X_{t_{k}}|\ge C h(t_{k+1})\, |\, \FF_{t_k} )&\ge \inf_{y \in \R^d}   \P^y( |X_{t_{k}}-y|\ge C h(t_{k+1}))\\
&=\inf_{y \in \R^d} \int_{|z-y| \ge C h(t_{k+1})} p(t_{k}, y,z)dz. 
  \end{align*}
We claim that for every $C>1$, 
\begin{equation}\label{sup-2-3}\sum_{k=1}^\infty \P^x ( |X_{t_{k+1}}-X_{t_{k}}|\ge C h(t_{k+1})|\mathscr{F}_{t_{k}})=\infty,\end{equation}
which implies the theorem.
In fact, by the second Borel-Cantelli lemma, $\P^x(\limsup \{|X_{t_{k+1}}-X_{t_{k}}|\ge C h(t_{k+1})\})=1$. Whence, for infinitely many $k\ge 1$,
  $|X_{t_{k+1}}-x|\ge C h(t_{k+1})/2  \\ $
  or $|X_{t_{k}}- x|\ge C h(t_{k+1})/2  \ge C h(t_{k})/2.$
Therefore, for all
$x\in \R^d$
 $$\limsup_{t\to \infty} \frac{ |X_t-x| }{h\big({t}\big)}=\limsup_{k\to \infty} \frac{|X_{t_k}-x|}{h(t_k)}\ge \frac{C}{2}, \quad \P^x\mbox{-a.e. }$$
Since the above holds for every $C>1$, the theorem follows. 

We now prove the claim \eqref{sup-2-3} by considering two cases separately. 

      \noindent
   {\it Case 1:}
Suppose 
    $  \int_0^\infty \frac{s ds}{\psi(s) \log\log s} =\infty$. \\
       If there exist infinitely many $k\ge 1$ such that $Ch(t_{k+1})\le a \Phi^{-1}(t_{k})$, then,  for infinitely many $k\ge 1$,
\begin{align*}&\P^x ( |X_{t_{k+1}}-X_{t_{k}}|\ge C h(t_{k+1})|\mathscr{F}_{t_{k}})\ge \inf_{y \in \R^d} \int_{|z-y| \ge a \Phi^{-1}(t_k)} p(t_{k}, y,z)dz\\
&=1- \sup_{y \in \R^d} \int_{|z-y| < a \Phi^{-1}(t_k)} p(t_{k}, y,z)dz>1/2.\end{align*}  
Thus, we get \eqref{sup-2-3}.

 If there is $k_0\ge3$ such that for all $k\ge k_0$, $Ch(t_{k+1})> a \Phi^{-1}(t_{k})$, then by Lemma \ref{l:int_outball} and Proposition \ref{l:LHK-psi}, 
   for all $k\ge k_0$  \begin{align*}  &\P^x ( |X_{t_{k+1}}-X_{t_{k}}|\ge C h(t_{k+1})|\mathscr{F}_{t_{k}})\ge c_1 \int_{C h(t_{k+1})}^{\infty} \frac{t_{k}}{r\psi(r)}\,d r  
\ge c_2 \frac{t_{k+1}}{\psi(h(t_{k+1}))}.
   \end{align*}
 Combining this with \eqref{e:cvi1} and the assumption that  $  \int_0^\infty \frac{s ds}{\psi(s) \log\log s} =\infty$, we also get \eqref{sup-2-3}.

       \noindent
   {\it Case 2:}
 We now assume that  $ \int_4^\infty \frac{s ds}{\psi(s) \log\log s}<\infty$. Then  for all $r>4$, 
$$\frac{1}{ \log\log r} \frac{r^2}{\Phi(r)}=\frac{1}{ \log\log r} \int_0^r \frac{s ds}{\psi(s)} \le
\frac{1}{ \log\log 4} \int_0^{\log\log 4} \frac{s ds}{\psi(s)}
 +\int_{\log\log 4}^\infty \frac{s ds}{\psi(s)\log\log s} 
=:
M <\infty, $$
and thus $\frac{r^2}{ \log\log r} <M\Phi(r).$ Thus, for any $s>4$ we have
\begin{equation}\label{e:Phi-inv}
 \Phi^{-1}(s) \le c_3 s^{1/2} (\log \log s)^{1/2}.
\end{equation}    
Also, using \eqref{e:fsecond7} to \eqref{d:Phi} we obtain
\begin{equation}\label{e:Phi-inv1}
\lim_{s \to \infty} \frac{\Phi^{-1}(s)}{s^{1/2}} = \infty.
\end{equation} 
Let $r=|x-y|$ and $\delta_1,C_5>0$ be the constants in \eqref{e:LHK-psi}. Also, let $C_6=C_6(d,\delta_1,C_5)>0$ be the constant $C_6$ in \eqref{3}. Now define $C_0=(2C_6)^{-1}$, $N=\ceil{C_0 \log k}$. Then, by \eqref{e:Phi-inv1} we have $\displaystyle \lim_{k \to \infty}\frac{\Phi^{-1}(t_k/N)}{(t_k/N)^{1/2}}=\infty$. Thus,  there exists $k_0 \in \N$ such that for any $k \ge k_0$, we have $N(k) \ge 3$ and
$$ \frac{\Phi^{-1}(t_k/N)}{(t_k/N)^{1/2}} \ge \frac{12C}{\delta_1 C_0^{1/2}}. $$
Then, there exists a constant $c_5>0$ such that for any $k \ge k_0$ and $C h(t_{k+1}) \le |x-y| \le 2C h(t_{k+1})$,
	\begin{equation}\label{2}
	p(t_k,x,y) \ge c_5  \Phi^{-1}(t_k)^{-d} k^{-1/2}.
	\end{equation}
Indeed, for $k \ge k_0$ we have
\begin{align*}
\frac{\delta_1}{3} \Phi^{-1}\left(\frac{t_k}{N}\right) = \frac{\delta_1}{3} \left(\frac{t_k}{N} \right)^{1/2} \frac{\Phi^{-1}(t_k/N)}{(t_k/N)^{1/2}} \ge \frac{4C}{C_0}2^{k/2} (\log k)^{-1/2} \ge \frac{2Ch(t_{k+1})}{N} \ge \frac{r}{N}.
\end{align*}
 Since we have \eqref{e:LHK-psi}, following the proof of Proposition \ref{l:LHK-exp1} we obtain
$$p(t_k,x,y) \ge c_6 \Phi^{-1}(t_k)^{-d} \exp(-C_6 N) \ge c_6 \Phi^{-1}(t_k)^{-d} \exp(-\frac{1}{2} \log k )= c_6 \Phi^{-1}(t_k)^{-d} k^{-1/2}. $$   
      By \eqref{2} and \eqref{e:Phi-inv} we have that for every  $k \ge k_0$, 
  \begin{align*} &  \P^x( |X_{t_{k+1}}-X_{t_{k}}|\ge C h(t_{k+1})| \sF_{t_k})\ge \inf_{y \in \R^d} \int_{Ch(t_{k+1}) \le |z-y| \le 2C h(t_{k+1})} p(t_{k}, y,z)dz \nn\\
&\ge c_6 k^{-1/2} \Phi^{-1}(t_k)^{-d} h(t_{k+1})^d \ge c_8 k^{-1/2}.
  \end{align*}
    Therefore, we conclude that 
        \begin{align*}  
& \sum_{k=k_0}^\infty \P^x( |X_{t_{k+1}}-X_{t_{k}}|\ge C h(t_{k+1})| \sF_{t_k})\ge c_5\sum_{k=k_0}^\infty  (\log k)^{-d/2} k^{-1/2}=\infty.
        \end{align*} 
        We have proved \eqref{sup-2-3}.
        \qed

\section{Examples}\label{s:Example}

In this section, we will use the notation $f(t)\simeq g(t)$ at $\infty$ (resp. $0$) if
 $\frac{f(t)}{g(t)}\to 1$ as $t\to \infty$ (resp. $t\to 0$).
We denote $\sR_0^\infty$ (resp. $\sR_0^0$) by the class of slowly varying functions at $\infty$ (resp. $0$).   For $\ell\in \sR_0^\infty$, we denote $\Pi_\ell^\infty$ (resp. $\Pi_\ell^0$) by the class of real-valued measurable function $f$  on $[c,\infty)$ (resp. $(0,c)$) 
such that for all $\lambda>0$
$$
\lim_{x\to\infty}\frac{f(\lambda x)-f(x)}{\ell(x)}=\log \lambda \quad \left(\text{resp.} \lim_{x\to 0}\frac{f(\lambda x)-f(x)}{\ell(x)}=\log \lambda\right).
$$
$\Pi_\ell^\infty$ (resp. $\Pi_\ell^0$) is called de Haan class at $\infty$ (resp. $0$) determined by $\ell$. 

For $\ell \in \sR_0^\infty$ (resp. $\sR_0^0$), we say $\ell_\#$ is de Bruijn conjugate of $\ell$ if $\ell(t)\ell_\#(t\ell(t)) \to 1, \ell_\#(t)\ell(t\ell_\#(t))\to 1$ as $t \to \infty$ (resp. $t\to 0$). Note that $|f|\in \sR_0^\infty$ if  $f\in \Pi^\infty_{\ell}$ (see \cite[Theorem 3.7.4]{BGT}).
\begin{remark}\label{l:example}{\rm
Let  $T\in(0,\infty)$.
\begin{enumerate}
\item[(1)] Let $ \ell\in \sR_0^0$ satisfying  $\int_0^1 \frac{\ell(s)}{s}ds<\infty$. Then $f(s):=\int_0^s \frac{\ell(t)}{t}dt\in \Pi^0_{\ell}$. If $\psi(s) \asymp \frac{s^2}{\ell(s)}$ for $s<T$, then for $s<T$
\begin{align*}
\Phi(s)\asymp \frac{s^2}{f(s)},\quad
\Phi^{-1}(s)\asymp s^{1/2}(1/f^{1/2})_\#(s^{1/2}), \quad
\scK_\infty^{-1}(s)\asymp s (1/f)_{\#}(s).
\end{align*}
If $g\in \Pi^0_{\ell} $ is differentiable and vanishes at $0$ and that $g'$ is monotone, then $g(s)\asymp f(s)$ for $s<T$. 
Indeed, by \cite[Theorem 3.6.8]{BGT}, $g'(s)\asymp \frac{\ell(s)}{s}$ for $s<T$. 
\item[(2)] Let $ \ell\in \sR_0^\infty$. Suppose that $\psi(s)\asymp \frac{s^2}{\ell(s)}$  for $s>T$.  If $\int_1^\infty \frac{\ell(t)}{t}dt<\infty $, then $\Phi(s)\asymp s^2$.
If $\int_1^\infty \frac{\ell(t)}{t}dt=\infty $, then for any  $f\in \Pi^\infty_{\ell} $, there exists $T_0>0$ such that  for $s>T_0$,
\begin{align*}
\Phi(s)\asymp \frac{s^2}{f(s)},\quad
\Phi^{-1}(s)\asymp s^{1/2}(1/f^{1/2})_\#(s^{1/2}), \quad
\scK_\infty^{-1}(s)\asymp s (1/f)_{\#}(s).
\end{align*}
Indeed by \cite[(1.5.8), Theorem 3.7.3]{BGT}, $
 f(s) \asymp \int_1^s \frac{\ell(t)}{t}dt
$ for $s>T_0$.
\end{enumerate}}
\end{remark}
\begin{remark}\label{r:ell}
{\rm Let $T\in(0,\infty)$ and $\ell\in \sR_0^\infty$(resp. $ \sR_0^0$).
 Suppose that $\ell$ satisfies 
\begin{align}\label{ellcheck}
\lim_{s\to \infty (\text{resp.}s\to 0)}\frac{\ell(s\ell^\gamma(s))}{\ell(s)} = 1  \text{ for some } \gamma \in \R. 
\end{align}
Then by \cite[Corollary 2.3.4]{BGT}, $(\ell^\gamma)_\# \simeq 1/\ell^\gamma$ at $\infty$ (resp. $0$).  Thus  we have the following:
\begin{enumerate}
\item[(1)] 
If $\int_0^1 \frac{\ell(s)}{s}ds<\infty$ and $f(s):=\int_0^s \frac{\ell(t)}{t}dt$ satisfies \eqref{ellcheck}  for $\gamma =1/2, 1 $, then
for $s\leq T$,
\begin{align*}
\Phi(s)\asymp \frac{s^2}{f(s)},\quad
\Phi^{-1}(s)\asymp s^{1/2}f^{1/2}(s^{1/2}), \quad
\scK^{-1}(s)\asymp s f(s).
\end{align*} 
\vspace{-0.5cm}
\item[(2)] If  $\int_1^\infty \frac{\ell(t)}{t}dt=\infty $ and $f$ is any function in $\Pi^\infty_{\ell} $ satisfying \eqref{ellcheck}  for $\gamma =1/2, 1 $, then, for $s>T$,
\begin{align*}
\Phi(s)\asymp \frac{s^2}{f(s)},\quad
\Phi^{-1}(s)\asymp s^{1/2}f^{1/2}(s^{1/2}), \quad
\scK_\infty^{-1}(s)\asymp s f(s).
\end{align*}
\end{enumerate}
}

\end{remark}

In the following corollary and examples $a_i=a_{i,L}$ or $a_i=a_{i,U} $ depending on whether we consider lower or upper bound.

\begin{corollary}\label{c:example}
Let  $T\in (0,\infty)$ and $\psi$ be a 
non-decreasing function that satisfies $L(\beta_1,C_L)$ and $U(\beta_2,C_U)$.
\begin{enumerate}
\item[(1)] Let $ \ell\in \sR_0^0$ satisfying  $\int_0^1 \frac{\ell(s)}{s}ds<\infty$ and $f(s):=\int_0^s \frac{\ell(t)}{t}dt$ satisfying \eqref{ellcheck} for $\gamma= 1/2,1$. If $\psi(s) \asymp \frac{s^2}{\ell(s)}$ for $s<1$, then for $t<T$\vspace{-0.2cm}
\begin{align*}
p(t,x,y) 
	\asymp t^{-d/2}f(t^{1/2})^{-d/2}
	\wedge \left(\frac{ t}{|x-y|^{d}\psi(|x-y|)} + t^{-d/2}f(t^{1/2})^{-d/2} \exp{\left(-\frac{a_1|x-y|^2}{tf(t/|x-y|)}\right)}\right). 
\end{align*}
\vspace{-0.15cm}
Furthermore, if $f(s^2)\asymp f(s)$
for $s<1$, then for  $t <T$,
\begin{align}\label{e:smallt}
p(t,x,y) 
	\asymp t^{-d/2}f(t)^{-d/2}\wedge \left(\frac{ t}{|x-y|^{d}\psi(|x-y|)}+ t^{-d/2}f(t)^{-d/2} \exp{\left(-a_2\frac{|x-y|^2}{tf(t)}\right)}\right). 
\end{align}
\vspace{-0.6cm}
\item[(2)]  Let $ \ell\in \sR_0^\infty$ satisfying $\int_1^\infty \frac{\ell(t)}{t}dt=\infty $. Suppose that $\psi(s)\asymp \frac{s^2}{\ell(s)}$  for $s>1$ and   $f\in \Pi^\infty_{\ell} $ satisfies \eqref{ellcheck} for $\gamma=1/2, 1$. Then for $t>T$,
\vspace{-0.1cm}
\begin{align*}
p(t,x,y) 
	\asymp  t^{-d/2}f(t^{1/2})^{-d/2}\wedge\left( \frac{ t\ell(|x-y|)}{|x-y|^{d+2}} + t^{-d/2}f(t^{1/2})^{-d/2} \exp{\Big(-\frac{a_3|x-y|^2}{tf(t/|x-y|)}\Big)}\right). 
\end{align*}
\vspace{-0.2cm}
Furthermore, if $f(s^2)\asymp f(s)$
for $s>1$, then for  $t >T$,
\begin{align}\label{e:larget}
p(t,x,y) 
	\asymp t^{-d/2}f(t)^{-d/2}\wedge \left(\frac{ t\ell(|x-y|)}{|x-y|^{d+2}} + t^{-d/2}f(t)^{-d/2} \exp{\left(-a_4\frac{|x-y|^2}{tf(t)}\right)}\right). 
\end{align}

\end{enumerate}
\end{corollary}

\pf Let $r=|x-y|$ and $\delta_1>0$ be the constant in Proposition \ref{l:LHK-psi}.

(1) Remark \ref{l:example}(1) implies that $\Phi$ satisfies $L_a(\delta, \wt C_L)$ for some $\delta>1$. Thus, by Theorem \ref{t:main}(1) and Remark \ref{r:ell}(1), we obtain the first inequality. The upper bound in the second inequality follows from Theorem \ref{t:exp-1} and $\Phi^{-1}(s)^2 \asymp sf(s^{1/2}) \asymp sf(s)$. 

For the lower bound, first choose  small $\theta>0$ such that $\frac{1}{2}+\theta(\frac{1}{2}-\frac{1}{\beta_1})=:\eps_1<1$. Note that $f(s) \asymp f(s^2)$ for $s<1$ implies $f(s^b)\asymp f(s)$ for all $b>0$ since $f$ is non-decreasing. Since the last term in the heat kernel estimate dominates other terms only in the case $\delta_1\Phi^{-1}(t) < r \le \delta_1 \frac{\Phi^{-1}(t)^{1+\theta}}{\psi^{-1}(t)^\theta}$, it suffices to show $f(t/r) \ge cf(t)$ for this case.
 Using \eqref{e:Phi} and $L(\beta_1,C_L)$ for $\psi$ we have $\Phi^{-1}(t)/\psi^{-1}(t) \le c_1 t^{\frac{1}{2}-\frac{1}{\beta_1}}$ for $t \le T$.
Thus we have $f(t/r)\geq f(c_2t^{1-\eps_1})\asymp f(t)$ for every $t\leq T$ and $\delta_1\Phi^{-1}(t)<r \le \delta_1 \frac{\Phi^{-1}(t)^{1+\theta}}{\psi^{-1}(t)^\theta}$.

(2) Remark \ref{l:example}(2) implies that $\Phi$ satisfies $L^a(\delta,\wt{C}_L)$ for some $\delta>1$. Note that $\psi(r) \asymp \frac{r^2}{\ell(r)}$ when $r>\delta_1 \Phi^{-1}(t)$ and $t >T$ since $r>\delta_1 \Phi^{-1}(t) \ge \delta_1 \Phi^{-1}(T)$. Since the second term in the heat kernel estimate dominates only in the case $r>\delta_1\Phi^{-1}(t)$, the first inequality and upper bound in the second one similarly follows from Theorem \ref{t:main}(2), Remark \ref{l:example}(2)
and Theorem \ref{e:exp-1}. Now choose small $\theta'>0$ such that $\frac{1}{\delta}+\theta'(\frac{1}{\delta}-\frac{1}{\beta_2})=:\eps_2<1$. Without loss of generality we can assume that $f$ is non-decreasing since $f(s) \asymp \int_1^s \frac{\ell(t)}{t}dt$. Now $f(s) \asymp f(s^2)$ for $s>1$ implies $f(s^b) \asymp f(s)$ for all $b>0$. Similarly, using $L^a(\delta,\wt{C}_L)$ for $\Phi$ and $U(\beta_2,C_U)$ for $\psi$ we have $\frac{\Phi^{-1}(t)}{\psi^{-1}(t)} \le c_3 t^{\frac{1}{\delta} - \frac{1}{\beta_2}}$ so $f(t/r) \ge f(c_4 t^{1-\eps_2}) \asymp f(t)$ for every $t \ge T$ and $r \le \delta_1 \frac{\Phi^{-1}(t)^{1+\theta'}}{\psi^{-1}(t)^{\theta'}}$. This finishes the proof.
\qed

We now provide some non-trivial examples where
\eqref{e:smallt}
or 
\eqref{e:larget} holds.
\begin{example}\label{ex:log} 
{\rm 
We assume that  $\psi:(0,\infty)\to(0,\infty)$ is a
non-decreasing function which 
satisfies $L(\beta_1,C_L)$ and $U(\beta_2,C_U)$.
Suppose  $\alpha >1 $ and  \begin{align*}
\psi(\lambda) \asymp 
\lambda^2(\log\frac{1}{\lambda})^\alpha, \quad 0<\lambda<1/2.
 \end{align*}
Let $\ell(s)=(\log 1/s)^{-\alpha}$ for $s \le 1/2$. Since $\alpha>1$, $f(s):=\int_0^s \frac{\ell(t)}{t}dt=(\alpha-1)^{-1}(\log 1/s))^{1-\alpha}<\infty$. Then $f$ satisfies \eqref{ellcheck} for all $\gamma\in \R$ and $f(s)\asymp f(s^2)$ for $s\leq 1/2$.
  Thus, by Corollary \ref{c:example}(1),
for $t<1/2$,
\begin{align*}
&p(t,x,y)\\
&\asymp t^{-\frac{d}{2}}\left(\log\tfrac{1}{t}\right)^{\tfrac{d(\alpha-1)}{2}}\wedge
\left(\frac{t}{|x-y|^{d}\psi(|x-y|)} + t^{-\frac{d}{2}}\left(\log\tfrac{1}{t}\right)^{\tfrac{d(\alpha-1)}{2}}\exp\left({-\tfrac{a_1|x-y|^2}{t}(\log\tfrac{1}{t})^{\alpha-1}}\right)\right).
\end{align*}

}\end{example}
\begin{example} 
{\rm 
We assume that  $\psi:(0,\infty)\to(0,\infty)$ is a
non-decreasing function which 
satisfies $L(\beta_1,C_L)$ and $U(\beta_2,C_U)$.
Suppose   $\alpha>1$ and 
\begin{align*}
\psi(\lambda) \asymp 
\lambda^2(\log\frac{1}{\lambda}) (\log\log\frac{1}{\lambda})^{\alpha} ,  \quad 0<\lambda<1/16.
 \end{align*}
 For $s\leq 1/16$, let $\ell(s)=(\log 1/s)^{-1}(\log\log 1/s)^{-\alpha}$. Then $\int_0^1 \frac{\ell(s)}{s}ds<\infty$ and  $f(s):=\int_0^s \frac{\ell(t)}{t}dt=(\alpha-1)^{-1}(\log\log 1/s)^{1-\alpha}$. Then $f$ satisfies \eqref{ellcheck} for all $\gamma\in \R$ and $f(s)\asymp f(s^2)$ for $s\leq 1/16$.
  Thus, by Corollary \ref{c:example}(1),
for $t<1/16$,
\begin{align*}
p(t,x,y)\asymp t^{-d/2}\left(\log\log\tfrac{1}{t}\right)^{\tfrac{d(\alpha-1)}{2}}&\wedge
\biggl(\frac{t}{|x-y|^d\psi(|x-y|)}\\
& + t^{-d/2}\left(\log\log\tfrac{1}{t}\right)^{\tfrac{d(\alpha-1)}{2}}\exp\left({-\tfrac{a_1|x-y|^2}{t}(\log\log\tfrac{1}{t})^{\alpha-1}}\right)\biggl).
\end{align*}

}\end{example}

\begin{example} \label{exm5.6}
\rm{
We assume that  $\psi:(0,\infty)\to(0,\infty)$ is a 
non-decreasing function which 
satisfies $L(\beta_1,C_L)$, $U(\beta_2,C_U)$, \eqref{e:intcon} and that 
$\psi(r){\1_{\{r>16\}}}\asymp r^{2}(\log r)^\beta{\1_{\{r>16\}}}$ for $\beta\in \R$. 
Let $\ell(s)=(\log r)^{-\beta}$. Then for $\beta\leq 1$, $\int^{\infty}_{1}\frac{\ell(s)}{s}=\infty$. For  $s>16$,   let 
\begin{align*}
f(s)=
\begin{cases}
\frac{(\log s)^{1-\beta}}{1-\beta} \quad & \text{if} \quad \beta<1,\\
\log\log s \quad & \text{if} \quad \beta=1.
\end{cases}
\end{align*}
Then, $f\in \Pi_{\ell}^\infty$ and $f$ satisfies  \eqref{ellcheck} for all $\gamma \in \R$ and $f(s)\asymp f(s^2)$ for $s>16$.
 On the other hand, if $\beta>1$, \eqref{e:fsecondp} holds which is equivalent to \eqref{e:fsecond2}. Thus, by Remark \ref{r:ell}(2), for $r>16$ 
\begin{align*}
\Phi(r)&\asymp 
\begin{cases}
\frac{r^2}{(\log r)^{1-\beta}}&\text{if}\quad \beta<1;\\
\frac{r^2}{\log \log r}&\text{if}\quad\beta=1;\\
r^2&\text{if}\quad\beta>1,
\end{cases} \quad\qquad\qquad
\Phi^{-1}(r)\asymp
\begin{cases}
r^{1/2}(\log r)^{(1-\beta)/2}&\text{if}\quad \beta<1;\\
r^{1/2}(\log \log r)^{1/2}&\text{if}\quad\beta=1;\\
r^{1/2}&\text{if}\quad\beta>1,
\end{cases}\\
\scK^{-1}(r)&\asymp 
\begin{cases}
r\big(\log (1+r)\big)^{1-\beta}&\text{if}\quad \beta<1;\\
r\log \log (1+r)&\text{if}\quad\beta=1;\\
r&\text{if}\quad\beta>1.
\end{cases}
\end{align*}
Thus, by the above estimates, Corollary \ref{c:example}(2) and Theorem \ref{t:main}, we have the following heat kernel estimates  for $t\ge 16$.

\vspace{3mm}

\noindent(i) If $\beta<1$:
\begin{align*}
p(t,x,y)&\asymp
\;t^{-d/2}(\log t)^{-\tfrac{d(1-\beta)}{2}}\\
&\wedge \left(\frac{t}{|x-y|^{d+2}(\log(1+|x-y|))^{\beta}}+ t^{-d/2}(\log t)^{-\tfrac{d(1-\beta)}{2}} \exp\left(-\tfrac{a_1|x-y|^2}{t(\log t)^{1-\beta}} \right)\right),
\end{align*}
(ii) If $\beta=1$:
\begin{align*}
p(t,x,y) &\asymp t^{-d/2}(\log \log t)^{d/2}\\
&\wedge \left(\frac{t}{|x-y|^{d+2}\log(1+|x-y|)}+ t^{-d/2}(\log \log t)^{d/2} \exp\left(-\tfrac{a_2|x-y|^2}{t\log\log t} \right)\right),
\end{align*}
(iii) If $\beta>1$:
\begin{align*}
p(t,x,y)
\asymp 
t^{-d/2}\wedge \left(\frac{t}{|x-y|^{d+2}(\log(1+|x-y|))^{\beta}}+ t^{-d/2} \exp\left(-a_3\tfrac{|x-y|^2}{t} \right)\right).
\end{align*}

}
\end{example}

\noindent
{\bf Acknowledgements:} 
We are grateful to Xin Chen, Zhen-Qing Chen, Takashi Kumagai and Jian Wang for giving helpful comments on earlier version of this paper.

\smallskip 
This research is  supported by the National Research Foundation of Korea(NRF) grant funded by the Korea government(MSIP) (No. 2016R1E1A1A01941893).

\end{doublespace}
\begin{singlespace}

\end{singlespace}

\vskip 0.1truein

\parindent=0em

{\bf Joohak Bae}

Department of Mathematical Sciences,

Seoul National University, Building 27, 1 Gwanak-ro, Gwanak-gu Seoul 08826, Republic of Korea

E-mail: \texttt{juhak88@snu.ac.kr}

\bigskip

{\bf Jaehoon Kang}

Research Institute of Basic Sciences,

Seoul National University, Building 501, 1 Gwanak-ro, Gwanak-gu Seoul 08826, Republic of Korea

E-mail: \texttt{jhnkang@snu.ac.kr}

\bigskip

{\bf Panki Kim}

Department of Mathematical Sciences and Research Institute of Mathematics,

Seoul National University, Building 27, 1 Gwanak-ro, Gwanak-gu Seoul 08826, Republic of Korea

E-mail: \texttt{pkim@snu.ac.kr}

\bigskip

{\bf Jaehun Lee}

Department of Mathematical Sciences,

Seoul National University, Building 27, 1 Gwanak-ro, Gwanak-gu Seoul 08826, Republic of Korea

E-mail: \texttt{hun618@snu.ac.kr}

\end{document}